\documentclass[10pt,english]{article}
\usepackage{mathptmx}
\usepackage[T1]{fontenc}
\usepackage{microtype}
\usepackage[utf8]{inputenc}
\usepackage[english]{babel}
\usepackage{amsmath}
\usepackage{amssymb}
\usepackage{bm}
\DeclareRobustCommand{\randombold}[1]{\bm{#1}}
\usepackage{amsthm}
\numberwithin{equation}{section}
\renewenvironment{cases}
 {\left\{\begin{array}{@{}l@{\qquad}l@{}}}
 {\end{array}\right.}
\usepackage{xurl}
\usepackage{setspace}
\usepackage{graphicx}
\usepackage[pass]{geometry}
\usepackage{afterpage}
\makeatletter
\newcommand{\restoretextgeometry}{%
  \restoregeometry
  \global\@colht=\textheight
  \global\@colroom=\textheight
  \global\vsize=\textheight
}
\makeatother
\AtBeginDocument{
  \setlength{\abovedisplayskip}{12pt}
  \setlength{\belowdisplayskip}{12pt}
  \setlength{\abovedisplayshortskip}{12pt}
  \setlength{\belowdisplayshortskip}{12pt}
}
\usepackage[pdfusetitle,
 bookmarks=true,bookmarksnumbered=false,bookmarksopen=false,
 breaklinks=true,pdfborder={0 0 1},backref=false,colorlinks=false]
 {hyperref}
\begin{document}
\frenchspacing
\title{The Urn of Hill, Lane and Sudderth}
\author{Simone Franchini{\normalsize\thanks{Correspondence: simone.franchini@istc.cnr.it}\thanks{CNR-ISTC, Via Gian Domenico Romagnosi 18, 00196 Rome, Italy}}}
\date{~}
\maketitle
\begin{abstract}
We review some facts, properties and applications of the urn of Hill,
Lane and Sudderth, a paradigmatic model for stochastic processes with
memory where the urn evolution is as follows: consider an urn of given
capacity, at each step a new ball, black or white, is added to the
urn with probability that is a function (urn function) of the fraction
of black balls. 
We relate the urn function to limiting compositions, growth trajectories,
rare fluctuations and free energy, and discuss its reconstruction along
nonstationary trajectories. Applications to reinforced walks, random trees and graphs, polymers,
economic and social dynamics, physical and biological growth, and machine
learning and language models illustrate coexistence, selection and sensitivity
to early events. We also consider
randomly chosen urn functions held fixed during growth, distinguishing
persistent variation between systems from fluctuations within each system.
For initial populations occupying a fixed fraction of capacity, we also
describe replica symmetry breaking through nontrivial overlap distributions
of trajectory mixtures and compare their free energies with averages over
frozen initial conditions and urn functions.

~

\noindent\textit{keywords: urn models, increasing returns, stochastic
approximation, lattice field theory}
\end{abstract}

\noindent\newpage{\large\tableofcontents{}}\newpage{}

\section{Introduction}
\begingroup
\clubpenalties 2 10000 0
\widowpenalty=10000
\displaywidowpenalty=10000
Many growing systems remember the effects of past events in their
composition, which in turn influences what is added next. Previous
adoptions can change the appeal of a technology, earlier steps can bias
a walk, and the relative sizes of competing biological populations can
affect their subsequent growth. A common question is how these local
responses determine the eventual composition: whether two types coexist,
whether one predominates, and how strongly the outcome depends on early
events. The route to that outcome also matters. Systems with the same
limiting composition can approach it at different rates and exhibit
different fluctuations. The urn of Hill, Lane and Sudderth (HLS) provides
a tractable setting in which to connect the rule for each addition to
these properties of the entire growth history.

The model starts from a prescribed initial composition. At each step,
one black or white ball is added until the urn reaches its capacity.
If $\psi$ is the current fraction of black balls, the urn function
$\pi(\psi)$ gives the probability that the next ball is black. The two
colors can therefore record any pair of outcomes whose probabilities
depend on the accumulated history through the current share. To obtain
an exact HLS representation, one identifies a count that increases by
zero or one at each step and verifies this dependence of its conditional
law. Some applications acquire this form after a change of variables
or by counting additions instead of elapsed physical time. Specifying
the count, the time variable and the initial configuration establishes
which results can be transferred to the original process.
Other applications use exact identities for count distributions,
time--reversal duality, or convergence to an HLS response; each
connection is stated with the observables and conditions it preserves.

The probabilistic foundations were developed by Hill, Lane and Sudderth
and by Blum and Brennan in connection with strong laws for dependent
sequences \cite{HLS,BB}. Arthur, Ermoliev and Kaniovskii studied
generalized urns as models of increasing returns and the emergence of
macroscopic structure from a sequence of individual choices
\cite{AEK,Arthur=000020Ermoliev=000020Kaniovski,Arthur}. These developments
belong to the broader theory of reinforced processes and stochastic
approximation, surveyed in \cite{Pemantle,MahmoudBook}. The asymptotics
of nonlinear feedback have been developed further in \cite{Gottfried_2}. Extensions to nonbinary increments are studied in
\cite{FranchiniNonbinary2025}. This review develops the tools of the
binary model and brings together their uses across different fields.
For each application, we identify the urn function and explain what
its theory says about compositions, growth trajectories and fluctuations.

\subsection{Theoretical framework}
\label{intro:theory}

\subsubsection{Equilibria and trajectories}

The shape of the urn function gives a first description of the dynamics.
Where $\pi(\psi)$ exceeds $\psi$, the expected change in the share is
positive; where it falls below $\psi$, that change is negative. The
intersections with the diagonal identify the fixed compositions. At an
isolated crossing, the share is attracted if its expected change is
positive to the left and negative to the right; the reverse signs give
repulsion. Touchpoints, endpoints and accessibility
require the qualifications discussed in Section~\ref{sec:definition}. For a continuous urn function, the share converges almost
surely to a possibly random fixed point \cite{HLS,Pemantle}. Which
points can be reached, and with what probabilities, also depends on
the response and the initial configuration. Corrective feedback can
stabilize coexistence, while amplification of an advantage can produce
competing attracting compositions. Changing the response can alter
their number or stability and can substantially slow the approach to
equilibrium \cite{Kazuaki,FB}.

Trajectories describe how those effects unfold during growth. The
continuous description in Section~\ref{sec:thermodynamic} retains both
the composition and the fraction of the final capacity already occupied.
Write $\tau$ for the continuous limit of $n/T$, where $n$ is the
current urn size, including the prescribed initial population, and
$T$ is the final capacity. For a prescribed macroscopic
initial population with share $\psi_0$ at saturation $\tau_0>0$, the
limiting share obeys
\begin{equation}
\partial_\tau\psi(\tau)=[\pi(\psi(\tau))-\psi(\tau)]/\tau,
\qquad \psi(\tau_0)=\psi_0.
\label{eq:intro-trajectory}
\end{equation}
When this problem has a unique solution, as it does for locally
Lipschitz responses, that solution describes the subsequent growth. On an interval without fixed
points, any antiderivative
\begin{equation}
\Pi(\alpha):=\int\frac{d\alpha}{\pi(\alpha)-\alpha}
\label{eq:intro-transform}
\end{equation}
gives the trajectory by inversion,
\begin{equation}
\psi(\tau)=\Pi^{-1}(\Pi(\psi_0)+\log(\tau/\tau_0)),
\label{eq:intro-solution}
\end{equation}
for as long as it remains in that interval \cite{FB}. With only
H\"older continuity, continuation through fixed points and possible
waiting times require separate treatment \cite{Franchini_URNS}. Starting instead with a fixed number of balls leaves a
different role for early fluctuations and the selection of the limiting
composition \cite{Franchini_URNS}. In that case the normalized
black-count path converges to a straight line with the limiting
share as its slope (Section~\ref{sec:trajectories}). Thus the initial share alone is insufficient: one must also specify
how the initial population scales with capacity. The trajectory plots throughout the review
show both the direction and the rate of change, including slow transients
near critical responses. They also distinguish the composition at final
capacity from the equilibrium approached under continued growth.

\subsubsection{Reconstructing the urn function}

The trajectory equation can also be used in reverse. On a nonstationary
segment where the limiting share is invertible as a function of the
growth variable, its slope determines the response that produced it.
Writing $\tau(\psi)$ for the inverse trajectory gives
\begin{equation}
\pi(\psi)=\psi+\tau(\psi)\biggl(\frac{d\tau(\psi)}{d\psi}\biggr)^{-1}.
\label{eq:intro-inverse}
\end{equation}
Section~\ref{sec:inverse} develops this reconstruction through the
transformed urn function. It recovers the original function on the
compositions visited by the trajectory. A stationary path identifies
only a fixed composition. Applying the same
idea to noisy observations requires an additional statistical estimation
procedure.

\subsubsection{Rare histories, entropy and free energy}

To describe fluctuations around these growth patterns, we use the
sample--path large--deviation principle (LDP) established for HLS urns
\cite{Franchini_URNS,FB}. It assigns a cost to a count trajectory by accumulating
the discrepancy between the local proportion of black additions along it
and the probability prescribed by the urn function. A positive cost
means that the corresponding probabilities decrease exponentially with
the final capacity. Minimizing the cost over paths satisfying a
constraint gives the exponential cost of the event and characterizes
its least costly routes \cite{Dembo_Zeitouni}. This formulation can
address a prescribed final composition or restrictions on the course
of growth. Zero--cost paths satisfy the differential equation in
Eq.~(\ref{eq:intro-trajectory}); they can also describe events whose
probabilities vanish more slowly than exponentially. Section~\ref{sec:ldp} explains the passage
from discrete histories to this variational description and the
conditions under which it applies.

For a H\"older--continuous response strictly between zero and one, and a
fixed finite initial history whose law is independent of $T$, the
variational formula is explicit. An admissible normalized black-count
path $\varphi$ is absolutely continuous, starts at zero, and has slopes
between zero and one almost everywhere; its share is
$\psi(\tau)=\varphi(\tau)/\tau$ for $\tau>0$. For an event $E$ with
limiting constraint $E^*$, let $Q(E^*)$ denote the admissible paths
satisfying that constraint. When the large--deviation lower and upper
bounds agree, the limiting logarithmic probability per addition is
\begin{equation}
\phi(E^*)=\sup_{\varphi\in Q(E^*)}\int_0^1d\tau\,
[L(\partial_\tau\varphi,\pi(\psi))-L(\partial_\tau\varphi,\partial_\tau\varphi)],
\label{eq:intro-event}
\end{equation}
where time--dependent quantities are evaluated at $\tau$ and
\begin{equation}
L(\alpha,\beta):=\alpha\log\beta+(1-\alpha)\log(1-\beta).
\label{eq:intro-local-weight}
\end{equation}
The arguments of $L$ are probabilities, with $0\log0=0$. The negative of the integrand is the cost of the additions
required along the path. The Mogulskii theorem and the Varadhan lemma
explain this variational form
\cite{Mogulskii1977,Varadhan1966,Dembo_Zeitouni}. Responses attaining zero
or one require the accessibility analysis in \cite{Franchini_URNS}.
For a prescribed macroscopic initial history, the optimization instead
runs over its continuations and accumulates cost only after that history.

The terminal entropy and the scaled cumulant--generating function
(CGF) provide complementary descriptions of the final composition.
The terminal entropy $\phi(x)$ records the limiting logarithmic
probability per addition of the exact lattice count with interior
share $x$, as defined in Eq.~(\ref{eq:terminal-entropy}). The CGF
describes the effect of reweighting histories to favor larger or
smaller final counts and has the interpretation of a free energy.
With $\psi_T$ the terminal share and $\beta$ now a real parameter,
its definition, in the conventions of Section~\ref{sec:entropy}, is
\begin{equation}
\zeta(\beta):=\lim_{T\to\infty}T^{-1}\log\mathbb E(\exp(\beta T\psi_T)).
\label{eq:intro-cgf}
\end{equation}
Under the preceding positivity and finite--history assumptions, it obeys
\begin{equation}
\pi(\partial_\beta\zeta(\beta))=
\frac{\exp(\zeta(\beta))-1}{\exp(\beta)-1},
\qquad \beta\ne0,
\label{eq:intro-cgf-equation}
\end{equation}
at its differentiability points \cite{Franchini_URNS}. Here $\zeta(0)=0$;
differentiability at zero is not assumed. This equation connects the response to terminal fluctuations on the
exponential scale. For affine urn functions, it has explicit solutions,
presented together with the Bagchi--Pal representation in
Section~\ref{app:classical}. For a concave, upper semicontinuous terminal entropy, Legendre--Fenchel
duality gives
\begin{equation}
\phi(x)=\inf_{\beta\in\mathbb R}\{\zeta(\beta)-x\beta\},
\qquad 0<x<1.
\label{eq:intro-duality}
\end{equation}
Otherwise the transform supplies the upper semicontinuous concave
envelope. Exact boundary counts also retain restrictions from the
initial history.
The same relations can also be used to construct a candidate urn function from
a prescribed entropy profile, subject to the admissibility conditions
given in the theory. This is a separate inverse problem from
reconstructing the response along a trajectory.

Some applications have a response that depends on the population size
and approaches a fixed urn function. For a fixed finite initial history
whose law is independent of capacity, uniform convergence over the whole
share range to a H\"older--continuous response bounded away from zero and
one preserves the path large--deviation costs and the limiting CGF.
Section~\ref{sec:asymptotic-urns} compares the probabilities of entire
histories. When the limit reaches zero or one, accessible histories
and relative transition probabilities must be checked separately.

Recent large--deviation studies extend the treatment of empirical memory,
the accumulated frequencies of visited states, to reinforced finite--state
chains \cite{BudhirajaWaterbury2022}, self--interacting Markov chains
with affine dependence on this empirical measure
\cite{BudhirajaWaterburyZoubouloglou2025}, and jump processes with that
affine structure \cite{BudhirajaCoghi2026}. Joint descriptions of
occupation, the fractions of time spent in each state, and transition
flows, the numbers of transitions per unit time, are also developed
under a separation of timescales: state dynamics relax faster than the
accumulated frequencies change. Suitable mixing conditions control
this loss of dependence on the starting state
\cite{CoghiBudhirajaGarrahan2026,CoghiGarrahanFlux2026}.
These results have their own assumptions; Section~\ref{sec:ldp} explains
their relation to the scalar urn.

The probability of an entire history also admits a formulation in the
language of lattice field theory (LFT). Its logarithm is a sum of
contributions from successive additions, which defines an action for
the sequence of binary variables. Summing over histories brings
together their individual weights and the entropy associated with
their multiplicity. This structure underlies the variational treatment
of fluctuations and connects the urn to action formulations used for
spin and neural systems
\cite{Franchini2024,BardellaFranchiniShort2024,BardellaFranchini2024}.
It also provides a point of contact with variational free--energy
principles in active inference \cite{Friston,Friston2,PPF}.
The connection is between variational formulations:
Section~\ref{sec:lft} shows explicitly how the action and entropy
enter the probabilities of reinforced growth histories.

\subsection{Applications}
\label{intro:applications}

\subsubsection{Classical urns, trees, graphs and reinforced walks}

The classical urns show how differently a composition can behave under
simple feedback. With fixed initial counts of both colors, P\'olya
reinforcement preserves a random limiting share, whereas the
opposite--color additions of the Friedman example restore balance.
P\'olya also describes proposed models of football passes and financial
default counts, with the corresponding initial masses
\cite{Yamamoto2025,HisakadoMori2020}. With one initial ball of each
color, its sequence law is also obtained by averaging independent
Bernoulli choices over a single uniformly drawn propensity
\cite{GelastopoulosSageRijt2025}.

After shifting both count and clock by one, descents in a permutation
grown by uniform insertion realize the Friedman urn. Precise
large--deviation estimates for this statistic are available \cite{BercuBonnefontRichou2024}.
Reversing histories with fixed endpoints relates OK Corral elimination
to Friedman growth, with the normalization in
Eq.~(\ref{eq:app-ok-corral-duality}) \cite{Fajolet2}.

The balanced Bagchi--Pal examples connect replacement matrices to
affine responses on the accessible compositions through the count of black extractions \cite{Bagchi-Pal}. The number
of black balls in the replacement urn is recovered by the corresponding
affine change of variables, provided every required replacement
remains possible and the initial composition is compatible with the
binary representation. For $\pi(\psi)=a+b\psi$ with $b<1$, the trajectory solution gives
attraction towards $a/(1-b)$ and relaxation exponent $1-b$, while the
path stays in the affine interval. For a fixed finite seed, the CGF
also has explicit solutions. For example, if $b<0$ and the response is
strictly between zero and one throughout the share interval,
\begin{equation}
\exp(-\zeta(\beta))=
\int_0^1[t^{-b}+(1-t^{-b})\exp(\beta)]^{a/b}\,dt.
\label{eq:intro-affine-cgf}
\end{equation}
Section~\ref{app:classical} gives the increasing and constant cases,
the branch conditions, and the accessible truncations needed for
subtractive replacements \cite{Franchini_URNS}.

Related counts appear in growing trees. In a recursive tree, adding
a vertex to a uniformly chosen parent changes the number of leaves by
zero or one, depending on whether that parent was already a leaf.
Leaf counts and counts of cherries (pairs of leaves with a common
parent) allow urn results to describe selected features of recursive
and Yule trees and preferential--attachment graphs
\cite{Bryc}.
Plane--oriented trees and plateau counts in Stirling permutations
illustrate the passage from a size--dependent rule to an affine HLS
response. General finite seeds in preferential attachment also give
finite--size corrections. Their limiting path LDP and CGF are those
of the corresponding affine urns, as detailed in
Section~\ref{app:classical} \cite{Bryc}.

In the elephant random walk, the count is the number of
positive steps: remembering a uniformly sampled past step makes the
next--step probability an affine function of the current share
\cite{ERW=000020shcutz=000020trimper,ERW=000020UM=000020Baur=000020Berton}.
The position and center of mass are observables of this history
\cite{Bercu}. Recent results refine Gaussian approximations in the diffusive
and critical regimes~\cite{FanHuMa2021}. In the superdiffusive regime,
they describe the density, moments and tails of the rescaled position's
limiting law \cite{GuerinLaulinRaschel2026,GuerinLaulinRaschelSimon2025}.
The full--memory Bernoulli walk remains affine, whereas applying a
nonlinear response directly to the current share gives further exact
HLS urns \cite{GutStadtmuller2026,MaulikRoySadhukhan2026}.

For a fixed sample of past steps, a response to the sampled positive
fraction is averaged over its sampling law. Independent extraction
with replacement gives a binomial average and an exact HLS
urn \cite{Fra2022,Podder}. Without replacement, the finite--size rule
is hypergeometric and converges uniformly to the same response.
A growing sample instead recovers the direct response when that response
is continuous. The transfer of path costs and the CGF uses the conditions
of Section~\ref{sec:asymptotic-urns}.

The same sampling construction gives genealogical counts in recursive
directed graphs: the count records vertices from which a specified
vertex can be reached \cite{Burghart2026}. A new vertex samples its
parents uniformly with replacement and is counted when at least one
parent is counted. In broadcasting, independent
errors affect the sampled colors before their majority is taken
\cite{BriendDevroyeLugosi2025}. With three parents sampled independently
and uniformly with replacement,
the final majority can retain information about the initial majority
even in a regime where the color share tends to one half.
Section~\ref{app:graphs} explains why this information threshold differs
from the change in stability of the balanced composition.

\subsubsection{Economic and social reinforcement}

In economic applications, reinforcement connects an early advantage
to its possible persistence \cite{Arthur}. Majority--based adoption
rules make this mechanism explicit: customers consult previous choices,
and sufficiently strong agreement with the sampled majority can create
competing attracting market shares \cite{FB}. In the example with three
independent consultations with replacement, let $p$ be the probability
of following the sampled majority. The balanced share loses stability
at $p=5/6$; above that threshold, the attracting shares are
\begin{equation}
\psi_\pm=(1\pm\sqrt{(6p-5)/(2p-1)})/2,\qquad 5/6<p\le1.
\label{eq:intro-majority}
\end{equation}
For $p<1$, both are interior: a persistent majority coexists with a
surviving minority. Section~\ref{app:economics} describes the branches
and relaxation, distinguishing a macroscopic starting composition from
basin selection by early fluctuations. Related generalized urn
models describe technological dynamics and innovation diffusion
\cite{Dosi=000020Ermoliev=000020Kaniovski,DMS}.
Information cascades also combine majority--following and independent
voters. Uniform sampling of a fixed number of previous votes without
replacement gives a response approaching the corresponding polynomial urn; Section~\ref{app:economics} states
the conditions for transferring its large deviations
\cite{HisakadoNakayamaMori2024}.
Faction sizes in a complete, structurally balanced growing network
follow the affine urn already described \cite{MohandasEtAl2025}.

Social influence can
also produce self--correction, as experiments show when a less appealing
option loses an artificial initial advantage. The measured response
curves and the marginal--majority model explain why positive social
influence alone does not determine whether lock--in occurs
\cite{VanR,Gelast}.

Another competition arises when new wealth is
assigned according to existing holdings. In the two--agent model with
a common power of wealth, stronger reinforcement can favor
concentration. Adding wages changes the feedback and can stabilize
the composition; subtracting the deterministic wage contribution
identifies the binary count to which the theory applies
\cite{Gottfried_2,Gottfried_1}. Normalizing nonlinear weights of the
two current shares gives further exact urn responses. Under their
reinforcement assumptions, these models distinguish a share tending
to one from a color receiving all subsequent allocations
\cite{RuszelThacker2024}.
The P\'olya construction for collective
decisions and the common--power model of persistent leadership provide
further examples \cite{Espinosa,Iyer=000020}.
For two--agent allocation with common positive count weights regularly
varying with positive index, the response approaches the corresponding
power urn.
The path--probability comparison in Section~\ref{app:economics},
with its fixed finite seed and support conditions, also gives the
same path LDP and CGF.

\subsubsection{Physical growth, polymers and biological systems}

Selection and coexistence also organize the physical and biological
applications. In irreversible magnetic growth at the mean--field level,
the current magnetization biases the next addition, and sufficiently
strong coupling permits two attracting polarized compositions. Counting
additions gives an urn description of their histories and rare
fluctuations \cite{KGW,Jack=000020LD}. At zero field, the transition at
coupling $J=1$ separates attraction to balance from two polarized
branches; at criticality, the trajectory relaxes logarithmically.
Whitelam and Jacobson examine branch selection and transient growth
\cite{WhitelamJacobson2021}. At zero field and $1/2<J<1$, Jack and
Harris obtain a bound on the cost of rare compositions through an
initial excursion followed by relaxation \cite{Jack=000020Harris};
Section~\ref{app:magnetic} distinguishes this bound from the full rate
function.

Random-walk range and polymers use the urn through distributions of
geometric observables. For a simple random walk on the
$d$--dimensional cubic lattice, count a black addition whenever a step
reaches a previously visited site; this count determines the range.
Averaging the next-return probability over walks with the same count
defines the microcanonical atmosphere. Its urn recursion reproduces
the count distribution exactly at every length
(Eq.~(\ref{eq:range_master})), without identifying the joint law of
successive returns with that of the urn. These distributions give the density
of states, and the entropy $\phi$ gives its exponential abundance
relative to all walks
\cite{FranchiniBalzanRANGE2018,Franchini=000020Range}.

Equation~(\ref{eq:stanley_cgf}) gives the Stanley model's free energy
in terms of the return-count CGF $\zeta$, linking the cost of geometric
constraints to the thermodynamic response. The same distribution gives
the survival probability averaged over independent traps in the
Rosenstock problem (Eq.~(\ref{eq:rosenstock_urn})); the exponential
scale alone does not resolve its finer decay. The zero-return sector
counts self-avoiding walks, and a convergent atmosphere at that boundary
determines the connective constant $\mu_d$
(Eq.~(\ref{eq:saw_connective})). Attractive contacts require an
additional geometric statistic: the partition function also involves
the conditional contact weights in
Eq.~(\ref{eq:isaw_conditional_partition}), whose self-avoiding limit
gives the interacting self-avoiding walk
\cite{Franchini=000020Range=000020Line}. Replacing the finite-length
atmospheres by an autonomous HLS response is a separate limit, subject
to the uniform convergence, regularity and boundary conditions detailed in
Section~\ref{app:polymers}.

In the model of competing
neurites, superlinear reinforcement of length can select a single
neurite that eventually receives all further growth \cite{Khanin}.
For two cell clones growing through births alone, normalizing the birth
rates gives the probability of the next newborn's type; differences in
fitness then determine the sign of the expected change in the share
\cite{BioKreger2024}. Ma et al. use the equal-rate P\'olya case
as a statistical model of composition variability in confined bacterial
communities \cite{MaRothschild2024}. The proposed controller of
cell--type proportions illustrates the use of feedback to stabilize a
prescribed mixture.

Autocatalytic growth supplies a further example: alongside the published
P\'olya treatment of handedness \cite{BioBarabas2009}, the nonlinear
rule developed here allows three attracting compositions, a balanced
one and two asymmetric ones. For example, cooperative order $q=8$, spontaneous-production parameter
$a=1/25$ and zero bias $h=0$ give three attractors separated by two
repellers (Eq.~(\ref{eq:bio-chiral})). Their distinct basins show how a growth
mechanism can sustain several persistent outcomes.

\subsubsection{Sequential allocation and admission rules}

In sequential allocation, the aim is to limit imbalance while
assignments remain random. Efron's biased coin favors the
underrepresented treatment; Wei's continuous response allows the
strength of this correction to vary with the relative imbalance
\cite{Efron1971,Wei1978}. For the smooth affine example, the theory
gives the composition trajectory and a variational formula for the
exponential cost of missing a prescribed balance tolerance. With bias
parameters proportional to the inverse sample size, the self-balancing
sampler of McNulty and Raban has the decreasing logistic response in
Eq.~(\ref{eq:app-balancing-sampler}) \cite{McNultyRaban2026}.
Urn allocation also belongs to the broader theory of adaptive clinical
designs \cite{Zhang2024Urns}; procedures that update estimates of treatment effects
retain additional statistical state \cite{YeHuMa2024,Zhang2026Designs}.

In the proposed memory model
of Section~\ref{app:memory}, feedback acts on the admission of new
items. Its two types are stored items classified as erroneous or
correct, and growth is measured by actual admissions. A filter that
accepts the two types with different probabilities can change the
probability of an erroneous addition, the equilibrium contamination
and its relaxation.
Let $p(\psi)$ be the candidate error probability and $\alpha,t$ the
acceptance probabilities of erroneous and correct candidates, respectively,
with $0<\alpha,t\le1$. The response among admitted items is
\begin{equation}
\pi(\psi)=\frac{\alpha p(\psi)}{\alpha p(\psi)+t[1-p(\psi)]}.
\label{eq:intro-admission}
\end{equation}
Thus $\alpha<t$ reduces the error probability where $0<p(\psi)<1$.
For affine candidate responses, the fixed-point condition reduces to
the quadratic equation given in Section~\ref{app:memory}.
For a fixed finite initial memory and a H\"older--continuous candidate
response strictly between zero and one, the path variational problem
describes the cost and least costly paths for exceeding a contamination
threshold.

\subsubsection{Machine learning and language models}
\label{intro:ml}

An HLS history can also be read as a binary sequence generator: the
empirical frequency of earlier symbols determines the next-symbol
probability. The negative log probability of a continuation is its
cumulative logarithmic prediction loss
(Eq.~(\ref{eq:ml-sequence-loss})), connecting the urn action to the
prediction of entire sequences. Integrating out the Beta--distributed
propensity of a Bernoulli generator gives the P\'olya predictor used in
controlled studies of prompting \cite{Genewein2025}. The same rule
describes the predictor for a new topic in the binary model of Cagnetta
and Christensen
\cite{CagnettaChristensen2026}. These identifications concern the data
generator or the ideal predictor against which a trained network is
compared. With fixed finite positive prior masses, a finite prompt
changes the limiting frequency distribution. At a prescribed final
share, however, the zero--cost count trajectory is straight, and
Eq.~(\ref{eq:ml-polya-path-cost}) penalizes departures from it.
Section~\ref{app:ml-prediction} thus distinguishes variability of the
final frequency from fluctuations along the sequence.

Repeated words within documents give another use of the same mechanism.
Grouping the vocabulary into two classes turns a Dirichlet/P\'olya
document generator into an exact binary urn, with initial masses
obtained by adding the prior parameters within each class
\cite{MadsenKauchakElkan2005}. This generator underlies the SPUD
retrieval model, whose search score combines document and collection
information \cite{CumminsPaikLv2015}. Unequal reinforcement, constant
within each class, instead gives the rational response in
Eq.~(\ref{eq:ml-document-reinforcement}) \cite{Cummins2017Burstiness}.
The zero--cost trajectories are those of clonal competition; from a
finite seed containing both classes, the share tends towards the more
strongly reinforced class. Grouped judgments of similarity
between word usages provide a further P\'olya representation in
inference of semantic change \cite{AidaEtAl2026}; here the count records
judgments rather than word occurrences (Section~\ref{app:ml-documents}).

Algorithmic decisions can in turn alter the data on which later
decisions depend. For one user and two items, the recommendation model
of Bellina et al.\ \cite{Bellina2023} gives the symmetric power response
in Eq.~(\ref{eq:ml-feedback-family}): sublinear reinforcement balances
the click shares, whereas superlinear reinforcement eventually confines
clicks to one item. In discovery-only predictive policing and the
two-arm Bernoulli bandit described in Section~\ref{app:ml-feedback},
counting recorded incidents or successful rewards gives the corresponding
response with linear weights \cite{Ensign2018,YangYanJiang2026}.
This change of clock connects these examples to the admission rule
above. Unequal linear weights produce a dominant share while additions
of the other type continue indefinitely, distinguishing dominance from
permanent fixation. From a fixed finite seed containing both types,
all these power responses have the same CGF in
Eq.~(\ref{eq:ml-feedback-cgf}). Thus this exponential description of
terminal fluctuations does not distinguish their different typical
behavior, which is resolved by trajectories and the reinforcement law.

Vocabulary growth requires a different binary observable: a black
addition records a previously unseen word. The novelty indicators of the
Pitman--Yor generators considered here, also used in authorship attribution
\cite{TaniRaffaelli2024}, have the finite-size response and HLS limit
in Eq.~(\ref{eq:ml-pitman-rule}). Their novelty share tends to zero,
although discoveries continue indefinitely. Urns that introduce new
colors upon discovery can also sustain a positive attracting novelty
share \cite{Tria2014}. In the example with reinforcement increasing
with the number of draws \cite{Bellina2025}, innovation can instead
stop with positive probability; its limiting response can have
competing attracting shares. Figure~\ref{fig:ml-innovation} compares
the responses and their growth trajectories. The finite-time novelty
probabilities distinguish a vanishing frequency of discoveries from
their eventual arrest. For the models
treated in Section~\ref{app:ml-innovation}, a uniform comparison of
history probabilities transfers the path LDP and CGF from the limiting
HLS urns to the original processes, including at the boundary where
the urn responses vanish. Sublinear
vocabulary growth and its finer fluctuations require their own analysis.

\subsubsection{Neural lattice field theory}

The HLS urn provides an exact reference to test the two--body
approximation of the neural lattice field theory (NLFT) action proposed in
\cite{BardellaFranchiniShort2024,BardellaFranchini2024,FranchiniBardella2026}.
Each lattice site represents a neuron at an observation time;
a binary variable indicates whether it is active. The action determines
the probabilities of complete recordings, and hence mean activity and
correlations. The variational description in Section~\ref{sec:lft}
combines the action with the entropy of the distribution over histories.
This provides the connection with free--energy principles and active
inference \cite{Friston,Friston2,PPF}.

For one neuron, each urn increment records activity in one time bin.
We use the affine urn $\pi(\psi)=1/2+b(\psi-1/2)$, with $0<b<1$:
$\psi$ is the fraction of past active bins and $\pi(\psi)$ the probability
of activity in the next bin. The slope $b$ measures how much this
probability changes when the past activity fraction changes.
Although the response is affine, the logarithmic history action contains
higher--order terms. Expanding the action around activity $1/2$, with
$b$ fixed, and retaining terms up to degree two gives the approximation
studied in Section~\ref{sec:lft-two-body}.

This construction determines the field $I^n$ acting on activity in
time bin $n$, and the temporal couplings $B^{nm}$ between activities
in bins $n$ and $m$, shown in Figure~\ref{fig:lft-two-body}. The resulting weights must be
normalized over complete histories; normalizing each transition
separately changes the model. An exact recursion gives the weights of finite histories.
The limiting CGF is determined by a variational principle and a
differential equation, solved numerically. Trajectories are compared
after the same prescribed initial history, occupying a fixed positive
fraction of the recording. For CGFs, the first bit is zero or one with
equal probability.

For small $b$, the approximation accurately reproduces the action per step
and the limiting free energy: Eq.~(\ref{eq:lft-remainder}) gives a uniform
error of order $b^3$, including the CGF error for every $\beta$ with the
stated first-bit law. At fixed $b$, the local action error is small
when the past activity fraction is close to $1/2$, and vanishes at
$1/2$. From initial share $1/2$ and at $\beta=0$, the two models have
exactly the same balanced limiting trajectory. Figure~\ref{fig:lft-comparison}
also shows close agreement at $b=1/2$: from initial share $0.9$ at
saturation $0.1$, the largest share difference is about $0.0052$;
the CGF error at $\beta=1$ is about $1.9\%$ relative to HLS.
Equation~(\ref{eq:lft-cgf-discrepancy}) quantifies the remaining discrepancy
for large positive $\beta$, which favors histories with many active bins.

\subsubsection{Persistent heterogeneity and quenched urns}

Persistent differences between growing systems can be represented by
choosing the urn function once and keeping it fixed throughout growth.
This is the quenched setting \cite{Bovier2006}; urns with a replacement
kernel chosen once provide a related construction
\cite[Section~6]{BandyopadhyayThacker2022}. The function's law is independent
of capacity, and the conditional theory's hypotheses hold almost surely.
Bold symbols denote random quantities: $\randombold{\pi}$ is the chosen
response, $\randombold{\zeta}$ its conditional CGF, and
$\randombold{\psi}_n$ the share under the joint choice of function and
growth history. The operators $\mathbb E_\pi,\mathbb P_\pi$ refer to
the sampled urn parameters; $\mathbb E,\mathbb P$ refer to growth
at fixed parameters. Applying the theory conditionally and then averaging
separates fluctuations within an urn from variation between urns
(Section~\ref{sec:quenched}). Even when each urn approaches an equilibrium,
those equilibria can differ. Averaging the conditional CGFs generally
differs from solving the model with the mean urn function.

In Section~\ref{app:quenched-affine}, two Bagchi--Pal matrices chosen
with equal probabilities share an attracting composition but have
different conditional fluctuations; their CGFs and quenched average
admit closed expressions. A continuous example draws the coefficients
of $\randombold{\pi}(\psi)=\randombold{a}+\randombold{b}\psi$
uniformly over the full admissible domain,
\begin{equation}
0<\randombold{a}<1,\qquad
-\randombold{a}<\randombold{b}<1-\randombold{a},
\label{eq:intro-random-domain}
\end{equation}
which has area one and includes both signs of the slope. Averaging
the CGFs gives
\begin{equation}
\mathbb E_\pi[\randombold{\zeta}(\beta)]
=\int_0^1 da\int_{-a}^{1-a}db\,\zeta(\beta),
\label{eq:intro-random-cgf}
\end{equation}
where $a,b$ are fixed integration variables and $\zeta$ is their
corresponding affine solution. The attracting share is
$\randombold{a}/(1-\randombold{b})$. From a fixed finite seed, the
variance of $\randombold{\psi}_n$ over coefficients and growth tends
to $3/4-\log2$. The averaged
CGF is differentiable at the origin but has divergent curvature;
the urn with mean coefficients instead has a Bernoulli CGF with finite
curvature (Figure~\ref{fig:quenched-affine-cgf}). Persistent variation
between equilibria and curvature of the averaged conditional CGF
therefore describe distinct effects.

Frozen biases can amplify differences between systems near criticality.
At the magnetic critical and autocatalytic tricritical parameters in
Sections~\ref{app:quenched-magnetic} and \ref{app:quenched-chiral},
magnetic growth and the illustrative autocatalytic model have equilibrium
displacements proportional to the cubic and fifth roots of a small
bias, respectively. With two attracting magnetic branches, the field
changes their positions and basins, while early growth still influences
branch selection. The cell--type controller in
Section~\ref{app:quenched-biological} distinguishes random targets from
random correction strengths: the former change limiting proportions;
the latter, towards a common target, change relaxation and finite--size
fluctuations. For an interior random target $\randombold{p}$ and finite
nonnegative strength $\randombold{\kappa}$ almost surely, the variance
of the shares over parameters and growth converges to the variance
of the sampled targets.

\subsubsection{Replica symmetry breaking in the HLS model}
\label{intro:rsb}

Section~\ref{app:rsb} compares the overlap laws and free energies
associated with random initial conditions and quenched random urn
functions. For a fixed Lipschitz urn function, initial shares
$\psi_0,\psi_0'$ at saturation $0<\tau_0<1$
determine the limiting trajectories $\psi,\psi'$. Encoding each binary
increment as the spin $2\sigma_n-1$ identifies its local mean with
$2\pi(\psi(\tau))-1$ along the first trajectory, and analogously along
the second. For two independent continuations under this same urn
function, consider the overlap over additions after initialization,
normalized by their number. The local means give
Eq.~(\ref{eq:rsb-overlap-limit}):
\begin{equation}
q_\pi(\psi_0,\psi_0')=\frac{1}{1-\tau_0}\int_{\tau_0}^{1}
[2\pi(\psi(\tau))-1][2\pi(\psi'(\tau))-1]\,d\tau.
\label{eq:intro-overlap}
\end{equation}
The initial spins contribute separately to the overlap of complete histories.

With a prescribed initial share, the conditional overlap concentrates
at one value. Drawing the share uniformly and freezing it preserves
this concentration within each system, while averaging systems gives
a distribution of diagonal overlaps. When initial shares are integrated
within a single path measure, independent histories can instead start
from different conditions. For a nonconstant response, their overlap
law is already nontrivial at a fixed urn function. This construction
applies also to nonlinear responses under the same concentration
assumptions. Affine urns provide simple examples: in the symmetric case,
the two laws reduce to a scaled square of one uniform variable
and a scaled product of two independent uniform variables on $[-1,1]$,
respectively. Both preparations have the same mean trajectory, equal
to $1/2$, despite their different overlap laws
(Figures~\ref{fig:rsb-trajectories} and~\ref{fig:rsb-overlap-densities}).
Repeating all three preparations for quenched random urn
functions then separates the overlap law within each realized urn
from its average over coefficients.

The free energy reveals a further difference between these preparations.
For affine responses with the same macroscopic initial population,
the conditional CGF is obtained through a one--dimensional quadrature
and a scalar equation. Averaging these conditional free energies over
frozen initial shares or urn coefficients can leave a differentiable
result even when the externally averaged overlap
law is broad. Averaging partition functions over the initial share
before taking the logarithm instead allows the field to select initial
conditions, producing a cusp at the origin in the total-count free
energy. The calculation also separates the contribution of the initial
population from that of future additions.

The concluding comparison in Section~\ref{app:rsb-kernel} relates this
decomposition into trajectory profiles to the representation of pure
states through local magnetizations and to the Blanket construction. It retains the component weights prescribed by the
urn dynamics and distinguishes overlap nonconcentration within a single
measure from variation between frozen environments.

\subsection{Structure of the paper}

Section~\ref{sec:theory} presents the HLS model and links convergence and
trajectories to entropy, large deviations and free energy. It also discusses
reconstruction of the urn function, limiting urn responses and quenched
randomness. Section~\ref{sec:applications} brings together applications to
classical urns and combinatorial growth, reinforced walks, economic and
social dynamics, physical growth and polymers, biological systems, and
machine learning and language models. The final part tests the two--body
NLFT approximation, examines quenched heterogeneity, and compares trajectory
overlaps and free energies under different initial preparations and urn
functions, concluding with a Kernel/Blanket interpretation of these results.

\par
\endgroup

\section{The HLS model}\label{sec:theory}

\subsection{Definition and basic properties}\label{sec:definition}

The urn model of Hill Lane and Sudderth (HLS model) was studied by
Hill, Lane and Sudderth \cite{HLS}, Blum and Brennan
in 1980 \cite{BB} and by Arthur, Ermoliev and Kaniovskii in 1983
\cite{AEK}. The model is as follows: consider an urn with total capacity
of $T$. The balls can be of two colors: say black and white. Then,
at every step $n$ a new ball is added, whose color is decided depending
on the fraction of black balls. Hereafter we call the fraction of
black balls the ``urn share'', denoted by $\psi_{n}$.
More precisely, the steps are controlled by the urn function \cite{Franchini_URNS,FB,Pemantle}
\begin{equation}
\pi:[0,1]\rightarrow[0,1]
\end{equation}
whose value gives the probability of adding a black ball. If
$\psi_{n}$ is the urn share at the step $n$, then, at the next step
$n+1$ a black ball will be added with a probability $\pi\left(\psi_{n}\right)$
and a white ball with probability $1-\pi\left(\psi_{n}\right)$
\cite{Pemantle=000020=0000202,Gouet,Kazuaki}. A detailed sample--path large deviation principle (LDP) \cite{Dembo_Zeitouni} for the
HLS urn is proved in \cite[Theorem~1]{Franchini_URNS}.
Related large--deviation and analytic results for specific urn families
are discussed in \cite{Fajolet-Analytic=000020Urns,Fajolet=0000203,Fajolet2,Bryc,Stochastic=000020urns};
further details are given in \cite{FranchiniPhD2015,FranchiniMS2011}.

Figure~\ref{fig:HLS}\label{call:fig:HLS} illustrates one step of this update rule.

\afterpage{%
\clearpage
\begin{figure}[p]
\centering{}{\Huge\includegraphics[scale=0.23]{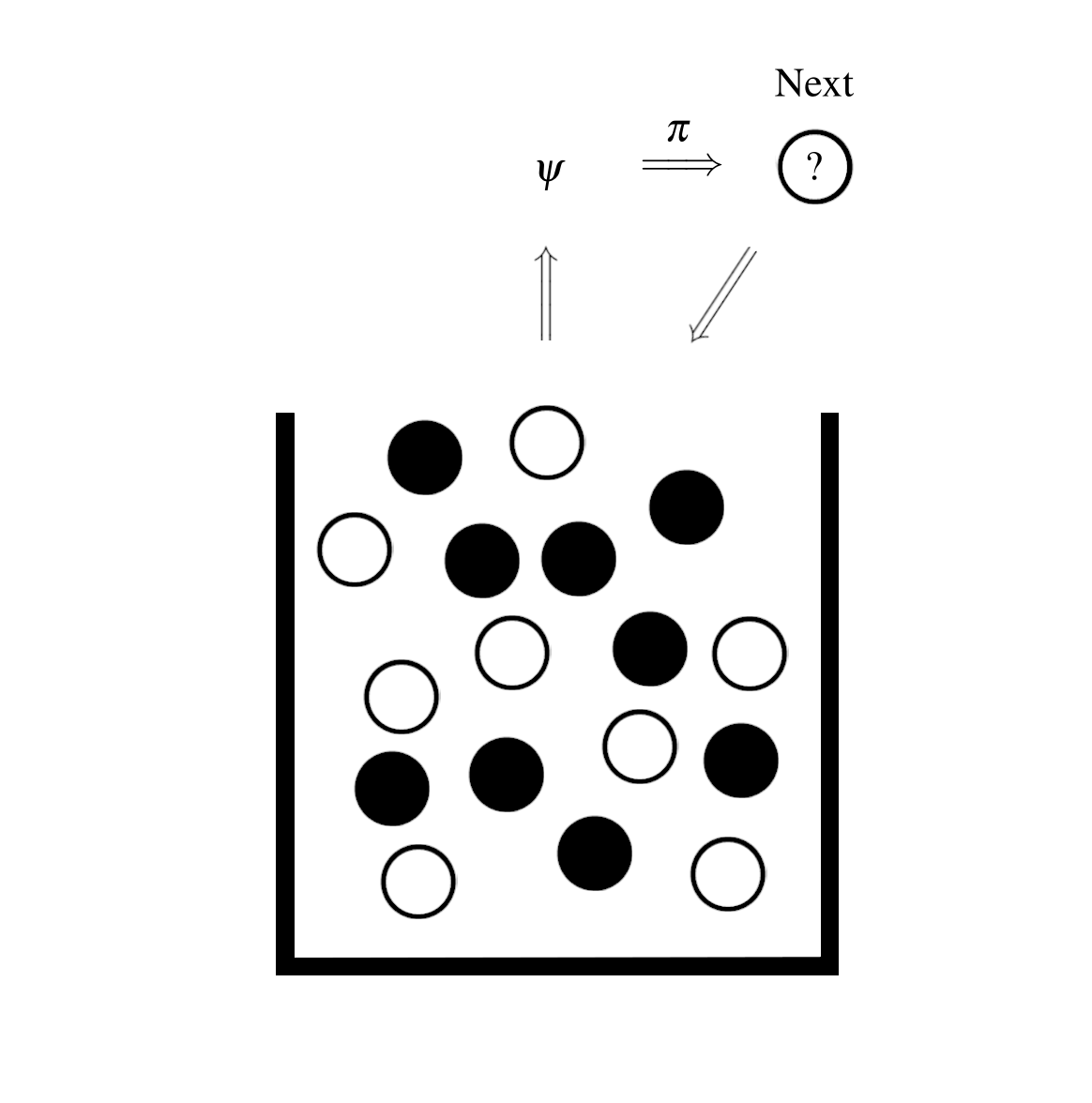}}\caption{Scheme for the step of the HLS model. The urn evolution is as follows:
for each step a new ball is added, whose color depends on the fraction
of black balls $\psi$ according to the urn function
$\pi$.}
\label{fig:HLS}
\end{figure}
\clearpage
}

\subsubsection{Basic symbols}

Let $n$ be the integer index that charts the sequence of additions;
hereafter we denote by $T$ the total capacity of the urn and by
\begin{equation}
S:=\{1,\ldots,T\} 
\end{equation}
the index set. The binary state space is
\begin{equation}
\Omega:=\{ 0,1\} 
\end{equation}
Then, the full urn history is registered in the binary vector that
records the outcome of each individual addition from first to last
\begin{equation}
\sigma:=\{ \sigma_{n}\in\Omega:\,n\in S\} \in\Omega^{S}
\end{equation}
The braces here denote an indexed family, retaining the order of its
entries. The history space is the product of the binary state spaces
labelled by the indices in $S$:
\begin{equation}
\Omega^{S}=\prod_{n\in S}\,\Omega^{\{n\}}
\end{equation}
where $\Omega^{\{n\}}$ is a copy of $\Omega$ labelled by $n$,
so that the $n$th spin belongs to that copy:
\begin{equation}
\sigma_{n}\in\Omega^{\{n\}}
\end{equation}

\subsubsection{Share of black balls}

We denote the (normalized) total number of black balls up to $n$
with 
\begin{equation}
\varphi:=\{ \varphi_{n}\in[0,1]:\,n\in S\} ,\qquad \varphi_{n}:=\frac{1}{T}\sum\limits_{1\leq m\leq n}\sigma_{m}
\end{equation}
and the average density of black balls up to $n$, or ``urn share'',
with 
\begin{equation}
\psi:=\{ \psi_{n}\in[0,1]:\,n\in S\} ,\qquad \psi_{n}:=\frac{1}{n}\sum\limits_{1\leq m\leq n}\sigma_{m}
\label{eq:urn-share}
\end{equation}
These, together with $\sigma$ itself, are alternative ways in which
we could encode our process history, the most convenient depending
on the form of the associated Lagrangian. When the action is sufficiently
regular, expansions in these fields must take account of their defining relations. A common habit
is to work in the variable $\varphi$ but this is not the only possibility,
for example in Bardella, Franchini et al. \cite{FranchiniBardella2026,BardellaFranchiniShort2024,BardellaFranchini2024}
the Taylor expansion is performed with respect to the spins $\sigma$ directly.
For HLS, the step probability depends on the accumulated share $\psi$.
Writing this share in terms of $\varphi$ introduces an explicit time
dependence in the Lagrangian, as explained in Section~\ref{sec:lft}.

\subsubsection{Urn function and the HLS process}

We will call ``urn function'' the map
\begin{equation}
\pi:[0,1]\rightarrow[0,1]
\end{equation}
or, equivalently, its indexed family of values
\begin{equation}
\pi:=\{ \pi(x)\in[0,1]:\,x\in[0,1]\} 
\end{equation}
and assume that it is Hölder--continuous. Conditional on the complete
past, the probability of the next color depends only on $\psi_n$.
For $j\in\Omega$ and $1\leq n<T$, this defines the HLS transition rule
\begin{equation}
\mathbb{P}(\sigma_{n+1}=j\mid \psi_{n}):=\pi(\psi_{n})\mathbb{I}(j=1)+(1-\pi(\psi_{n}))\mathbb{I}(j=0)\label{eq:definition}
\end{equation}
The distribution of the first spin is part of the initial condition.
An urn already containing several balls is represented by an assigned
initial history; the transition rule then applies to the subsequent
additions.

\subsubsection{Strong convergence}

We deduce the stochastic approximation equation for $\psi$
\cite[Section~2.4]{Pemantle}. Start from
\begin{equation}
\mathbb{E}(\sigma_{n+1}\mid\psi_{n})=\pi(\psi_{n})\label{eq:martin}
\end{equation}
The definitions give the identity
\begin{equation}
\sigma_{n+1}=\psi_{n}+(n+1)(\psi_{n+1}-\psi_{n})
\end{equation}
Substituting we find the stochastic approximation equation
\begin{equation}
\mathbb{E}(\psi_{n+1}-\psi_{n}\mid\psi_{n})=\frac{\pi(\psi_{n})-\psi_{n}}{n+1}
\end{equation}
For continuous $\pi$, the urn share converges almost surely to a
random point in the fixed-point set \cite{HLS,Pemantle}
\begin{equation}
C:=\{ \psi\in[0,1]:\,\pi(\psi)=\psi\} 
\end{equation}
This set need not be finite: for the Polya rule $\pi(\psi)=\psi$,
every share is fixed. At isolated crossings, a change of sign of
$\pi(\psi)-\psi$ from positive to negative is attracting, and the
reverse change is repelling. Points where the urn function touches the diagonal without crossing it
(touchpoints), and endpoints, require separate consideration; the initial history also affects which limiting shares
are accessible. See \cite{FB,Franchini_URNS,Pemantle,Pemantle=000020=0000202} for details.
Stochastic approximation also gives convergence and fluctuation results
for normalized nonlinear reinforcement and generalized elephant walks
\cite{RuszelThacker2024,MaulikRoySadhukhan2026}. The corresponding
responses are discussed in Sections~\ref{app:classical} and
\ref{app:economics}, with their model-specific assumptions.

\subsection{Thermodynamic limit}\label{sec:thermodynamic}

Here we review the main properties for the HLS model in the limit
of infinite capacity, that we identify with the analogue of the thermodynamic
limit. In particular we will describe how to find the zero--cost trajectories
of the process for any initial urn configuration, and their scaling
in that limit. We then recall the properties of the cumulant
generating function and its Legendre transform.

\subsubsection{Continuous embedding}

To study the thermodynamic limit, we first embed the urn history
in a space of continuous functions. Let $Q$ be the set of absolutely
continuous functions from $[0,1]$ to $[0,1]$ satisfying
\begin{equation}
\varphi(0)=0,\qquad \partial_\tau\varphi(\tau)\in[0,1]
\quad\text{almost everywhere}.
\end{equation}
These paths are $1$--Lipschitz \cite[Section~2.1]{Franchini_URNS}.
The embedding is provided by the map
\begin{equation}
\varphi(\sigma):=\{ \varphi(\tau\mid \sigma)\in[0,1]:\,\tau\in[0,1]\} 
\end{equation}
the actual interpolating function is 
\begin{multline}
T\varphi(\tau\mid \sigma):=T\varphi_{\,\lfloor \tau\,T\rfloor}+(\tau\,T-\lfloor \tau\,T\rfloor )\sigma_{\,\lfloor \tau\,T\rfloor+1}\\
=\lfloor \tau\,T\rfloor \psi_{\,\lfloor \tau\,T\rfloor}+(\tau\,T-\lfloor \tau\,T\rfloor )\,\sigma_{\,\lfloor \tau\,T\rfloor+1}\label{eq:embedding}
\end{multline}
where $\left\lfloor \tau\,T\right\rfloor $ is the integer part of $\tau T$.
Here $0\leq\tau<1$, $\varphi_0=0$, and the term involving $\psi$ is taken
as zero when $\lfloor\tau T\rfloor=0$; at the endpoint,
$\varphi(1\mid\sigma)=\varphi_T$. Inside each interpolation interval
the derivative is the outcome at step $\lfloor\tau T\rfloor+1$.

\subsubsection{Urn saturation and scaling limit }

We call ``urn saturation'' the ratio between the current number of
balls and the total capacity of the urn,
\begin{equation}
\tau_{n}:=n/T
\end{equation}
and use it as the time variable of the scaled process. In the scaling
limit, the capacity tends to infinity while the saturation tends to a
prescribed value,
\begin{equation}
\lim_{T\rightarrow\infty}\,\tau_{n}=:\tau\in[0,1]
\end{equation}
Here $n$ grows with $T$, with $n/T$ tending to the specified $\tau$.
For $\tau>0$, whenever
the limits exist, we denote the limiting normalized number of black balls,
or urn trajectory, and the limiting urn share by
\begin{equation}
\lim_{T\rightarrow\infty}\ \varphi_{n}=:\varphi(\tau),\qquad \lim_{T\rightarrow\infty}\ \psi_{n}=:\psi(\tau)
\end{equation}
These two functions satisfy the following identities, with derivatives
understood almost everywhere:
\begin{equation}
\varphi(\tau)=\tau\,\psi(\tau),\qquad \partial_{\tau}\varphi(\tau)=\psi(\tau)+\tau\,\partial_{\tau}\psi(\tau)
\end{equation}
By continuity of the urn function, the conditional mean increment satisfies
\begin{equation}
\lim_{T\rightarrow\infty}\,\mathbb{E}(\sigma_{n+1}\mid\psi_{n})=\pi(\psi(\tau))
\end{equation}

\subsubsection{Zero--cost trajectories}\label{sec:trajectories}

\begin{samepage}
Zero--cost trajectories are those for which the action difference
in Eq.~(\ref{eq:VARPRINZ}) vanishes. Along these trajectories
the derivative agrees with the limiting conditional mean increment,
\begin{equation}
\mathbb{E}(\sigma_{n+1}\mid\psi_{n})\rightarrow\partial_{\tau}\varphi(\tau),\qquad \pi(\psi_{n})\rightarrow\pi(\psi(\tau))
\end{equation}
\end{samepage}
so Eq. (\ref{eq:martin}) gives the differential equation
\cite[Lemma~18]{Franchini_URNS}
\begin{equation}
\partial_{\tau}\varphi(\tau)=\pi(\psi(\tau))
\label{eq:zero-cost}
\end{equation}
Writing this in terms of $\psi$ and prescribing the share $\psi_0$
at a positive saturation $\tau_0$ gives the Cauchy problem
\begin{equation}
\partial_{\tau}\psi(\tau)=\frac{\pi(\psi(\tau))-\psi(\tau)}{\tau},\qquad \psi(\tau_{0})=\psi_{0}\label{eq:Caucy}
\end{equation}
On an interval containing no fixed point, introduce the ``transformed'' urn function
\begin{equation}
\Pi(\alpha):=\int\frac{d\alpha}{\pi(\alpha)-\alpha}
\end{equation}
the solution to the Cauchy problem is readily found to be
\cite[Section~IV]{FB}
\begin{equation}
\psi(\tau)=\Pi^{-1}(\Pi(\psi_{0})+\log(\tau/\tau_{0}))
\end{equation}
as long as the trajectory remains in that interval. Stationary solutions
satisfy $\psi(\tau)=\psi_0$ with $\pi(\psi_0)=\psi_0$.
From this formula we can also find the terminal point
\begin{equation}
\psi(1)=\Pi^{-1}(\Pi(\psi_{0})-\log(\tau_{0}))
\end{equation}
provided the trajectory stays in the same interval up to saturation.
For a prescribed macroscopic initial history with $\tau_0>0$, this
Cauchy problem describes the deterministic scaling limit when its
solution is unique, as holds for locally Lipschitz urn functions.
With only Hölder continuity, a trajectory may reach a fixed point in
finite logarithmic time; continuations and possible waiting times must
then be treated separately \cite[Corollaries~6--7]{Franchini_URNS}.
Figure~\ref{fig:Urnfuntraj0T2}\label{call:fig:Urnfuntraj0T2} illustrates an urn function and its zero--cost share trajectories.

A fixed finite initial history is a different regime. The almost-sure
limit of the share gives a straight limiting count trajectory with
that slope. Other zero--cost trajectories may still describe events
whose probabilities vanish more slowly than exponentially in $T$.
Zero cost alone therefore does not imply a positive limiting probability
\cite{HLS,Franchini_URNS,FB}.

\afterpage{%
\clearpage
\begin{figure}[p]
{\Huge\includegraphics[width=\linewidth]{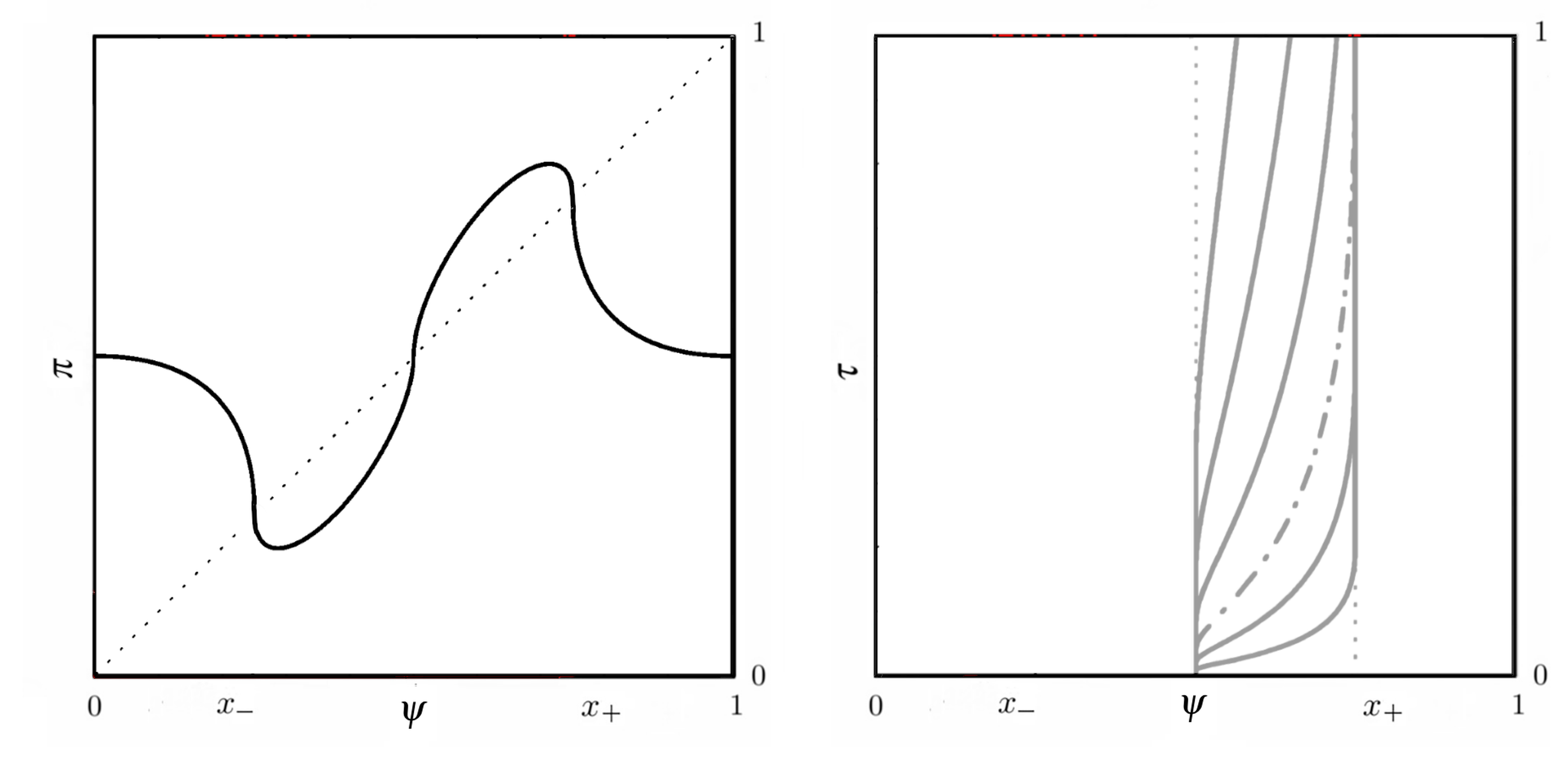}}\caption{On the left: an urn function, from Franchini 2017 \cite{Franchini_URNS}.
On the right: zero--cost share trajectories for the same function,
with saturation $\tau$ increasing upwards.
The outer fixed shares $x_-$ and $x_+$ are attracting, whereas the
central fixed share is repelling. In the right panel, vertical portions
represent a constant share before departure from the central point or
after arrival at $x_+$. They illustrate the waiting periods and arrivals
at finite saturation discussed in Section~\ref{sec:trajectories}.}
\label{fig:Urnfuntraj0T2}
\end{figure}
\clearpage
}

\subsubsection{Urn function from the trajectory}\label{sec:inverse}

\noindent From Eq. (\ref{eq:Caucy}), the fundamental equation for
the trajectory is
\begin{equation}
\Pi(\psi(\tau))-\Pi(\psi_{0})=\log(\tau/\tau_{0})
\end{equation}
On an invertible, nonstationary trajectory segment, we can merge the initial conditions into a single constant
\begin{equation}
\Pi(\psi)-\Pi_{0}^{*}=\log\tau(\psi),\qquad \Pi_{0}^{*}:=\Pi(\psi_{0})-\log\tau_{0}
\end{equation}
Thus a nonstationary trajectory represented by the pairs
\begin{equation}
\lim_{T\rightarrow\infty}(\tau_{n},\,\psi_{n})=(\tau,\,\psi(\tau)) 
\end{equation}
can equivalently be represented, after taking the logarithm of the saturation
and interchanging the coordinates, by
\begin{equation}
\lim_{T\rightarrow\infty}(\psi_{n},\,\log\tau_{n})=(\psi,\,\log\tau(\psi))=(\psi,\,\Pi(\psi)-\Pi_{0}^{*}) 
\end{equation}
So we can find the transformed urn function apart from a constant
offset by plotting the logarithm of the saturation versus the share. From this we can recover the derivative of the transformed
function, and then use
\begin{equation}
\pi(\psi)=\psi+\biggl(\frac{d\,\Pi(\psi)}{d\psi}\biggr)^{-1}
\end{equation}
to recover the urn function. The derivative can be rewritten as
\begin{equation}
\frac{d\,\Pi(\psi)}{d\psi}=\frac{d\,\log\tau(\psi)}{d\psi}=\frac{1}{\tau(\psi)}\biggl(\frac{d\tau(\psi)}{d\psi}\biggr)
\end{equation}
Substituting this expression in the previous formula
\begin{equation}
\pi(\psi)=\psi+\tau(\psi)\biggl(\frac{d\tau(\psi)}{d\psi}\biggr)^{-1}
\end{equation}
we recover the urn function on the shares visited by the trajectory,
in terms of the inverse trajectory and its derivative. Estimation from noisy data
requires an additional statistical procedure.

\subsubsection{Entropy density and cumulant--generating function}\label{sec:entropy}

The scaled cumulant--generating function (CGF) is defined as
\begin{equation}
\zeta(\beta):=\lim_{T\rightarrow\infty}\,\frac{1}{T}\log\mathbb{E}(\exp(\beta T\psi_{T}))
\label{eq:cgf}
\end{equation}
and is related to the free energy, up to sign and normalization conventions.
For an initial history of fixed length whose law does not depend on $T$,
it satisfies the following equation at its differentiability points
\cite[Eq.~(2.31) and Section~3.3.1]{Franchini_URNS}:
\begin{equation}
\pi(\partial_{\beta}\zeta(\beta))=\frac{\exp(\zeta(\beta))-1}{\exp(\beta)-1},
\qquad \beta\ne0.
\label{eq:cgf-equation}
\end{equation}
Here $\zeta(0)=0$, but differentiability at zero is not assumed.
The affine case is solved explicitly in
Section~\ref{app:classical}; Section~\ref{app:quenched-affine}
then computes quenched averages for random coefficients.
For broader classes of self-interacting dynamics, spectral methods couple
tilted operators, evolution operators modified by the exponential
weights used in generating functions, to the evolution of the empirical
memory \cite{CoghiGarrahanSpectral2026}. They provide a framework for computing
fluctuations, while the equation above is the specific binary urn relation.
The closely related entropy density
of a terminal share $x\in[0,1]$ is defined by
\begin{equation}
\phi(x):=\lim_{T\rightarrow\infty}\,\frac{1}{T}\log\,\mathbb{P}(\psi_{T}=\lfloor xT\rfloor /T)
\label{eq:terminal-entropy}
\end{equation}
whenever the limit exists. For $0<\pi<1$ and $0<x<1$, this local
limit agrees with the terminal large-deviation entropy
\cite[Theorem~4]{Franchini_URNS}. At the endpoints, the initial history matters
even in this positive case: one initial black ball makes a final
black count of zero impossible, although arbitrarily small terminal
shares remain accessible. The entropy of an exact boundary count
must therefore be distinguished from the limit through interior shares.

When the terminal entropy is concave and upper semicontinuous,
Legendre--Fenchel duality recovers it as
\cite{Dembo_Zeitouni,Franchini_URNS}
\begin{equation}
\phi(x)=\inf_{\beta\in\mathbb{R}}\{ \zeta(\beta)-x\beta\} 
\label{eq:entropy-duality}
\end{equation}
For the point-probability definition above, this statement applies
at interior shares; at an endpoint it requires agreement with the
corresponding limiting entropy. Otherwise the transform supplies
the upper semicontinuous concave envelope of the terminal entropy.
Urn functions attaining zero or one also require the accessibility
conditions discussed in Section~\ref{sec:ldp}.
At differentiable dual points $\beta=-\partial_x\phi(x)$, so the substitutions
\begin{equation}
\partial_{\beta}\zeta(\beta)=x,\qquad \zeta(\beta)=\phi(x)-x\partial_{x}\phi(x)
\end{equation}
in the equation for the generating function give \cite{Franchini_URNS,FranchiniMS2011}
\begin{equation}
\pi(x)=\frac{\exp(\phi(x)-x\partial_{x}\phi(x))-1}{\exp(-\partial_{x}\phi(x))-1}
\end{equation}
This relation holds at those dual points where the denominator is
nonzero. It gives a differential equation for the entropy, or a
candidate urn function for a prescribed entropy profile. In the latter
direction one must still check that the candidate takes values in
$[0,1]$ and satisfies the regularity and boundary conditions of the
inverse result in \cite[Corollary~10]{Franchini_URNS}. A flat zero-entropy region
gives a zero numerator and denominator and does not determine $\pi$
through this formula.

\subsection{Large Deviations and Lattice Field Theory}\label{sec:ldp}\label{sec:asymptotic-urns}

We follow the sample--path LDP proof in
\cite[Theorem~1 and Section~3.1]{Franchini_URNS}, based on the
Mogulskii theorem and the Varadhan lemma
\cite{Mogulskii1977,Varadhan1966,Dembo_Zeitouni}.
The direct argument below is for $0<\pi<1$, starting from a single
spin with a fixed law. An assigned initial history of fixed length
changes only finitely many factors and gives the same LDP
\cite[Section~3.1.5]{Franchini_URNS}. The cases where $\pi$ attains
zero or one require the accessibility treatment in
\cite[Corollary~2 and Sections~3.1.3--3.1.5]{Franchini_URNS}.

Recent large-deviation results place this construction within a wider
theory of empirical memory. Reinforced finite-state chains and
self-interacting Markov chains have dynamical variational costs
\cite{BudhirajaWaterbury2022,BudhirajaWaterburyZoubouloglou2025};
their binary affine cases include HLS responses. General stochastic
approximation schemes are treated in \cite{HultEtAl2025}.
For continuous-time jump processes, occupation measures and transition
fluxes are studied jointly in \cite{BudhirajaCoghi2026,CoghiBudhirajaGarrahan2026,CoghiGarrahanFlux2026}.
The rate retains the evolution of the memory along a rare history.
These results have their own hypotheses: the full probabilistic LDPs in
\cite{BudhirajaWaterburyZoubouloglou2025,BudhirajaCoghi2026}
require affine dependence on the empirical measure, while
\cite{CoghiBudhirajaGarrahan2026,CoghiGarrahanFlux2026}
use a separation of timescales with corresponding mixing conditions.

\subsubsection{Events}

An event is a collection of urn histories
\begin{equation}
E\subseteq\Omega^{S}
\end{equation}
For example we can consider the following:
\begin{equation}
E=\{ \psi_{T}\in[x,y]\} 
\end{equation}
This is the set of urn histories whose terminal share lies between $x$
and $y$, with $x<y$. Given any event
we can compute its image on the $1-$Lipschitz embedding space
\begin{equation}
Q(E):=\{\,\varphi(\sigma)\in Q:\,\sigma\in E\}
\end{equation}
The image of any event lies in the $1-$Lipschitz function space
\begin{equation}
Q(E)\subseteq Q(\Omega^{S})\subseteq Q,\qquad \forall E\subseteq\Omega^{S}
\end{equation}
We denote the scaling limit of an event by a star
\begin{equation}
E^{*}:=\lim_{T\rightarrow\infty}E
\end{equation}
which is assumed here to have non--empty interior. For example, 
\begin{equation}
E:=\{ \psi_{T}\in[x,y]\} ,\qquad E^{*}:=\{ \psi(1)\in[x,y]\} 
\end{equation}
The case $x=y$ concerns point probabilities for the terminal share and requires a
local estimate in addition to the contraction principle
\cite[Theorem~4 and Section~3.2.1]{Franchini_URNS}. Here we focus
on intervals with $x<y$.
Conditioning on a history constraint can itself generate dependence on
the past, as illustrated by the Doob constructions in
\cite{CoghiGarrahanDoob2025}; their transition probabilities may also
depend on the time remaining to the prescribed endpoint.

\subsubsection{Scaling limit of the entropy density}

Suppose we are interested in computing the entropy density of a given
event
\begin{equation}
\phi(E^{*}):=\lim_{T\rightarrow\infty}\,\frac{1}{T}\log\,\mathbb{P}(\sigma\in E)\label{eq:entropydensity}
\end{equation}
where $E^{*}$ is a suitable scaling limit for the event $E$. When the
large--deviation lower and upper bounds agree, this limit is given by
the variational problem \cite[Theorem~1 and Section~3.1.2]{Franchini_URNS}
\begin{equation}
\phi(E^{*})=\sup_{\varphi\in Q(E^{*})}\{ \Phi(\varphi)-\Phi_{0}^{*}(\varphi)\} \label{eq:VARPRINZ}
\end{equation}
where the first functional is the scaled action associated to the
specific urn process, that actually contains the dependence on the
urn function
\begin{equation}
\Phi(\varphi):=\int_{0}^{1}d\tau\,L(\partial_{\tau}\varphi(\tau),\pi(\psi(\tau)))\label{eq:scaledaction}
\end{equation}
and the second is the mixing term inherited from the independent
Bernoulli reference process:
\begin{equation}
\Phi_{0}^{*}(\varphi):=\int_{0}^{1}d\tau\,L(\partial_{\tau}\varphi(\tau),\partial_{\tau}\varphi(\tau))\label{eq:mogulskii}
\end{equation}
The term $-\Phi_0^*(\varphi)$ is the binary mixing entropy: it accounts,
on the exponential scale, for how many sequences of black and white
additions have a count path close to $\varphi$. As shown below,
$\Phi_0^*$ differs from the nonnegative Mogulskii rate functional by
the constant $\log2$. The local function in both expressions is
\begin{equation}
L(\alpha,\beta):=\alpha\log\beta+(1-\alpha)\log(1-\beta)\label{eq:scaler}
\end{equation}
Here both arguments are in $[0,1]$; the second argument of $L$ is a
probability, rather than the tilting parameter of $\zeta$.
We use $0\log0=0$; a positive coefficient multiplying $\log0$ gives
$-\infty$. The difference $\Phi-\Phi_0^*$ is nonpositive: its negative
is the integral of the relative entropy between two Bernoulli laws,
with parameters given by the path derivative and the urn probability.
It vanishes precisely when Eq.~(\ref{eq:zero-cost}) holds almost everywhere.
A way to prove this result is to perform a change of measure, followed
by a continuous embedding, then compute the candidate action by the
scaling rules introduced above and, finally, apply the Varadhan lemma and the
Mogulskii theorem.

\subsubsection{Path integral formulation}\label{sec:lft}

As preparation we will recast the HLS process in the language of a
Lattice Field Theory (LFT) \cite{KERNEL=000020THEO,Franchini2024,BardellaFranchiniShort2024,BardellaFranchini2024}.
The average of any function of the process history can be written as
\cite{Franchini_URNS,FranchiniMS2011}
\begin{equation}
\mathbb{E}(\mathcal{O}(\sigma)):=\sum\limits_{\sigma\in\Omega^{S}}\mathcal{O}(\sigma)\,\frac{\exp(\mathcal{A}(\sigma))}{\displaystyle\sum\limits_{\sigma'\in\Omega^{S}}\exp(\mathcal{A}(\sigma'))}
\label{eq:path-gibbs}
\end{equation}
For the action defined below, the denominator is one because the exponential
of the action is already the normalized path probability. Retaining the denominator
makes the Gibbs form explicit.
This representation also provides
the connection with variational free--energy principles and active
inference \cite{Friston,Friston2,PPF}: the variational description
involves both the action and the entropy of the distribution over
histories. The lattice field formulation makes these ingredients
explicit in a reinforced process, alongside the corresponding
descriptions of neural and spin networks cited above. This is a
connection between variational formulations, not an identification
of the urn with a particular active-inference algorithm.
The action function is defined from its Lagrangian
\begin{equation}
\mathcal{A}(\sigma):=\log\mathbb{P}(\sigma_1)+\sum\limits_{1\leq n<T}\mathcal{L}(\sigma_{n+1},\psi_{n})
\label{eq:path-action}
\end{equation}
where $\mathbb P(\sigma_1)$ denotes the probability assigned to the first
spin of the history. From the transition rule in Eq. (\ref{eq:definition})
we can deduce that the probability of the event $E$ is provided by
the formula
\begin{equation}
\mathbb{P}(\sigma\in E)=\sum\limits_{\sigma\in E}\,\mathbb{P}(\sigma_1)\prod_{1\leq n<T}\,\pi(\psi_{n})^{\sigma_{n+1}}(1-\pi(\psi_{n}))^{1-\sigma_{n+1}}
\end{equation}
Therefore, we could actually interpret the logarithms of the weights
of the individual steps as the analogue Lagrangian for the process,
and the logarithm of the full weight as the analogue action. For the
HLS model the Lagrangian is
\begin{equation}
\mathcal{L}(\sigma_{n+1},\psi_{n})=\sigma_{n+1}\log\pi(\psi_{n})+(1-\sigma_{n+1})\log(1-\pi(\psi_{n}))
\label{eq:path-lagrangian}
\end{equation}
When expressed through the count trajectory $\varphi$, the limiting
Lagrangian depends explicitly on $\tau$ through the share
$\psi(\tau)=\varphi(\tau)/\tau$.
Section~\ref{sec:lft-two-body} uses this action to test the two--body
approximation of neural LFT against an affine HLS process.

\subsubsection{Change of measure}

The first step is to perform a change of measure
\cite[Section~3.1.2]{Franchini_URNS}
\begin{equation}
\mathbb{P}(\sigma\in E)=\sum\limits_{\sigma\in E}\exp(\mathcal{A}(\sigma))
=\sum\limits_{\sigma\in\Omega^{S}}\exp(\mathcal{A}(\sigma))\mathbb{I}(\sigma\in E)
\end{equation}
\begin{equation}
\mathbb{P}(\sigma\in E)=2^{T}\,\mathbb{E}_{0}(\exp(\mathcal{A}(\sigma))\mathbb{I}(\sigma\in E))\label{eq:stt}
\end{equation}
where the last average is defined as
\begin{equation}
\mathbb{E}_{0}(\mathcal{O}(\sigma)):=\frac{1}{2^{T}}\sum\limits_{\sigma\in\Omega^{S}}\mathcal{O}(\sigma)
\end{equation}
and is taken according to i.i.d. steps distributed as follows:
\begin{equation}
\mathbb{P}_{0}(\sigma_{n}=1):=1/2,\qquad \forall n\in S
\end{equation}

\subsubsection{Scaling of the action}

The candidate scaling limit is found by the coarse--grained substitutions
\begin{equation}
\sigma_{n+1}\rightarrow\partial_{\tau}\varphi(\tau),\qquad \psi_{n}\rightarrow\psi(\tau)
\end{equation}
in the previous definition of the Lagrangian. The first substitution
means averaging increments over short macroscopic intervals; individual
binary increments need not converge to the path derivative. The sum
is replaced by an integral,
\begin{equation}
\frac{1}{T}\sum\limits_{n\in S}\rightarrow\int_{0}^{1}d\tau
\end{equation}
The action divided by the terminal size then converges to
Eq.~(\ref{eq:scaledaction}) along histories allowed by the initial
law whose interpolations converge uniformly to $\varphi$
\cite[Lemmas~13--14]{Franchini_URNS}:
\begin{equation}
\mathcal{A}(\sigma)/T\rightarrow\Phi(\varphi)
\end{equation}
where the discrete Lagrangian is expressed through the same function $L$ as
\begin{equation}
\mathcal{L}(\sigma_{n+1},\psi_{n})=L(\sigma_{n+1},\pi(\psi_n))
\end{equation}
with $L$ as in Eq.~(\ref{eq:scaler}). The next step specifies the
uniform control needed to use this passage inside the path average.

\subsubsection{Varadhan lemma}

We apply the Varadhan lemma \cite{Varadhan1966}; see
\cite[Theorem~4.3.1]{Dembo_Zeitouni}. The quantity of interest is the limit in Eq.~(\ref{eq:entropydensity}).
For the positive Hölder urn functions considered here, the scaled
action has the following uniform approximation over histories allowed
by the initial law \cite[Lemma~14]{Franchini_URNS}:
\begin{equation}
\lim_{T\rightarrow\infty}\sup_{\sigma\in\Omega^S:\,\mathbb P(\sigma)>0}
|\mathcal{A}(\sigma)/T-\Phi(\varphi(\sigma))|=0
\end{equation}
The limiting action is continuous in the uniform norm
\cite[Lemma~13]{Franchini_URNS}: for any two paths $\varphi,\varphi'\in Q$,
\begin{equation}
|\,\Phi(\varphi)-\Phi(\varphi')|\leq\delta(\epsilon),\qquad \Vert \varphi-\varphi'\Vert \leq\epsilon
\end{equation}
where $\epsilon$ is a small positive number, $\delta$ is a continuous
function vanishing at zero and
\begin{equation}
\Vert \varphi-\varphi'\Vert :=\sup_{\tau\in[0,1]}|\varphi(\tau)-\varphi'(\tau)|
\end{equation}
Conditioning the reference law on an allowed initial history of fixed
length leaves its path LDP unchanged. Applying the Varadhan lemma
\cite[Section~4.3]{Dembo_Zeitouni}, for events whose upper and lower
bounds agree, therefore gives
\begin{equation}
\phi(E^{*})-\log2=\sup_{\varphi\in Q(E^{*})}\{ \Phi(\varphi)-\Phi_{0}(\varphi)\} \label{eq:VARADA}
\end{equation}
where $\Phi_0$ is the nonnegative rate functional of
the independent Bernoulli reference process and is found
by Mogulskii's theorem \cite{Mogulskii1977}.

\subsubsection{Mogulskii theorem}

For the independent Bernoulli increments used here, the Mogulskii
theorem \cite[Theorem~5.1.2]{Dembo_Zeitouni} gives the rate functional
\begin{equation}
\Phi_{0}(\varphi):=\int_{0}^{1}d\tau\,L_{0}(\partial_{\tau}\varphi(\tau))\label{ratefunction-1-1}
\end{equation}
where $L_{0}$ is the Mogulskii Lagrangian, defined as the Legendre
transform 
\begin{equation}
L_{0}(\alpha):=\sup_{\beta\in\mathbb{R}}\{ \alpha\beta-\zeta_{0}(\beta)\} 
\end{equation}
of the cumulant generating function of the increment
\begin{equation}
\zeta_{0}(\beta):=\log(\mathbb{E}_{0}(\exp(\beta\sigma_{1})))
\end{equation}
For $0\leq\alpha\leq1$, the binary Lagrangian is
\begin{equation}
L_{0}(\alpha)=\log2+\alpha\log\alpha+(1-\alpha)\log(1-\alpha)
\end{equation}
It is infinite outside this interval, since the reference increments
are zero or one.
The $\log2$ constant comes from the normalization of the Bernoulli law
with probability $1/2$. Cancelling it against the same term in Eq. (\ref{eq:VARADA})
gives \cite{FB,Franchini_URNS}
\begin{equation}
L_{0}^{*}(\alpha):=L_{0}(\alpha)-\log2
=L(\alpha,\alpha)=\alpha\log\alpha+(1-\alpha)\log(1-\alpha)
\end{equation}
which gives the integrand in Eq.~(\ref{eq:mogulskii}). The star here
denotes this removal of the normalization constant, not another
Legendre transform. The nonnegative rate functional is $\Phi_0$;
the urn rate is $\Phi_0^*-\Phi$.

Some applications have a time-dependent response $\pi_n(\psi_n)$
approaching a fixed urn function. Let $\pi$ satisfy the regularity
conditions above. Starting from a fixed finite initial history, suppose
that both responses are bounded away from zero and one uniformly for
all sufficiently large $n$, and that
\begin{equation}
\lim_{n\to\infty}\sup_{0\leq j\leq n}
|\pi_n(j/n)-\pi(j/n)|=0.
\label{eq:asymptotic-kernels}
\end{equation}
Condition on an admissible prefix of fixed length $m$ extending beyond
the exceptional initial steps. The logarithms are uniformly Lipschitz
on the common probability range. Summing the one-step
differences and using their vanishing Ces\`aro average gives
\begin{equation}
\sup_{\sigma}\Bigl|\sum\limits_{m\leq n<T}
[L(\sigma_{n+1},\pi_n(\psi_n))-L(\sigma_{n+1},\pi(\psi_n))]\Bigr|=o(T).
\label{eq:asymptotic-comparison}
\end{equation}
The supremum is over the common continuations of that prefix.
Their probabilities differ by factors between $\exp(-o(T))$ and
$\exp(o(T))$, so the path LDP and the limiting CGF are the same.
Each conditional HLS law has the same rate, independent of its finite
prefix; the finite mixture over admissible prefixes therefore preserves
these limits.
When a response reaches zero or one, uniform convergence alone does
not give this comparison: the supports and relative transition
probabilities must be checked, as in the applications below.

\begingroup
\clubpenalty=10000
\widowpenalty=10000
\displaywidowpenalty=10000
\subsection{Quenched random urn functions}\label{sec:quenched}

A random urn function can describe persistent differences between
otherwise identical growing systems. We choose the whole function once
and keep it fixed throughout the evolution \cite[Chapter~6]{Bovier2006}.
Urn schemes with a replacement kernel sampled once and a conditional
quenched law are formulated in \cite[Section~6]{BandyopadhyayThacker2022};
their two--color reduction is discussed in Section~\ref{app:quenched-affine}.
Bold symbols denote random coefficients and quantities determined by the
chosen urn: $\randombold{\pi}$ is the chosen urn function,
while $\randombold{\Phi}$, $\randombold{\phi}$ and
$\randombold{\zeta}$ are its action, entropy and scaled
cumulant--generating function. The random increment and share under the
joint choice of urn function and growth history are denoted by
$\randombold{\sigma}_n$ and $\randombold{\psi}_n$.
An ordinary $\psi$ denotes a fixed share argument or a deterministic
trajectory at fixed parameters; prescribed histories and paths retain
$\sigma,\varphi$. The operators $\mathbb E_\pi,\mathbb P_\pi$ refer to
the one-time choice of the function or of its defining parameters,
including a sampled replacement matrix. The operators $\mathbb E,\mathbb P$
describe growth conditional on that choice. For a random initial share,
$\mathbb E_{\psi_0},\mathbb P_{\psi_0}$ refer to its sampling law, as in
Section~\ref{app:rsb}. The reference product law retains
$\mathbb E_0,\mathbb P_0$. In applications, unsubscripted operators refer
to the model law specified locally. Conditional on the chosen function,
the transition rule is
\begin{equation}
\mathbb P(\randombold{\sigma}_{n+1}=1\mid\randombold{\sigma}_1,\ldots,\randombold{\sigma}_n,\randombold{\pi})
=\randombold{\pi}(\randombold{\psi}_n).
\label{eq:quenched-rule}
\end{equation}
The law used to choose the function is independent of the final
capacity $T$. For the applications below, almost every realization
is H\"older--continuous and strictly between zero and one. We use a
fixed finite initial history unless a macroscopic initial history is specified.
The preceding results then apply conditionally,
with exactly the same initial--condition and event assumptions
\cite{HLS,Franchini_URNS}.

For a fixed realization $\pi$ of the urn function, the deterministic
zero--cost share trajectories, starting at a positive saturation $\tau_0$, satisfy
\begin{equation}
\partial_\tau\psi(\tau)
=[\pi(\psi(\tau))-\psi(\tau)]/\tau,
\qquad \psi(\tau_0)=\psi_0.
\label{eq:quenched-trajectory}
\end{equation}
Thus the equilibria and their stability are determined separately
for each realization. The action $\randombold{\Phi}$ is obtained
by inserting $\randombold{\pi}$ in Eq.~(\ref{eq:scaledaction}), whereas
the mixing term $\Phi_0^*$ is unchanged. Whenever the conditional
large--deviation bounds for the specified event agree,
\begin{equation}
\randombold{\phi}(E^*)=
\sup_{\varphi\in Q(E^*)}
\{\randombold{\Phi}(\varphi)-\Phi_0^*(\varphi)\}.
\label{eq:quenched-event}
\end{equation}
For an assigned macroscopic initial history ending at $\tau_0>0$, the space
$Q$ is restricted to its continuations, with
$\varphi(\tau_0)=\tau_0\psi_0$, and both action integrals run from
$\tau_0$ to $1$. The agreement of the large--deviation bounds is required within that
space of continuations. The function $\randombold{\phi}$ describes
probabilities conditional on the chosen urn function and initial history.

The conditional scaled cumulant--generating function is
\begin{equation}
\randombold{\zeta}(\beta)=\lim_{T\to\infty}T^{-1}
\log\mathbb E(\exp(\beta T\randombold{\psi}_T)\mid\randombold{\pi}).
\label{eq:quenched-cgf-definition}
\end{equation}
For a fixed finite initial history it obeys Eq.~(\ref{eq:cgf-equation})
at its differentiability points \cite[Section~3.3.1]{Franchini_URNS}:
\begin{equation}
\randombold{\pi}(\partial_\beta\randombold{\zeta}(\beta))
=(\exp(\randombold{\zeta}(\beta))-1)/(\exp(\beta)-1),
\qquad \beta\ne0.
\label{eq:quenched-cgf}
\end{equation}
Here $\randombold{\zeta}(0)=0$; differentiability at zero is not
required. Where the urn function is invertible, the same identity gives
\begin{equation}
\partial_\beta\randombold{\zeta}(\beta)=
\randombold{\pi}^{-1}((\exp(\randombold{\zeta}(\beta))-1)/(\exp(\beta)-1)),
\qquad \beta\ne0.
\label{eq:quenched-cgf-inverse}
\end{equation}
The inverse is taken on the image of the relevant branch. This
algebraic step applies to strictly increasing or strictly decreasing
functions; existence and selection of a solution retain the conditions
of the deterministic problem. Constant responses are treated directly
with Eq.~(\ref{eq:quenched-cgf}).

The quenched average of the free energy \cite[Chapter~6]{Bovier2006},
with the convention of Section~\ref{sec:entropy}, is
\begin{equation}
\mathbb E_\pi[\randombold{\zeta}(\beta)]
=\lim_{T\to\infty}\mathbb E_\pi[T^{-1}
\log\mathbb E(\exp(\beta T\randombold{\psi}_T)\mid\randombold{\pi})].
\label{eq:quenched-free-energy}
\end{equation}
The outer expectation $\mathbb E_\pi$ averages over the sampled urn
parameters; the inner expectation averages only over growth. The exchange
of limit and expectation follows by dominated convergence
\cite[Chapter~1]{Kallenberg2021}, using the conditional limit and the bound
\begin{equation}
\min(0,\beta)\leq T^{-1}
\log\mathbb E(\exp(\beta T\randombold{\psi}_T)\mid\randombold{\pi})
\leq\max(0,\beta).
\label{eq:quenched-cgf-bound}
\end{equation}
Thus one first determines the free energy for each fixed function
and then averages it. Substituting the mean urn function into
Eq.~(\ref{eq:cgf-equation}) does not in general perform this calculation.

For observables, the two sources of variability are recorded by the
law of total variance \cite[Chapter~8]{Kallenberg2021}.
With $\operatorname{Var}$ denoting variance under the indicated law,
\begin{equation}
\operatorname{Var}(\randombold{\psi}_n)
=\mathbb E_\pi[\operatorname{Var}(\randombold{\psi}_n\mid\randombold{\pi})]
+\operatorname{Var}(\mathbb E[\randombold{\psi}_n\mid\randombold{\pi}]).
\label{eq:quenched-total-variance}
\end{equation}
The variance on the left is over both the sampled parameters and growth.
The conditional variance in the first term is over growth at fixed
parameters; $\mathbb E_\pi$ then averages it over the urns.
The last variance is over the parameter law, since its argument is
the conditional mean. If the share
converges almost surely to the unique attracting composition
$\randombold{\psi}_*$ of the chosen function, bounded convergence gives
\begin{equation}
\lim_{n\to\infty}\operatorname{Var}(\randombold{\psi}_n)
=\operatorname{Var}(\randombold{\psi}_*).
\label{eq:quenched-limit-variance}
\end{equation}
The variance on the right is over the sampled urn parameters.
Consequently, fluctuations within each system can vanish while
differences between systems persist. Section~\ref{app:quenched}
develops this description for four explicit families.
\par\endgroup

\section{Applications}\label{sec:applications}
\begingroup
\clubpenalty=10000
\widowpenalty=10000
\displaywidowpenalty=10000

The same urn rule can describe different physical quantities. In each
example we identify the binary increment $\sigma_n$, the current share
$\psi_n$ and the quantity counted by $T\varphi_n=n\psi_n$. In direct
representations, the conditional probability of the next increment depends on the history
only through this share. Where a change of variables or of the counting
clock is required, we specify the transformation and the admissible
initial configuration. Connections through count distributions,
time--reversal duality and limiting urn functions are specified separately.

The examples are grouped by field, while their urn functions identify
the common mathematical structure. The regularity and initial-condition
assumptions of Section~\ref{sec:theory} remain attached to the results
in which they are used. Discontinuous urn functions are identified
explicitly. Published models are recalled with their sources; illustrative
rules are distinguished from established descriptions of physical systems.
Section~\ref{app:quenched} then considers persistent differences
between systems through a randomly chosen, frozen urn function.

\subsection{Classical urns, trees, graphs and reinforced walks}
\label{app:classical}

In
one application the black count is a number of balls, in another it
counts leaves of a growing tree, and in another positive steps of a
walk. 

\subsubsection{Polya, Friedman and balanced replacement urns}

In the classical Polya urn, a ball is sampled uniformly, replaced,
and accompanied by another ball of the same color. Here $\sigma_n$
records the color added and the physical black fraction is $\psi_n$.
The urn function is simply
\begin{equation}
\pi(\psi)=\psi.
\end{equation}
Every share is a fixed point. With fixed initial counts of both
colors, the beta--binomial law gives only polynomially small local
probabilities in the interior, hence
\begin{equation}
\phi(x)=0,\qquad 0<x<1.
\label{eq:app-polya-entropy}
\end{equation}
The exact endpoints are inaccessible with this seed, although their
neighborhoods also have zero exponential cost. The limiting share
remains random \cite{MahmoudBook,Pemantle}. A prescribed macroscopic
seed is different: Eq.~(\ref{eq:Caucy}) gives the single trajectory
$\psi(\tau)=\psi_0$ for $0<\tau_0\leq\tau$, rather than a random
limiting share selected during microscopic growth.

With one initial ball of each color, the sequence of added colors also
has the law of independent Bernoulli draws conditional on a single
uniformly chosen propensity \cite{GelastopoulosSageRijt2025}. This
identity concerns the joint law after averaging the propensity.
The same beta--binomial construction models default counts; the
temporal model in \cite{HisakadoMori2020} reduces to P\'olya when past
observations receive equal weights, with positive initial pseudocounts,
that is, starting weights assigned to the two outcomes.
Financial applications are reviewed in \cite{HisakadoFinance2025}.

In Yamamoto's model of football passes \cite{Yamamoto2025}, each of
$M$ possible edges starts with weight one and receives reinforcement
$\alpha>0$ whenever selected. An edge versus its complement is therefore
P\'olya with initial masses $1/\alpha$ and $(M-1)/\alpha$; grouping
edges adds their initial masses. General reinforcement gives real
masses, while the binary initial histories used above require integer
rescaled counts. The identification concerns the proposed pass model.

Friedman's adverse--campaign urn adds a ball of the opposite color
to that sampled. With the same interpretation of the share,
\begin{equation}
\pi(\psi)=1-\psi.
\label{eq:app-friedman-urn}
\end{equation}
The feedback now favors the minority color and the share approaches
$1/2$. For a fixed finite nonempty seed its CGF is
\begin{equation}
\zeta(\beta)=\log[(\exp(\beta)-1)/\beta],\qquad \beta\ne0,
\label{eq:app-friedman-cgf}
\end{equation}
with $\zeta(0)=0$ by continuity \cite{Franchini_URNS,BercuBonnefontRichou2024}.

The OK Corral elimination model is linked to this urn by exact
time reversal \cite[Section~8]{Fajolet2}. Starting with groups of
sizes $m,n\geq1$, a uniformly chosen shooter removes a member of
the opposite group.
For $1\leq s\leq m$, consider a Friedman urn initialized at time $s$
with $s$ black balls and no white balls, and let $\mathbb P$ denote
its probability law. The probability that the OK Corral model started
from $(m,n)$ ends at $(s,0)$ is
\begin{equation}
\frac{s}{m+n}\mathbb P((m+n)\psi_{m+n}=m).
\label{eq:app-ok-corral-duality}
\end{equation}
The Friedman urn makes $m+n-s$ additions. Ratios of opposite
transition probabilities telescope to $s/(m+n)$, so reversed
histories with fixed endpoints have exactly the Friedman bridge law.

More general Friedman and Bagchi--Pal urns add several balls
at a time. Their physical composition must therefore be distinguished
from the binary count used in the HLS representation
\cite{MahmoudBook,Bagchi-Pal,Franchini_URNS}.

Let the replacement matrix have entries $a_{ij}$: after a black
extraction, $a_{11}$ black and $a_{12}$ white balls are added; after
a white extraction, the additions are $a_{21}$ and $a_{22}$. Negative
entries denote removals, when these are allowed by the configuration.
The balanced case has
\begin{equation}
a_{11}+a_{12}=a_{21}+a_{22}=M,\qquad M>0.
\end{equation}
For a tenable urn with $a_{11}\ne a_{21}$, recording black
extractions gives an exact binary representation. If the initial physical counts are $B_0,W_0$,
the corresponding initial HLS clock and black count satisfy
\begin{equation}
m=(B_0+W_0)/M,\qquad
m\psi_m=(B_0-a_{21}m)/(a_{11}-a_{21}).
\end{equation}
To use precisely the binary initial histories defined above, $m$
must be a positive integer and $m\psi_m$ an integer between $0$ and
$m$. Each subsequent black extraction adds one to this transformed
count. At HLS time $n$, the physical black count is
\begin{equation}
B_n=a_{21}n+(a_{11}-a_{21})n\psi_n,
\end{equation}
while the physical total is $Mn$. Thus the probability of the next
black extraction is
\begin{equation}
\pi(\psi)=a_{21}/M+(a_{11}-a_{21})\psi/M.
\end{equation}
This is a finite--time change of variables, not a limiting
approximation. Arbitrary seeds need not satisfy the integer
conditions.
Identical matrix rows give deterministic physical additions and are
separate from the invertible transformation above.

For the symmetric Friedman matrix, with diagonal entries $a_{11}$
and off--diagonal entries $a_{12}$, this gives
\begin{equation}
\pi(\psi)=a_{12}/(a_{11}+a_{12})
+(a_{11}-a_{12})\psi/(a_{11}+a_{12}).
\end{equation}
For nonnegative entries with $a_{11}+a_{12}>0$ and
$a_{11}\ne a_{12}$, the physical seed
$B_0=W_0=a_{11}+a_{12}$ corresponds to $m=2$ and $m\psi_m=1$.
For example, diagonal entries $2$ and off--diagonal entries $1$
give $\pi(\psi)=(1+\psi)/3$.

Subtractive replacement urns may instead produce an affine function
only on the accessible interval. Outside that interval it can be
extended by truncation at zero and one. The tenable Bagchi--Pal
example with diagonal entries $-1$, off--diagonal entries $2$, and
physical seed $(1,1)$ has
\begin{equation}
\pi(\psi)=\min\{1,\max\{0,2-3\psi\}\}.
\end{equation}
The transformed process stays in the region where the physical
counts are nonnegative. These replacement urns also arise in the
analysis of data structures and branching constructions
\cite{Bagchi-Pal,Fajolet-Analytic=000020Urns,Fajolet2,Stochastic=000020urns}.
Quantitative concentration estimates for balanced replacement urns
are given in \cite{Dasgupta2024}; affine multicolor urns with multiple
drawings, beyond the binary representation here, are studied in
\cite{SparksKubaBalajiMahmoud2025}.

These affine rules share one calculation. Write
$\pi(\psi)=a+b\psi$ on the accessible interval. For $b\ne1$,
Eq.~(\ref{eq:Caucy}) gives
\begin{equation}
\psi(\tau)=a/(1-b)+[\psi_0-a/(1-b)](\tau_0/\tau)^{1-b},
\qquad 0<\tau_0\leq\tau.
\label{eq:app-affine-trajectory}
\end{equation}
The formula applies while the trajectory remains in that interval;
$\varphi(\tau)=\tau\psi(\tau)$. For $b<1$, the attracting share is
$a/(1-b)$. Adverse--campaign Friedman has $b=-1$ and attractor $1/2$;
the symmetric replacement matrix has
$b=(a_{11}-a_{12})/(a_{11}+a_{12})$ and the same attractor when
$a_{12}>0$. The Bagchi--Pal example has $b=-3$: its displacement
from $1/2$ decays as $(\tau_0/\tau)^4$.

Figure~\ref{fig:review-classical-urns}\label{call:fig:review-classical-urns} contrasts P\'olya neutrality with convergence towards $1/2$ for the symmetric replacement matrix with diagonal entries $2$ and off--diagonal entries $1$. In the latter model, $\psi$ records black extractions and the physical black fraction is $(1+\psi)/3$.

\afterpage{%
\clearpage
\begin{figure}[p]
\centering
\includegraphics[width=\linewidth]{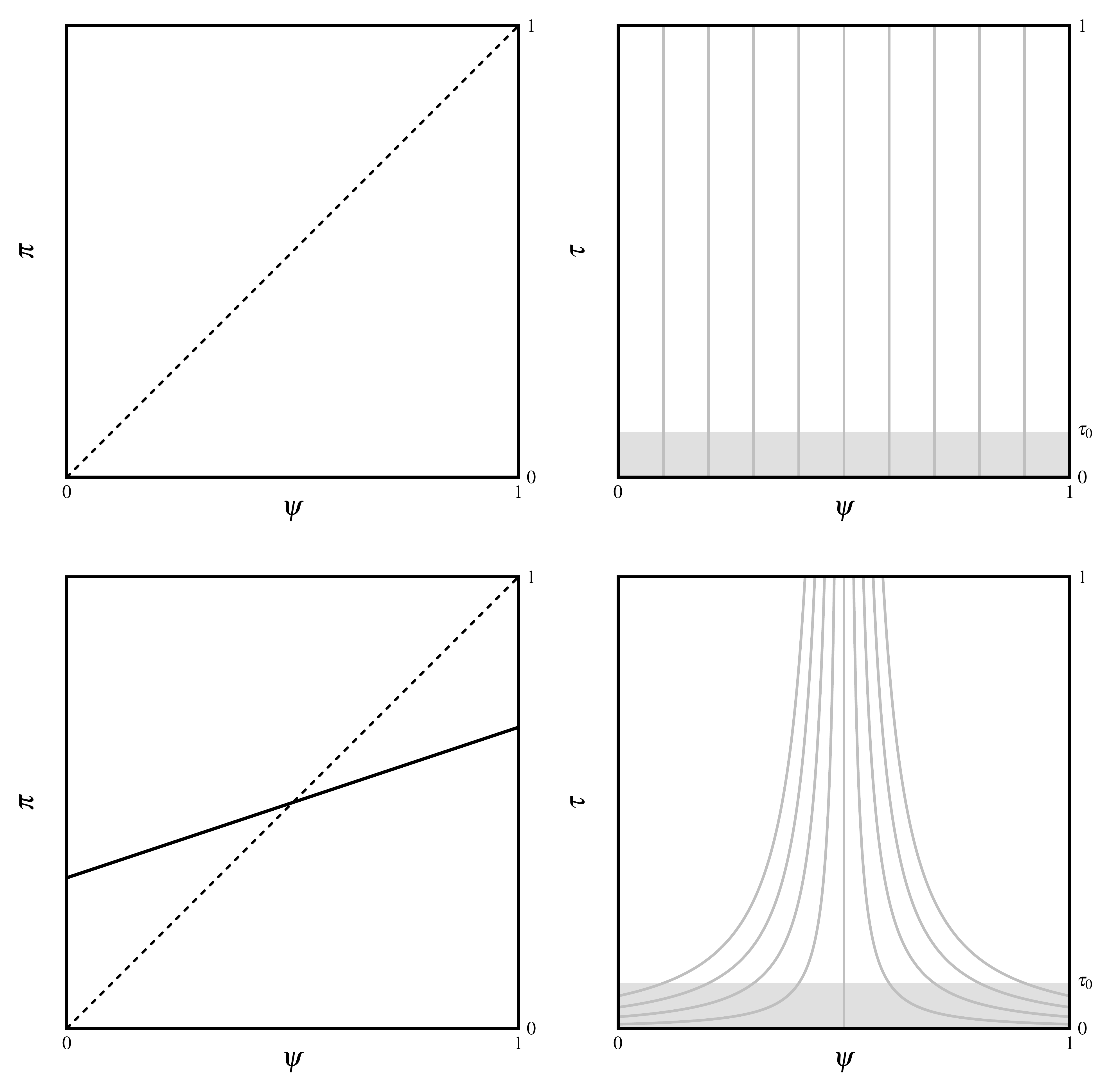}
\caption{Urn functions (left) and zero--cost share trajectories (right) for P\'olya's rule (top) and the symmetric replacement urn with diagonal entries $2$ and off--diagonal entries $1$ (bottom). For P\'olya, the response coincides with the dashed diagonal and shares remain constant. The affine response $\pi(\psi)=(1+\psi)/3$ instead attracts towards $1/2$. In both rows, gray curves pass through $\psi_0=0.1,0.2,\ldots,0.9$ at $\tau_0=0.1$. The lower trajectories follow Eq.~(\ref{eq:app-affine-trajectory}) with $a=b=1/3$ and remain distinct at $\tau=1$. Axes cover $[0,1]$, with saturation increasing upwards. Within the shaded strip $0\leq\tau<\tau_0$, the curves show backward continuations of the same equations.}
\label{fig:review-classical-urns}
\end{figure}
\clearpage
}

For a fixed finite initial history, the equation for $\zeta$ in
Section~\ref{sec:entropy}, wherever differentiable, becomes
\begin{equation}
[a+b\partial_\beta\zeta(\beta)](\exp(\beta)-1)=\exp(\zeta(\beta))-1,
\qquad \beta\ne0.
\label{eq:app-affine-cgf}
\end{equation}
This equation is linear in $\exp(-\zeta)$. The following integral
formulas solve the affine model and its Bagchi--Pal representation
\cite{Franchini_URNS}, including the corrigendum cited there. For $0<b<1$ and $0<a<1-b$, the positive
and negative half--axes have the following integral forms. For $\beta>0$,
\begin{equation}
\exp(-\zeta(\beta))
=\frac{\exp(-\beta)(1-\exp(-\beta))^{1/b}}{b}
\int_0^1 t^{a/b}(1-t\exp(-\beta))^{-1-1/b}\,dt.
\label{eq:app-affine-cgf-positive}
\end{equation}
Exchanging the colors gives, for $\beta<0$,
\begin{equation}
\exp(-\zeta(\beta))
=\frac{(1-\exp(\beta))^{1/b}}{b}
\int_0^1 t^{(1-a-b)/b}(1-t\exp(\beta))^{-1-1/b}\,dt.
\label{eq:app-affine-cgf-negative}
\end{equation}
For $b<0$ and $0<a<1-b$, both signs of $\beta$ are covered by
\begin{equation}
\exp(-\zeta(\beta))
=\int_0^1[t^{-b}+(1-t^{-b})\exp(\beta)]^{a/b}\,dt.
\label{eq:app-affine-cgf-decreasing}
\end{equation}
When the affine response leaves $[0,1]$, this last formula refers
to its truncation at zero and one, with the accessible histories
of the initial seed. It includes the tenable subtractive
Bagchi--Pal example above: $a=2,b=-3$ gives
\begin{equation}
\exp(-\zeta(\beta))
=\int_0^1[t^3+(1-t^3)\exp(\beta)]^{-2/3}\,dt.
\label{eq:app-bagchi-pal-cgf}
\end{equation}
All these formulas have the continuous value $\zeta(0)=0$.
For $b<0$, boundedness at zero selects the solution of the linear
equation. For $b>0$, the boundary conditions at positive and
negative infinity select the two branches; the value at zero alone
does not do so. In the increasing case these conditions are
\begin{equation}
\lim_{\beta\to+\infty}[\zeta(\beta)-\beta]=\log(a+b),
\qquad \lim_{\beta\to-\infty}\zeta(\beta)=\log(1-a).
\label{eq:app-affine-cgf-boundaries}
\end{equation}

These CGFs count black extractions in the HLS representation.
The physical black count follows from the exact affine change of
variables, so its CGF per HLS step is
\begin{equation}
\lim_{T\to\infty}T^{-1}\log\mathbb E(\exp(\beta B_T))
=a_{21}\beta+\zeta((a_{11}-a_{21})\beta).
\label{eq:app-bagchi-pal-physical-cgf}
\end{equation}
Thus the same solution determines both observables, with their
different normalizations.

For $b=0$ and $0<a<1$, the subsequent choices are independent and
\begin{equation}
\zeta(\beta)=\log(1-a+a\exp(\beta)).
\label{eq:app-bernoulli-cgf}
\end{equation}
The interior entropy is then
\begin{equation}
\phi(x)=-x\log(x/a)-(1-x)\log[(1-x)/(1-a)],\qquad 0<x<1.
\label{eq:app-bernoulli-entropy}
\end{equation}
Endpoint probabilities still depend on which colors the initial law
allows. In general the path variational principle retains the
accessible histories; Legendre duality alone does not recover a
nonconcave entropy or remove restrictions caused by zeros of $\pi$.
Section~\ref{app:quenched-affine} uses these solutions for a
Bagchi--Pal matrix chosen once and kept fixed, and calculates the
quenched average for two explicit laws of the coefficients.

\subsubsection{Growing trees and permutation statistics}

Tree growth supplies observables with binary increments even though
the full tree has a much larger state space. In a uniformly recursive tree, leaves
are vertices without children. Each new vertex attaches to a uniformly
selected existing vertex. The number of leaves increases by one if the chosen vertex is not a
leaf, and remains unchanged otherwise. With $n$ vertices and leaf
count $n\psi_n$, the exact urn function is
\begin{equation}
\pi(\psi)=1-\psi.
\end{equation}
This identifies the entire leaf--count process with an HLS urn,
without identifying the tree itself with a binary history
\cite{Bryc}.
Equation~(\ref{eq:app-affine-trajectory}) has $b=-1$: the limiting
leaf share is $1/2$, with deterministic relaxation exponent $2$.

Uniform insertion into a random permutation gives the same urn.
If $D_n$ counts adjacent pairs in decreasing order (descents), the
count $D_n+1$ in clock $n+1$
follows Eq.~(\ref{eq:app-friedman-urn}), starting from one ball of
each color. Its CGF is Eq.~(\ref{eq:app-friedman-cgf}); sharp large
deviations and concentration bounds are obtained in
\cite{BercuBonnefontRichou2024}.

In plane--oriented recursive trees, the children of each vertex are
ordered. A new vertex is inserted into a position chosen uniformly
from all positions before, between or after the children of existing
vertices. The leaf count $n\psi$ increases with probability
\begin{equation}
\pi_n(\psi)=1-n\psi/(2n-1).
\label{eq:app-plane-tree-urn}
\end{equation}
In a Stirling permutation of $\{1,1,\ldots,n,n\}$, the entries
between equal symbols are larger than those symbols. A plateau is
an adjacent equal pair. Inserting the next equal pair into a uniformly
chosen gap makes the plateau count increase with probability
\begin{equation}
\pi_n(\psi)=1-n\psi/(2n+1).
\label{eq:app-stirling-urn}
\end{equation}
The respective seeds are the single--vertex tree and the word $11$.
Dividing the leaf and plateau counts by the shifted totals $n-1/2$
and $n+1/2$, respectively, removes time dependence from the urn
functions. These totals are half--integral; for the tree, the first
deterministic insertion precedes this rescaling. In the ordinary clock, both kernels converge
uniformly at order $1/n$ to $\pi(\psi)=1-\psi/2$.
The path LDP and CGF are those of this limiting urn \cite{Bryc};
this result includes the boundary $\pi(0)=1$. Its fixed share is
$2/3$ and its relaxation exponent is $3/2$ in
Eq.~(\ref{eq:app-affine-trajectory}).

Under preferential attachment, each new vertex connects by one edge
to an existing vertex chosen with probability proportional to its
weight $d+\rho$, where $d$ is the degree and $\rho>-1$.
Leaves have degree one. Suppose the
initial graph has no isolated vertices, $V_0$ vertices and $E_0$ edges.
Set the initial HLS clock to
\begin{equation}
m=(2E_0+\rho V_0)/(2+\rho),
\end{equation}
then increase it by one with each new vertex and define $\psi_n$
as the leaf count divided by this shifted clock. If $m$ is a
positive integer and the initial leaf count is at most $m$, the
result is the autonomous HLS function
\begin{equation}
\pi(\psi)=1-(1+\rho)\psi/(2+\rho).
\end{equation}
For example, a cycle used as the initial graph has $E_0=V_0$ and
no leaves, so these conditions hold for every $\rho>-1$. The
leaf--count process has an exact HLS representation in this clock
\cite{Bryc}.
Here $b=-(1+\rho)/(2+\rho)$, so the attracting leaf share in the
shifted clock is $(2+\rho)/(3+2\rho)$ and its relaxation exponent
is $(3+2\rho)/(2+\rho)$.
For a general fixed seed with no isolated vertices, retain instead
the vertex count $n$ as clock. The exact kernel is
\begin{equation}
\pi_n(\psi)=1-\frac{(1+\rho)n\psi}
{(2+\rho)n+2(E_0-V_0)}.
\label{eq:app-preferential-finite}
\end{equation}
It converges uniformly on accessible shares, at order $1/n$, to the
same affine function. The seed offset does not change its path LDP
or CGF \cite{Bryc}, even when the shifted initial clock above is
not an admissible integer history.

Another exact example is the number of cherries, or pairs of leaves
with a common parent, in a Yule tree. Each step splits a uniformly
chosen leaf into two daughters. If the tree has $n$ leaves and
$n\psi_n$ cherries, a new cherry is gained unless the selected leaf
already belongs to a cherry. Hence
\begin{equation}
\pi(\psi)=\max\{0,1-2\psi\}.
\end{equation}
Only shares up to $1/2$ are accessible. Starting from one leaf and
no cherries corresponds to an initial white ball. This example
connects the urn directly to a shape statistic of a branching tree
\cite{Bryc}.
On its accessible interval, $b=-2$ gives attracting cherry share
$1/3$ and relaxation exponent $3$ in
Eq.~(\ref{eq:app-affine-trajectory}).
Quantitative convergence results for degree statistics of recursive,
plane recursive and $d$--ary trees, using multicolor urns, are
developed in \cite{NeriPinto2026}.

\subsubsection{Elephant random walks}

In the classical Elephant Random Walk, the next step repeats a
uniformly sampled past step with probability $p$, and reverses it
otherwise. Identify a positive step with $\sigma_n=1$ and a negative
step with $\sigma_n=0$. The positive fraction is $\psi_n$, and the
walk's position is $2n\psi_n-n$. Its exact urn function is
\begin{equation}
\pi(\psi)=(1-p)+(2p-1)\psi,\qquad 0\leq p\leq1.
\label{eq:app-erw-urn}
\end{equation}
The first step has its separately specified initial law. This
representation connects the memory rule to classical urn methods
\cite{ERW=000020shcutz=000020trimper,ERW=000020UM=000020Baur=000020Berton,FB}.
For $p<1$, the affine slope is $b=2p-1$ and the attracting positive
share is $1/2$. The scaled displacement along the trajectory is
\begin{equation}
2\varphi(\tau)-\tau
=(2\psi_0-1)\tau_0^{2(1-p)}\tau^{2p-1}.
\end{equation}
At $p=1$ the scaling trajectory preserves the initial composition;
a mixed microscopic seed instead gives the Polya behavior above.
With a single initial step, every subsequent step has its sign.
The center of mass is a different observable of the same history,
not a different urn function:
\begin{equation}
\mathcal O(\sigma)=T^{-1}\sum\limits_{1\leq n\leq T}(2n\psi_n-n).
\end{equation}
Its scaling limit is
\begin{equation}
\mathcal O(\sigma)/T\longrightarrow
\int_0^1[2\varphi(\tau)-\tau]\,d\tau,\qquad T\to\infty.
\end{equation}
For a fixed microscopic initial law and a terminal-displacement interval
$[u,v]\subset[-1,1]$, the path principle gives, whenever the
large--deviation upper and lower bounds coincide,
\begin{equation}
\lim_{T\to\infty}T^{-1}\log\mathbb P(2\psi_T-1\in[u,v])=
\sup_{\varphi\in Q:\,2\varphi(1)-1\in[u,v]}
\{\Phi(\varphi)-\Phi_0^*(\varphi)\}.
\label{eq:app-erw-contraction}
\end{equation}
For the scaled center of mass, replace the endpoint constraint by
the same interval constraint on the integral above. Both functionals
are continuous in the uniform norm, and $Q$ is restricted to
accessible histories. A prescribed macroscopic seed
instead fixes $\varphi(\tau_0)=\tau_0\psi_0$ and uses the actions
integrated only from $\tau_0$ to $1$. For the center of mass, its
assigned prefix also contributes to the integral from $0$ to $\tau_0$
and must be retained. For $0<p<1$ the urn
probabilities stay positive; degenerate endpoint rules require the
accessibility treatment already stated in Section~\ref{sec:ldp}.
Fluctuations of the center of mass and position have been studied in
\cite{Bercu,Jack=000020Harris}. Quantitative Gaussian approximations,
moderate deviations and local limits refine the diffusive and critical
results \cite{FanHuMa2021}. In the superdiffusive regime, distributional
fixed--point equations characterize the limiting law, including its
density, moments and tails
\cite{GuerinLaulinRaschel2026,GuerinLaulinRaschelSimon2025}.
Recent preprints further study local limits and the dependence of the
density's shape on the initial bias \cite{TrinasValle2026,Pemmasani2026}.
A prescribed initial run also affects subsequent return times
\cite{Fang2026}.

The finite--time characteristic functions have a polynomial recurrence
\cite{GuerinLaulinRaschel2025}. A generating--function approach using
Schwarz--Christoffel mappings has also been proposed for the Laplace
transform \cite{GonzalezRangel2026}. These transform methods complement
the affine solution above; the CGF used here retains normalization by
the number of steps, whereas the superdiffusive limit law has its own
transform after rescaling the position.

The Bernoulli elephant random walk instead has steps zero and one:
a uniformly sampled past step is repeated with probability $p$, and
zero is added otherwise \cite{GutStadtmuller2026}. With full memory,
this is the affine urn $a=0,b=p$. For $0<p<1$, conditional on a first
step equal to one, the count divided by $n^p$ converges to a
Mittag--Leffler law. A zero first step is absorbing. The share tends
to zero in either case, although the positive count can keep growing;
Eq.~(\ref{eq:app-affine-trajectory}) gives the corresponding
macroscopic trajectory.

For multiple extractions, let $k$ past steps be sampled independently
with replacement. If $k$ is odd and the walker follows the majority
with probability $p$, the urn function is
\begin{equation}
\pi(\psi)=(1-p)+(2p-1)\sum\limits_{(k+1)/2\leq j\leq k}
\binom{k}{j}\psi^j(1-\psi)^{k-j}.
\end{equation}
In particular, three extractions give
\begin{equation}
\pi(\psi)=(1-p)+(2p-1)(3\psi^2-2\psi^3).
\end{equation}
These are exact nonlinear HLS urns for fixed $k$, including pure
majority at $p=1$ \cite{Fra2022}. The same sampling mechanism
also describes adoption decisions based on a sample of past customers.

More generally, let $f:[0,1]\rightarrow[0,1]$ specify the response
to the positive fraction in the sample. The replacement rule of
\cite{Podder} gives
\begin{equation}
\pi(\psi)=(1-p)+(2p-1)\sum\limits_{0\leq j\leq k}
f(j/k)\binom{k}{j}\psi^j(1-\psi)^{k-j}.
\label{eq:app-sample-response}
\end{equation}
Thus the sample rule is averaged over a binomial distribution.
Majority is recovered by taking $f(u)=\mathbb I(u>1/2)$ for odd
$k$. More general multiple--drawing urns admit representations by
labelled directed acyclic graphs \cite{MaillerSteiner2025}.

Maulik and collaborators \cite{MaulikRoySadhukhan2026} apply $f$
directly to the current share. Their generalized ERW and their minimal
walk with steps zero and one have, respectively,
\begin{equation}
\pi(\psi)=1-p+(2p-1)f(\psi),\qquad
\pi(\psi)=q+(p-q)f(\psi).
\label{eq:app-direct-responses}
\end{equation}
In the minimal walk, a tentative one, drawn with probability $f(\psi)$,
produces a one with probability $p$; a tentative zero does so with
probability $q$. Both parameters lie in $(0,1)$. These exact HLS rules
have convergence and fluctuation results under the regularity and
stability assumptions of that work. The direct response differs from
its finite--sample binomial average.

Without replacement, for $n\geq k$ and integer $n\psi$, the exact
conditional probability is instead
\begin{equation}
\pi_n(\psi)=1-p+(2p-1)\sum\limits_{0\leq j\leq k}
f(j/k)\frac{\binom{n\psi}{j}\binom{n(1-\psi)}{k-j}}{\binom{n}{k}}.
\label{eq:app-hypergeometric-response}
\end{equation}
Impossible samples have zero binomial coefficient. Here $\pi$ denotes
the binomial response in Eq.~(\ref{eq:app-sample-response}), with the
same fixed sample size $k$. Couple the two samples by using the same indices until the sample
with replacement first repeats an index. Before that collision the
responses agree. Bounding its probability by the sum over pairs of
draws gives, for fixed $k$,
\begin{equation}
\max\limits_{0\leq j\leq n}|\pi_n(j/n)-\pi(j/n)|
\leq |2p-1|\frac{k(k-1)}{2n}.
\label{eq:app-sampling-bound}
\end{equation}
For $0<p<1$, Section~\ref{sec:asymptotic-urns} transfers the path LDP
and CGF to this limit. For pure majority, $p=1$ and odd $k$, the same
conclusion holds from a common seed with at least $(k+1)/2$ steps of each
sign. The logarithm of the ratio between the two conditional probabilities
of the same allowed addition is bounded by a constant depending on
$k$ divided by the current count of the sign added. Summing over
successive additions gives a uniform $O(\log T)$ comparison of path
log--probabilities. The seed condition ensures that both sampling
schemes allow exactly the same histories (common support).

If $k_n\to\infty$ and $f$ is continuous, either sampling scheme
converges uniformly to the direct ERW response above; without replacement
one requires $k_n\leq n$. The sample fraction has variance at most
$1/(4k_n)$, and uniform continuity of $f$ gives the claim. No condition
$k_n^2/n\to0$ is needed. For $0<p<1$, the same LDP and CGF follow when
the limiting response satisfies the regularity assumptions of
Section~\ref{sec:ldp}.

\subsubsection{Genealogy and information transmission in growing graphs}
\label{app:graphs}

In Burghart's recursive directed graph \cite{Burghart2026}, each new
vertex chooses $k$ earlier vertices uniformly with replacement, with
arrows directed from new to old. Fix a vertex $i\geq2$ and let
$n\psi_n$ count the vertices from which $i$ can be reached, including
$i$ itself. At $n=i$ the urn contains one black and $i-1$ white balls.
A new vertex is counted precisely when at least one chosen parent is
counted, giving
\begin{equation}
\pi(\psi)=1-(1-\psi)^k.
\label{eq:app-genealogy-urn}
\end{equation}
This is Eq.~(\ref{eq:app-sample-response}) with $p=1$ and
$f(u)=\mathbb I(u>0)$. For $k=1$ it is P\'olya; for $k>1$ every
interior zero--cost share increases. For $k>1$ and a macroscopic seed
$0<\psi_0<1$ at $\tau_0>0$, Eq.~(\ref{eq:Caucy}) integrates to
\begin{equation}
[1-\psi(\tau)]^{1-k}-1
=[(1-\psi_0)^{1-k}-1](\tau/\tau_0)^{k-1}.
\label{eq:app-genealogy-trajectory}
\end{equation}
The share approaches one as the growth ratio increases without bound;
at $\tau=1$ it remains below one. Burghart also studies seeds whose
size grows with the final graph, requiring their own scaling limits.
For a fixed microscopic seed, the zero and one boundary values retain
the accessibility conditions of Section~\ref{sec:ldp}.

Broadcasting assigns binary colors to growing vertices. In the model
of Briend, Devroye and Lugosi \cite{BriendDevroyeLugosi2025}, start with
$k$ colored roots, with $k$ odd. Each new vertex samples $k$ parents
uniformly with replacement, independently transmits each sampled color
correctly with probability $p\in[1/2,1]$, and takes their majority.
Writing $\psi_n$ for the black fraction gives the exact rule
\begin{equation}
\pi(\psi)=\sum\limits_{(k+1)/2\leq j\leq k}\binom{k}{j}
[1-p+(2p-1)\psi]^j[p-(2p-1)\psi]^{k-j}.
\label{eq:app-broadcast-urn}
\end{equation}
Repeated selections have independent transmission errors.
The sampling is with replacement, as in Eq.~(\ref{eq:app-sample-response}).
Here, however, $p$ is the probability of correctly transmitting each
sampled color before taking the majority; in
Eq.~(\ref{eq:app-sample-response}), it acts on the response to the sample.
There is no subsequent reversal of the majority.
For $k\geq3$
and $1/2<p<1$, it differs from reversing the majority after sampling;
the two rules coincide at $p=1/2$ and $p=1$.
For $k=3$, $\partial_\psi\pi(1/2)=3(2p-1)/2$: the central share attracts below
$p=5/6$ and becomes unstable above it, when two attracting shares appear:
\begin{equation}
\psi_{\pm}=[1\pm\sqrt{(6p-5)/(2p-1)^3}]/2,
\qquad 5/6<p\leq1.
\label{eq:app-broadcast-fixed}
\end{equation}
Their slope is $6(1-p)$, so Eq.~(\ref{eq:app-local-relaxation}) gives
the same local relaxation exponent $6p-5$ as for the DEK outer branches.
The published reconstruction threshold differs from this equilibrium
threshold: the final majority can retain information about the initial
majority even when $\psi_n\to1/2$, in particular for $k=3$ and
$2/3<p<5/6$. For $k=3$, the threshold $p=2/3$ for reconstruction
by the final majority corresponds to $\partial_\psi\pi(1/2)=1/2$
\cite{BriendDevroyeLugosi2025}. The classical elephant random walk
has this same slope at its diffusive threshold $p=3/4$
\cite{FanHuMa2021}. This comparison concerns the response near
the balanced share. For $k=1$, the same color count is
the affine urn in Eq.~(\ref{eq:app-erw-urn}), studied on recursive trees
in \cite{AddarioBerry2022,AlthausHartungSteiner2026}. Preferential
attachment extensions also track attachment weights
\cite{DesmaraisHolmgrenWagner2023,AlthausHartungSteiner2026}.

\subsection{Economics and social influence}\label{app:economics}

Arthur's increasing returns theory (IRT) provides a direct interpretation
of the HLS mechanism. A choice changes the environment in which later
choices are made: an additional user can make a technology more useful,
more familiar, or easier to adopt. A small initial advantage may therefore
persist and grow. The urn records successive choices through $\sigma_n$,
while $\psi_n$ is the fraction choosing the first alternative. Its function
$\pi$ specifies how this accumulated experience affects the next choice.
This connection motivated the generalized urn models of Arthur, Ermoliev
and Kaniovskii \cite{AEK,Arthur=000020Ermoliev=000020Kaniovski} and the
subsequent analysis of competing technologies and path dependence
\cite{Arthur,Arthur=000020book,ArthNat}. The associated convergence theory
was developed alongside these applications
\cite{Ermoliev=000020Arthur1,Ermoliev-Arthur2,Ermoliev-Arthur3}.
The initial share specifies the advantage already accumulated when
observation begins. The saturation specifies how much further growth
remains. Keeping these two quantities separate is useful when comparing
a new market with an established one: identical current shares need
not imply identical remaining opportunities for fluctuations to change
the eventual outcome.

A particularly transparent example is majority adoption. Each new
customer consults an odd number $k$ of previous customers, sampled
independently with replacement, and adopts the product used by the
majority of this sample.
Figure~\ref{fig:Maj1}\label{call:fig:Maj1} sketches the three--customer sampling rule.

\penalty0
\afterpage{%
\clearpage
\begin{figure}[p]
\centering{}{\Huge\includegraphics[scale=0.23]{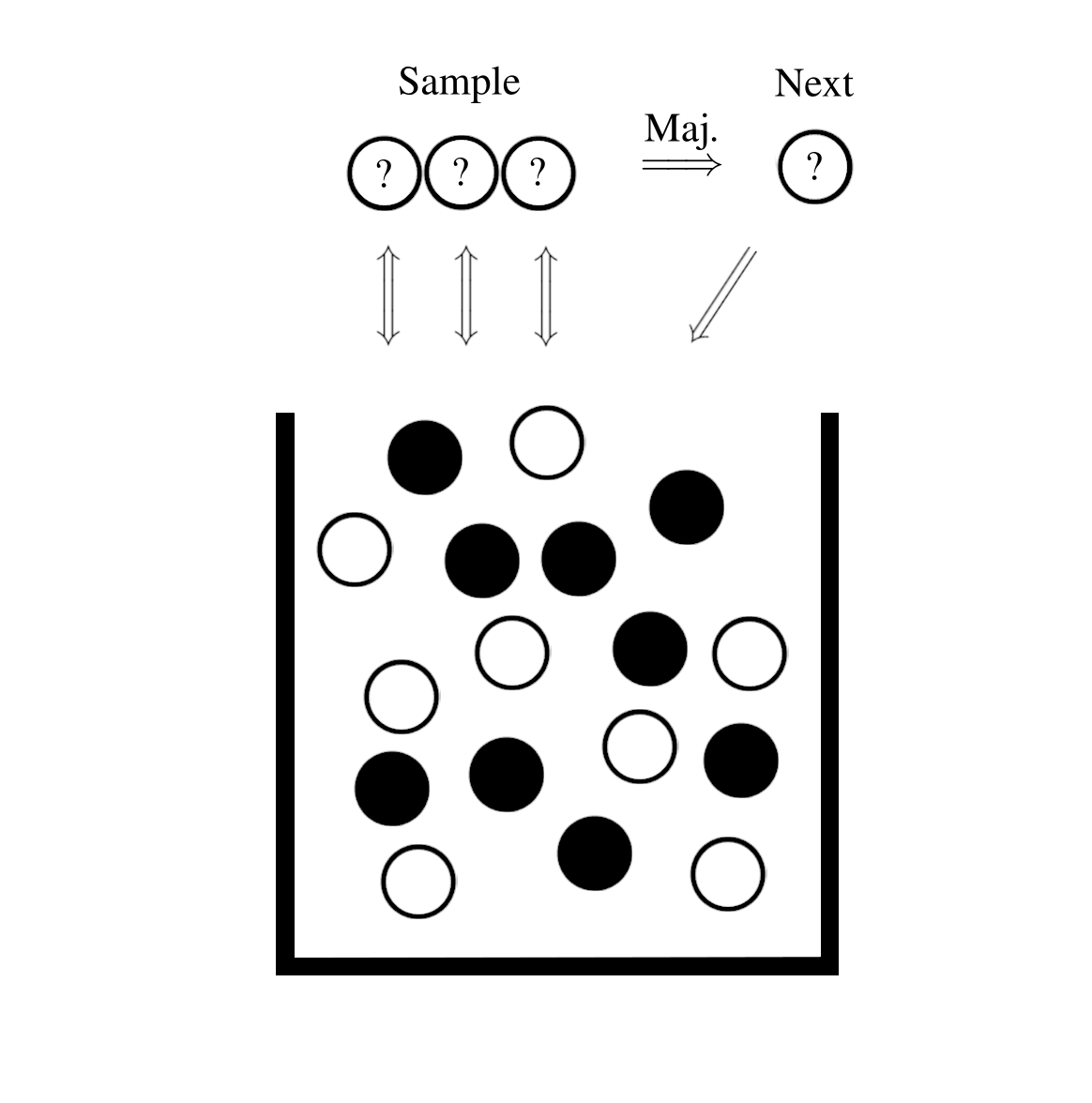}}\caption{Two competing products gain customers through a majority
mechanism \cite{AEK,FB,ArthNat}: each new customer samples an odd number of previous
customers (at least three), independently with replacement, and adopts
the product chosen by the majority. In figure: Scheme for the step
of the IRT model.}
\label{fig:Maj1}
\end{figure}
\clearpage
}

 Conditional on the current share, the number
of sampled users of the first product is binomial. Thus
\begin{equation}
\pi(\psi)=\sum\limits_{(k+1)/2\leq j\leq k}\binom{k}{j}\psi^j(1-\psi)^{k-j}.
\label{eq:app-majority}
\end{equation}
For three consultations this becomes
\begin{equation}
\pi(\psi)=3\psi^2-2\psi^3.
\label{eq:app-majority-three}
\end{equation}
Figure~\ref{fig:Maj2}\label{call:fig:Maj2} shows how the three sampled choices produce this HLS response.

\afterpage{%
\clearpage
\begin{figure}[p]
\begin{centering}
{\Huge\includegraphics[angle=-90,scale=0.24]{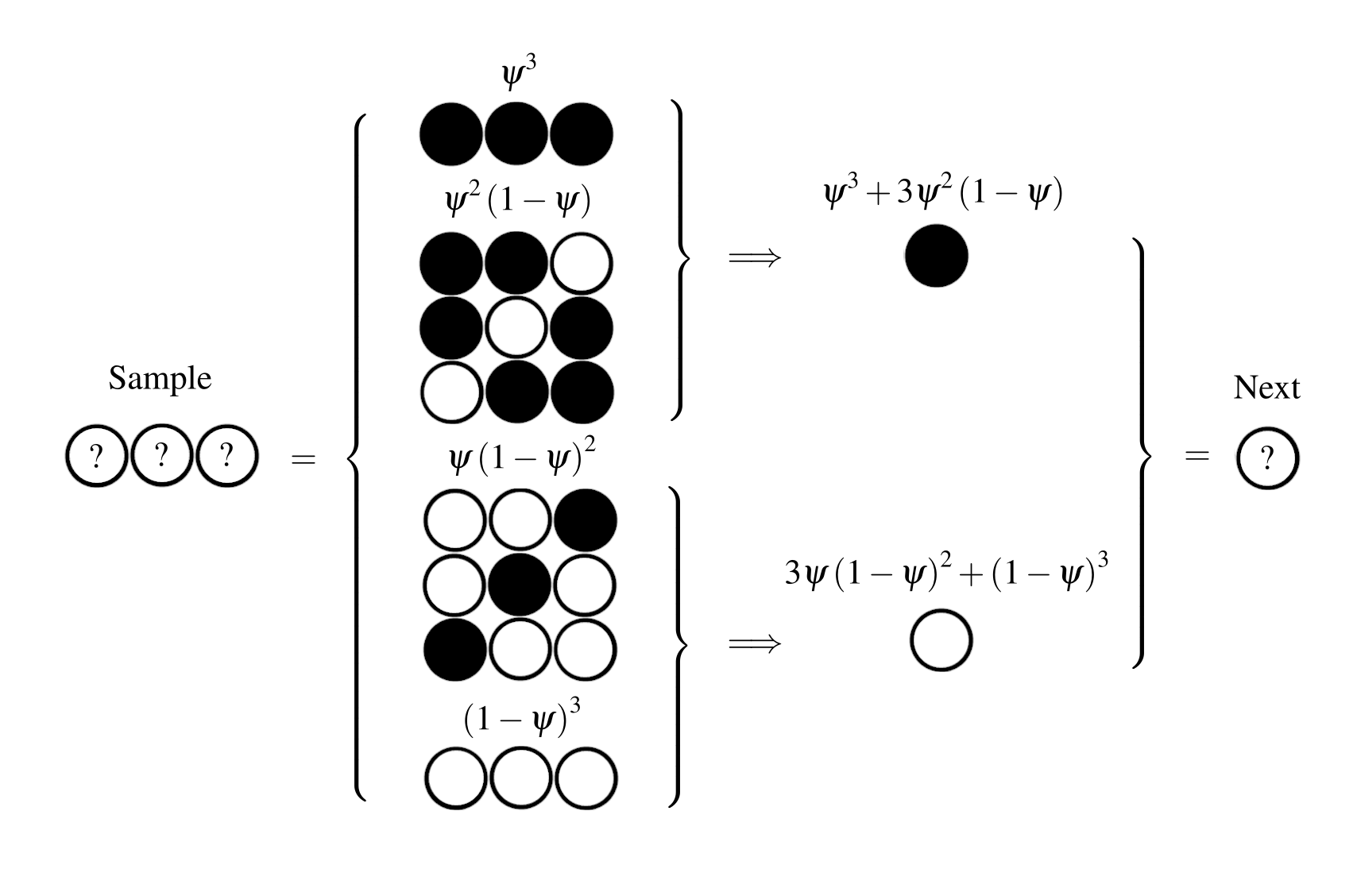}}{\Huge\par}
\par\end{centering}
\begin{centering}
~
\par\end{centering}
\begin{centering}
~
\par\end{centering}
\centering{}~\caption{Reduction of three-customer majority sampling to the HLS rule. The eight ordered samples are grouped by majority color. Independent sampling with replacement gives the black-addition probability in Eq.~(\ref{eq:app-majority-three}).}
\label{fig:Maj2}
\end{figure}
\clearpage
}

The curve lies below the diagonal for $0<\psi<1/2$ and above it for
$1/2<\psi<1$. Both products being initially present, the limiting share
is zero or one: the initial fluctuations help select which product
eventually dominates. The relevant feedback concerns the composition
of the market, rather than a fixed preference for either product.

Dosi, Moneta and Stepanova \cite{DMS} use this majority rule for additions.
Figure~\ref{fig:Maj3}\label{call:fig:Maj3} shows the urn function and
trajectories of the corresponding HLS growth process.

\afterpage{%
\clearpage
\begin{figure}[p]
{\Huge\includegraphics[width=\linewidth]{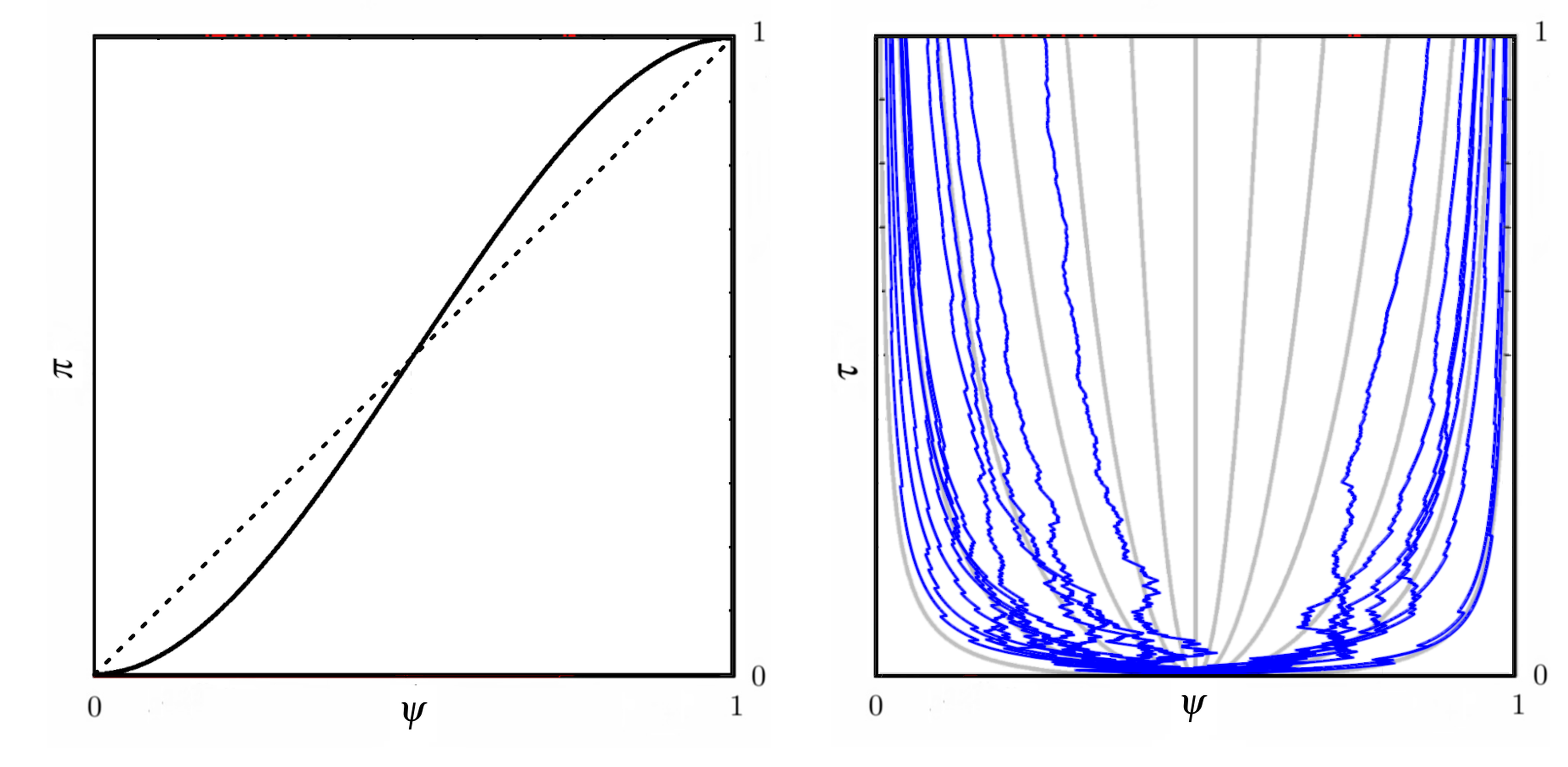}}\caption{Urn function (left) and zero--cost share trajectories (right, gray) for the IRT model. The blue curves are
simulated trajectories of the addition process in Dosi et al. 2019
\cite{DMS}. Both colors are present initially.}
\label{fig:Maj3}
\end{figure}
\clearpage
}

\begin{samepage}
The Dosi--Ermoliev--Kaniovsky (DEK) model allows customers to disagree with the sample
\cite{Dosi=000020Ermoliev=000020Kaniovski,FB}. Let $p$ be the probability
of following its majority, with the opposite product chosen otherwise.
The resulting function is\par
\end{samepage}
\begin{equation}
\pi(\psi)=1-p+(2p-1)
\sum\limits_{(k+1)/2\leq j\leq k}\binom{k}{j}\psi^j(1-\psi)^{k-j}.
\label{eq:app-dek}
\end{equation}
For $k=1$ the response is affine,
\begin{equation}
\pi(\psi)=1-p+(2p-1)\psi.
\end{equation}
For the same initial history, this is exactly the classical elephant
random walk response in Eq.~(\ref{eq:app-erw-urn}), with customer
choices replacing step signs. For $k=3$ the response is
\begin{equation}
\pi(\psi)=1-p+(2p-1)(3\psi^2-2\psi^3).
\end{equation}
At $p=1/2$ the choices are independent; at $p=1$ we recover majority
adoption. For three consultations the additional fixed shares are
\begin{equation}
\psi_{\pm}=(1\pm\sqrt{(6p-5)/(2p-1)})/2,
\qquad 5/6<p\leq1.
\label{eq:app-dek-fixed}
\end{equation}
The central share has slope $\partial_\psi\pi(1/2)=3(2p-1)/2$; at the two
outer shares $\partial_\psi\pi(\psi_{\pm})=6(1-p)$. Thus $1/2$ attracts for
$p<5/6$, whereas the two outer shares attract for $p>5/6$.
They remain interior when $p<1$: dominance need not mean monopoly.
Figure~\ref{fig:Urna}\label{call:fig:Urna} illustrates the two interior attractors and their zero--cost trajectories.

\afterpage{%
\clearpage
\begin{figure}[p]
{\Huge\includegraphics[width=\linewidth]{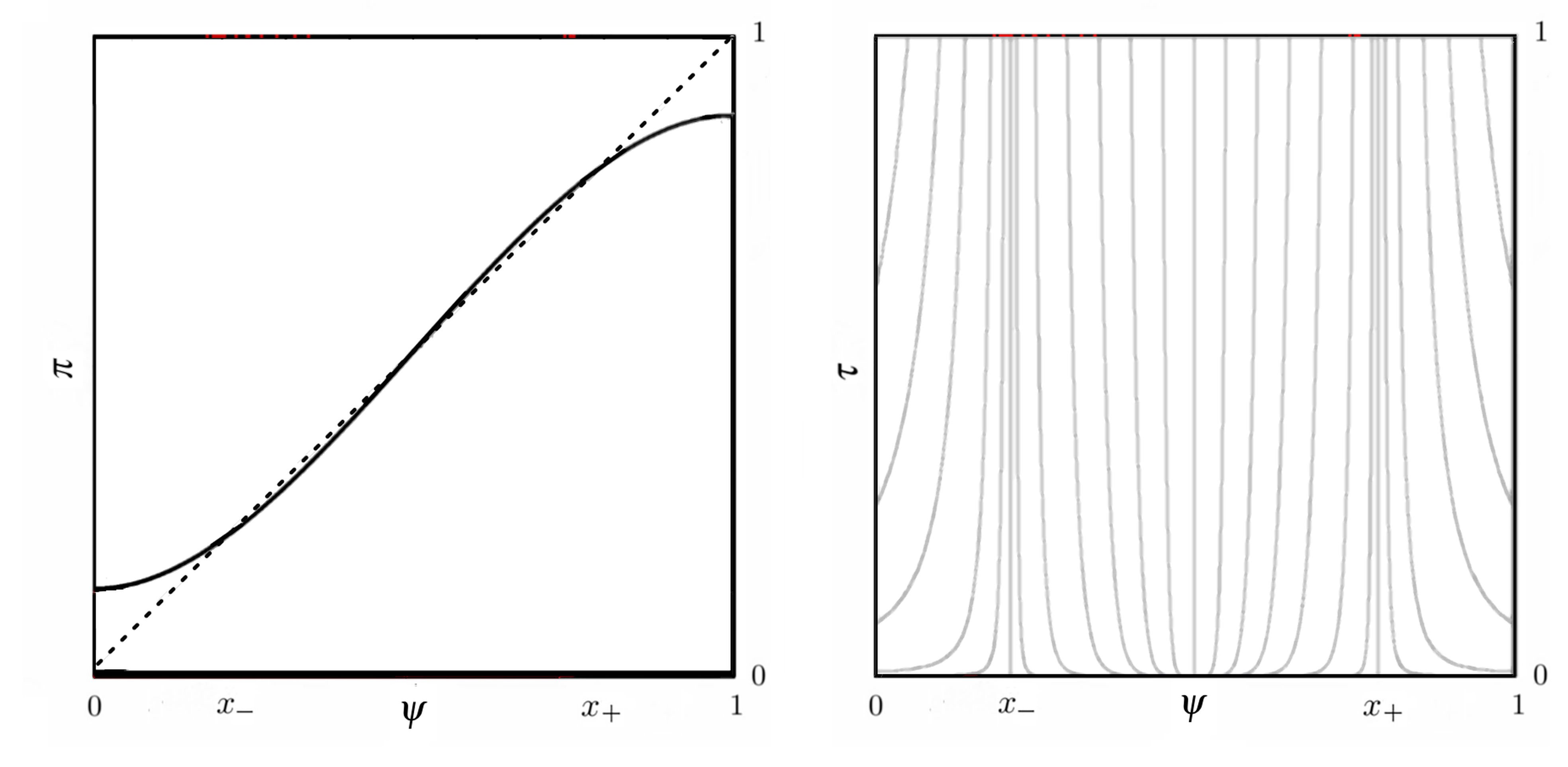}}\caption{Urn function (left) and zero--cost share trajectories (right) of the DEK model. From Franchini
\& Balzan 2023 \cite{FB}. Sample trajectories are in Figure~\ref{fig:Traj}.
The axis labels $x_-$ and $x_+$ in this figure and Figure~\ref{fig:Traj}
denote the attracting shares $\psi_-$ and $\psi_+$ in
Eq.~(\ref{eq:app-dek-fixed}).}
\label{fig:Urna}
\end{figure}
\clearpage
}

For any attracting fixed share $\psi_*$ with $\partial_\psi\pi(\psi_*)<1$,
linearizing Eq.~(\ref{eq:Caucy}) gives, for a small displacement,
\begin{equation}
\psi(\tau)-\psi_*\simeq(\psi_0-\psi_*)
(\tau_0/\tau)^{1-\partial_\psi\pi(\psi_*)}.
\label{eq:app-local-relaxation}
\end{equation}
The DEK relaxation exponent is $(5-6p)/2$ below the transition and
$6p-5$ on either stable branch above it. At $p=5/6$ the linear term
vanishes and the drift is $-4(\psi-1/2)^3/3$. A macroscopic initial
share determines the basin of its trajectory; a microscopic seed
can select either branch through early fluctuations.
Figure~\ref{fig:Traj}\label{call:fig:Traj} compares microscopic and macroscopic initial conditions.

\afterpage{%
\clearpage
\begin{figure}[p]
{\Huge\includegraphics[width=\linewidth]{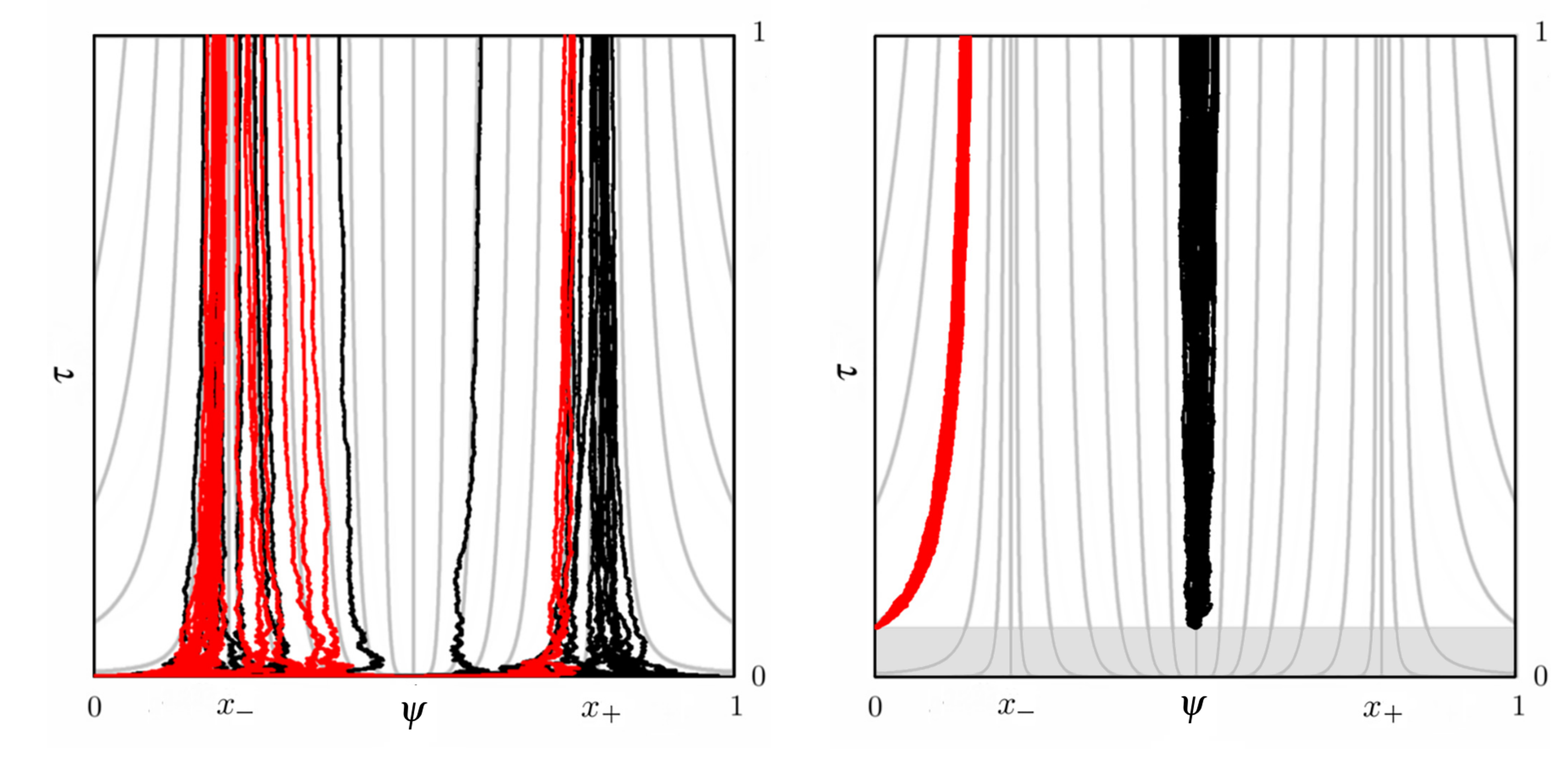}}\caption{Effect of initial conditions on the DEK model: microscopic
(left) and macroscopic (right) initial saturation.
Both panels use $k=3$, $p=7/8$ and $T=2^{14}$.
Red and black histories start from $\psi_0=0$ and $\psi_0=1/2$,
respectively. The initial saturation is $2/T$ on the left and $1/13$
on the right; gray curves are zero--cost trajectories.
A prescribed macroscopic
initial share selects a unique zero--cost trajectory in the scaling limit.
From Franchini \& Balzan 2023 \cite{FB}}
\label{fig:Traj}
\end{figure}
\clearpage
}

\begin{samepage}
The terminal
entropy and path costs in \cite{FB} distinguish these initializations.
The nonlinear function above is inserted in the action of
Eq.~(\ref{eq:scaledaction}), retaining the endpoint constraint rather
than inferring a nonconcave entropy solely from its generating function.\par
\end{samepage}

Information cascades give a related rule \cite{HisakadoNakayamaMori2024}.
With probability $p$, a voter follows the majority of a fixed odd
sample of size $k$; otherwise the first alternative is chosen
independently with probability $q$. Uniform sampling with replacement
gives
\begin{equation}
\pi(\psi)=(1-p)q+p\sum\limits_{(k+1)/2\leq j\leq k}
\binom{k}{j}\psi^j(1-\psi)^{k-j}.
\label{eq:app-independent-herders}
\end{equation}
The uniform-network version in that work samples distinct predecessors.
Its hypergeometric kernel converges to this response with error at most
$pk(k-1)/(2n)$, by the collision comparison underlying
Eq.~(\ref{eq:app-sampling-bound}). For $0\leq p<1$ and $0<q<1$,
Section~\ref{sec:asymptotic-urns} transfers the LDP and CGF from the
same compatible finite seed.

In the complete, structurally balanced network of
\cite{MohandasEtAl2025}, every pair of members is linked: links are
positive within either faction and negative between the two factions.
Each newcomer contacts a uniformly chosen member and joins that
member's faction with probability $p$, the other faction otherwise. The faction count follows exactly the affine urn
in Eq.~(\ref{eq:app-erw-urn}); the initial faction sizes specify its seed.

In experiments on social influence, the response can be measured instead
of prescribed by a sampling rule. The first alternative may be a product,
an answer to a question, or a cultural item. Its popularity is $\psi$,
and $\pi(\psi)$ is the probability that the next participant chooses it.
Van de Rijt \cite{VanR} studies situations in which an initially popular
but less appealing option loses its advantage. Positive influence alone
does not prevent this correction: even when $\pi$ increases with $\psi$,
the response may remain below the current popularity. The sign of
$\pi(\psi)-\psi$ directly describes the resulting movement of the share.
The empirical curves are specific to the alternatives and the information
shown to participants; there is no single fitted function common to all
the experiments.
Figure~\ref{fig:Vander}\label{call:fig:Vander} compares experimental trajectories under different initial--information conditions.

\afterpage{%
\clearpage
\begin{figure}[p]
{\Huge\includegraphics[width=\linewidth]{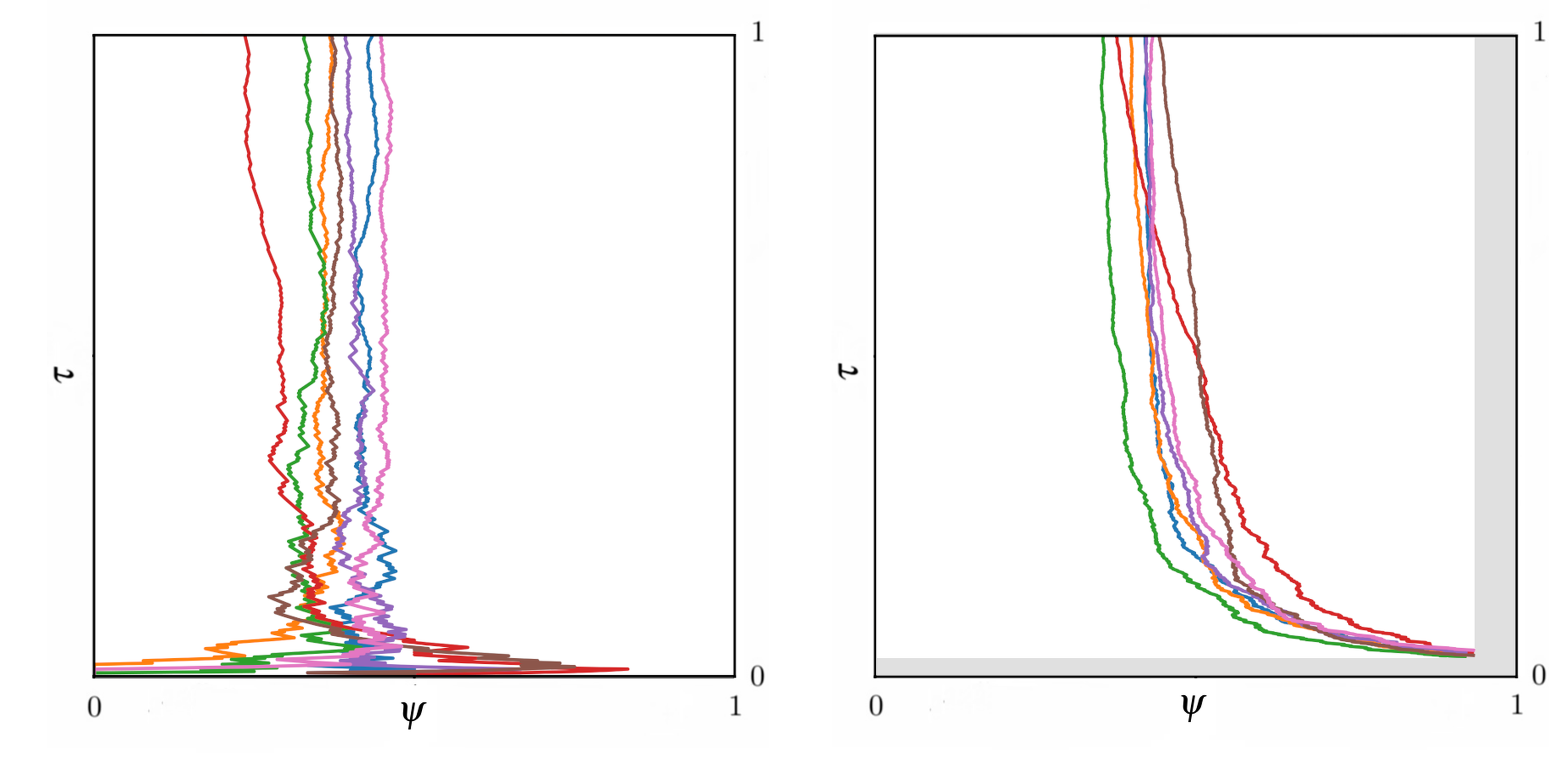}}\caption{Experimental trajectories under different initial-information
conditions, illustrating the role of the initial advantage discussed
alongside Figure~\ref{fig:Traj}. The panels show the proportion of A
answers versus urn saturation $\tau$ for the seven questions in the
van de Rijt experiments \cite{VanR}. Left: about 530 participants,
with both answer counts starting at zero and the true counts displayed.
Right: about 3500 participants, with an artificial advantage for A,
summarized as about 110 versus 10 in \cite{Gelast}. The initial counts
vary between questions. Figure provided by Alexandros Gelastopoulos
(IAS Toulouse).}
\label{fig:Vander}
\end{figure}
\clearpage
}

\begin{samepage}
Gelastopoulos and collaborators \cite{Gelast} distinguish a gradual
response to popularity from the extra effect of holding a majority.
Their influence curve can be written
\begin{equation}
\pi(\psi)=g(\psi)+(M/2)\operatorname{sgn}(\psi-1/2),
\label{eq:app-majority-effect}
\end{equation}
where $g$ is continuous at $1/2$, $M$ is the size of the jump and
$\operatorname{sgn}(0)=0$. The curve must take values in $[0,1]$.
The jump represents a change in behaviour when an option becomes the
more popular one, even by a small margin. If $d=1-2g(1/2)\geq0$ measures
the first option's disadvantage when the options are equally popular,
the condition $M>d$ is sufficient, under the assumptions of their
model, for a positive probability that the first option's limiting
share exceeds $1/2$ (lock-in).
\end{samepage}
Examples of measured response curves appear in Figure~\ref{fig:Gelast}\label{call:fig:Gelast}.

For comparison, the same work considers the continuous logistic response
\begin{equation}
\pi(\psi)=\frac{1}{1+[(1+d)/(1-d)]\exp(b(1-2\psi))},
\label{eq:app-social-logistic}
\end{equation}

\afterpage{%
\clearpage
\begin{figure}[p]
\centering{}{\Huge\includegraphics[angle=-90,scale=0.26]{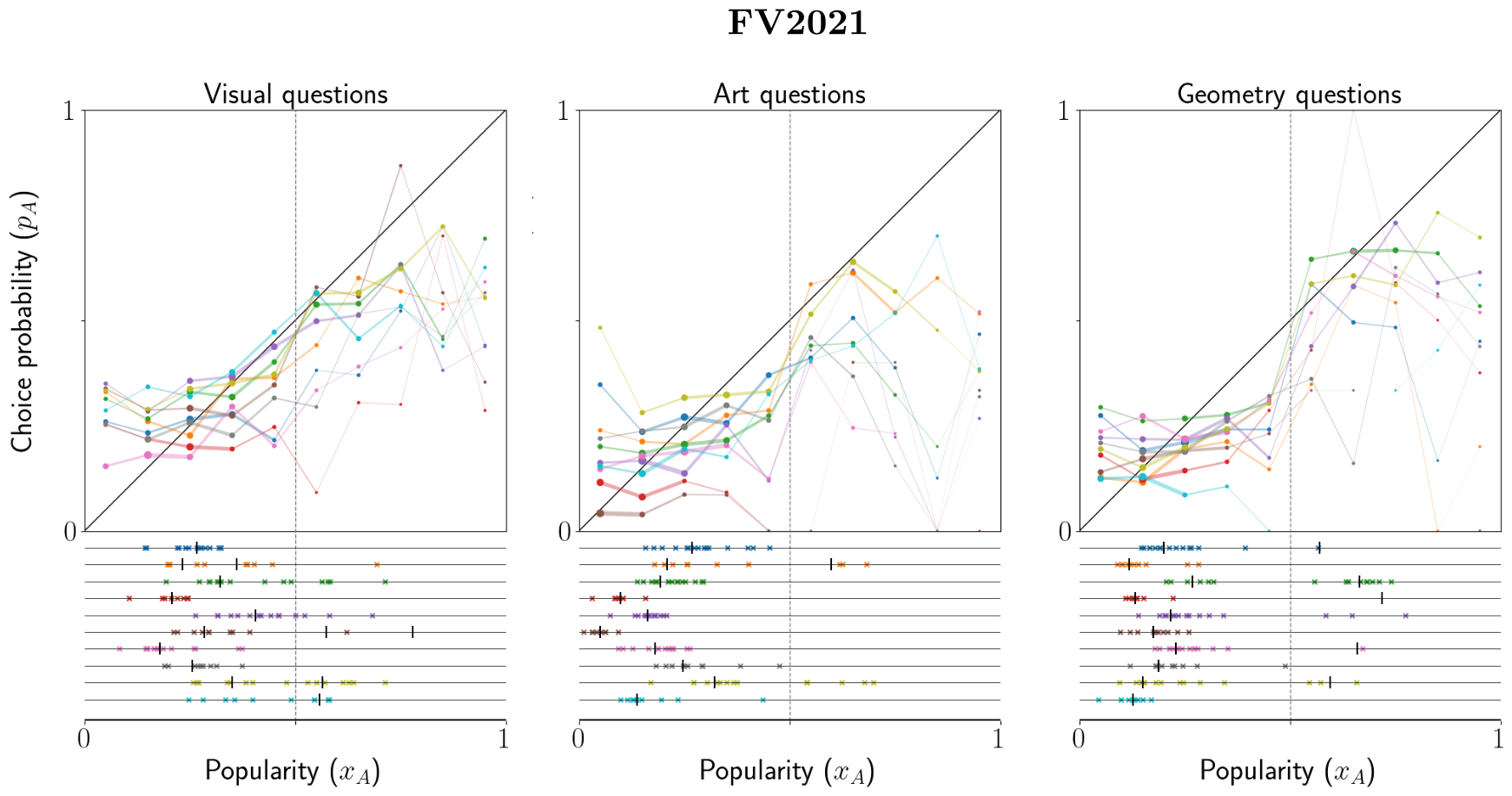}}\caption{Examples of urn functions from actual experiments, from the 2024 preprint version of \cite{Gelast}.}
\label{fig:Gelast}
\end{figure}
\clearpage
}

with $0\leq d<1$ and $b\geq0$. At a fixed share,
\begin{equation}
\partial_\psi\pi(\psi_*)=2b\psi_*(1-\psi_*).
\end{equation}
Equation~(\ref{eq:app-local-relaxation}) therefore applies when this
quantity is below one. For $d=0$, the balanced share changes stability
at $b=2$. With $b<2$, the response has a unique attracting fixed
point for every allowed $d$. To compare this model with observed
growth, Section~\ref{sec:inverse} reconstructs the response along a
smooth zero--cost scaling trajectory:
\begin{equation}
\pi(\psi(\tau))=\psi(\tau)+\tau\partial_\tau\psi(\tau).
\end{equation}
For the logistic family, the logarithm of the ratio of the reconstructed
probability to its complement (the log odds) is affine in $\psi$, with
slope $2b$ and intercept $-b-\log[(1+d)/(1-d)]$. This is a testable response shape, not a
formula fitted to every dataset. The discontinuous majority jump
does not satisfy the continuity assumption of the path LDP; its
published lock--in criterion remains a separate result.

Wealth accumulation gives another direct interpretation of reinforcement.
For two agents, suppose the next unit of capital is assigned with weights
proportional to $a_1 w^q$ and $a_2 w^q$, evaluated at their respective
wealths, where $a_1,a_2>0$ and $q>0$. With $\psi$ the first agent's wealth
share, homogeneity cancels the total wealth and gives
\begin{equation}
\pi(\psi)=\frac{a_1\psi^q}{a_1\psi^q+a_2(1-\psi)^q}.
\label{eq:app-capital}
\end{equation}
The coefficients describe unequal investment abilities; $q$ controls
the dependence of returns on existing wealth. The symmetric choice
$a_1=a_2$ gives the familiar power reinforcement function. Its linear
case is the classical P\'olya rule. General nonlinear feedback, monopoly
and the dependence on initial wealth are studied by Gottfried and
Grosskinsky \cite{Gottfried_2}; a broader account and applications appear
in Gottfried's thesis \cite{Gottfried_3}. The cancellation above uses
the common power.
For $q\ne1$, the interior fixed share and its slope are
\begin{equation}
\psi_*=[1+(a_1/a_2)^{1/(q-1)}]^{-1},\qquad \partial_\psi\pi(\psi_*)=q.
\end{equation}
With positive initial wealth for both agents, sublinear reinforcement
gives coexistence and relaxation exponent
$1-q$; superlinear reinforcement makes this share a basin boundary
between the two attracting monopolies. For $q=1$, unequal abilities
favor the more able agent; equal abilities give the neutral Polya urn.

For two colors and one addition, the normalized feedback of
\cite{RuszelThacker2024} is
\begin{equation}
\pi(\psi)=f(\psi)/[f(\psi)+f(1-\psi)].
\label{eq:app-normalized-feedback}
\end{equation}
Here $f$ is continuous and strictly increasing, with $f(0)=0$ and
$f(1)=1$; $f(u)=u^q$ recovers equal-ability power reinforcement.
Under their differentiability and reinforcement assumptions, the
authors distinguish dominance, a share tending to one, from fixation,
when every later addition has that color; their published fixation
theorem requires a stronger reinforcement bound.
For superlinear count weights, tails of the losing counts and the
time to monopoly are studied in \cite{GottfriedGrosskinsky2026}.

The model with wages \cite{Gottfried_1} divides each added unit into
a deterministic salary part and a randomly allocated capital part.
Let $r\in[0,1)$ be the salary fraction and $\gamma\in[0,1]$ the first
agent's share of salaries. The two wealth increments are
\begin{equation}
r\gamma+(1-r)\sigma_{n+1},\qquad r(1-\gamma)+(1-r)(1-\sigma_{n+1}).
\end{equation}
The wealth itself consequently does not have a binary increment.
For the two-agent power model an exact change of variables recovers
one. Writing $w_n$ for the first agent's wealth when total wealth is $n$,
the transformed urn share is
\begin{equation}
\psi_n=(w_n/n-r\gamma)/(1-r).
\label{eq:app-wage-share}
\end{equation}
Its count $n\psi_n$ increases by $\sigma_{n+1}$ at the next allocation.
The corresponding function is
\begin{equation}
\pi(\psi)=
\frac{a_1[r\gamma+(1-r)\psi]^q}
{a_1[r\gamma+(1-r)\psi]^q+
a_2[r(1-\gamma)+(1-r)(1-\psi)]^q}.
\label{eq:app-wages}
\end{equation}
This representation follows directly from the wealth update. It requires
an initial condition compatible with the transformed urn: its black and
white counts must be nonnegative integers for the binary initial
histories used here. In particular the initial wealth share must lie between
$r\gamma$ and $1-r(1-\gamma)$. Initial wealth proportional to salary
shares satisfies the nonnegativity condition. The physical wealth share
is recovered as $r\gamma+(1-r)\psi_n$.
For equal abilities and equal salary shares, $\gamma=1/2$, the
balanced state has
\begin{equation}
\pi(1/2)=1/2,\qquad \partial_\psi\pi(1/2)=q(1-r).
\end{equation}
Wages therefore make it attracting when $q(1-r)<1$, with exponent
$1-q(1-r)$ in Eq.~(\ref{eq:app-local-relaxation}). In particular,
for $q=1$ and $r>0$ the rule is affine and
Eq.~(\ref{eq:app-affine-trajectory}) gives the exact relaxation
exponent $r$. This quantifies how a salary component can oppose
the amplification of initial wealth differences.

Figure~\ref{fig:review-capital-wages}\label{call:fig:review-capital-wages} compares the two cases at the same superlinear exponent.

\afterpage{%
\clearpage
\begin{figure}[p]
\centering
\includegraphics[width=\linewidth]{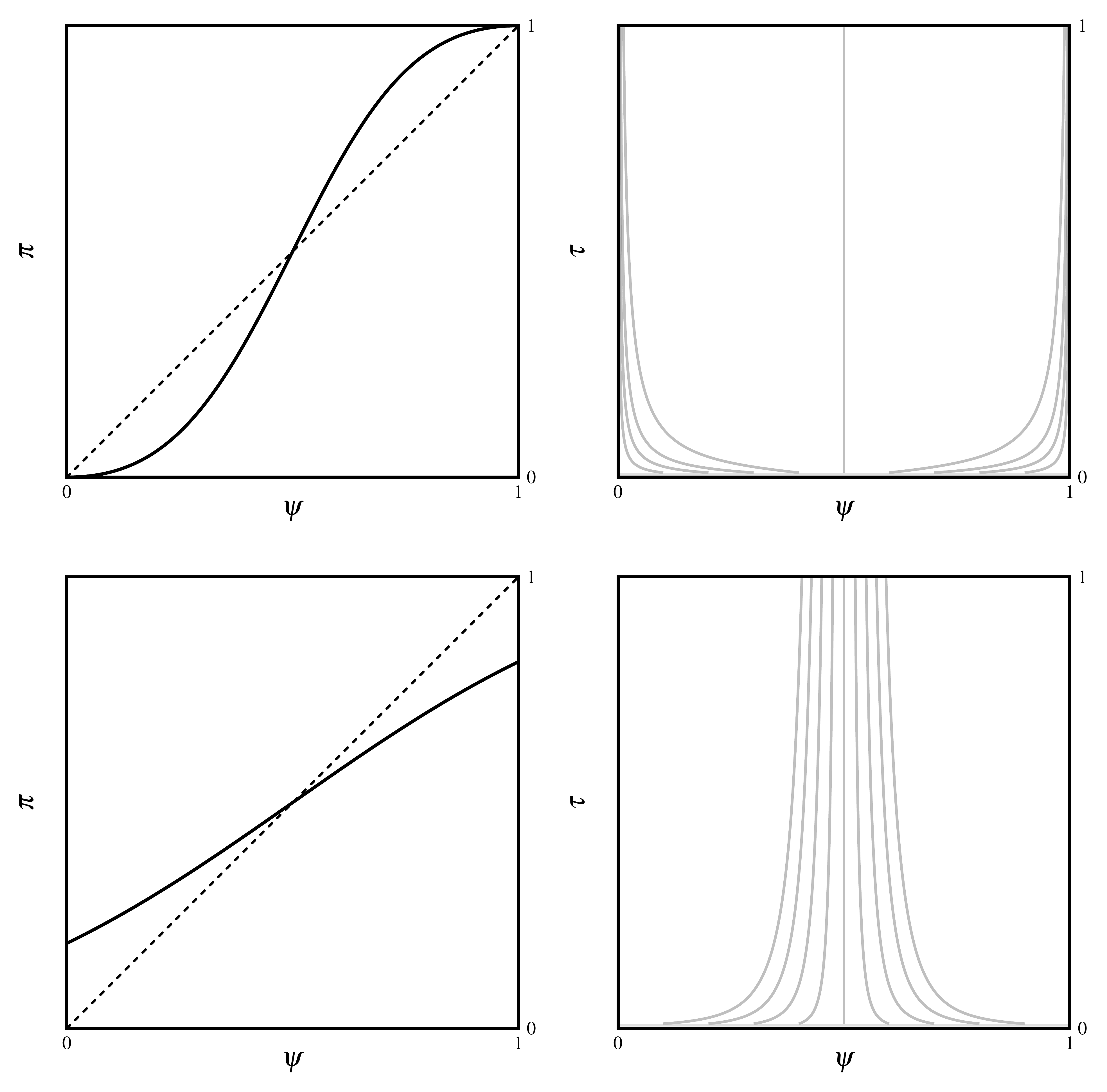}
\caption{Urn functions (left) and zero--cost share trajectories (right) for equal abilities and reinforcement exponent $q=2$. Without salaries (top), the endpoints attract; with $r=0.65$ and $\gamma=1/2$ (bottom), $1/2$ attracts. In the lower row the physical wealth share is $0.325+0.35\psi$, as in Eq.~(\ref{eq:app-wage-share}). Gray trajectories start from $\psi_0=0.1,0.2,\ldots,0.9$ at $\tau_0=0.01$. All axes span $[0,1]$, with saturation increasing upwards. The shaded strip marks $0\leq\tau<\tau_0$. The diagonal is dashed.}
\label{fig:review-capital-wages}
\end{figure}
\clearpage
}

Espinosa and Horna \cite{Espinosa} give an explicit P\'olya construction
for sequential collective decisions, with urn function
\begin{equation}
\pi(\psi)=\psi.
\end{equation}
The reinforcement scale is absorbed into the initial composition,
whose rescaled counts must be nonnegative integers. This neutral
urn has a continuum of fixed points, as discussed in
Section~\ref{sec:definition}; the law of the limiting share depends
on the initial composition.

\begin{samepage}
Iyer \cite{Iyer=000020} studies persistent leadership among growing
agents. In the two-agent balls-in-bins model with a common power
reinforcement, the probability of the next allocation is the symmetric
case of Eq.~(\ref{eq:app-capital}), with $a_1=a_2$. This count process
is an exact HLS urn.\par
\end{samepage}
More generally, let the common count weight $f(j)$ be finite and
positive on accessible counts, with a fixed finite seed, and suppose
\begin{equation}
\lim_{j\to\infty}f(\lfloor cj\rfloor)/f(j)=c^q,
\qquad c>0,\qquad q>0.
\label{eq:app-regular-variation}
\end{equation}
Then the kernel
\begin{equation}
\pi_n(\psi)=f(n\psi)/[f(n\psi)+f(n(1-\psi))]
\label{eq:app-regular-weight-kernel}
\end{equation}
converges uniformly on accessible shares to the same power response.
There is also a direct comparison of path probabilities, including
the boundaries. The weight divided by $(j+1)^q$ is slowly varying: its ratio at a
large count and at any fixed positive multiple of that count tends
to one.
The comparison power urn uses positive pseudocounts $j+1$; an accessible
zero count requires $f(0)>0$. Each count level is selected at most once
per color. Uniform convergence controls the slowly varying factors for count
fractions bounded away from zero. Near zero, Potter bounds control
ratios of these factors by a constant times arbitrarily small positive
powers of the larger shifted count divided by the smaller one. These
estimates control the sums of log weights in the numerators and,
since $q>0$, the denominators that normalize the two choices. Path log--probabilities thus differ uniformly by $o(T)$.
The fixed pseudocount shift vanishes after scaling, giving the same
path LDP and CGF as the limiting power urn. This comparison concerns the balls-in-bins allocation
rule, rather than Iyer's more general waiting-time construction.

\subsection{Irreversible magnetic growth}\label{app:magnetic}

An irreversible cluster that grows by adding one of two species gives a
physical realization of the urn rule. In the mean--field magnetic
interpretation, $\sigma_n=1$ denotes a positive spin and $\sigma_n=0$
a negative spin. The share is $\psi_n$ and the magnetization per
particle is $2\psi_n-1$. Time in the urn counts additions; the waiting
times of a continuous--time growth model are a separate observable.

Klymko, Garrahan and Whitelam study this irreversible growth process
\cite{KGW}. Normalizing the two addition rates gives the urn function
\begin{equation}
\pi(\psi)=\frac{1+\tanh(J(2\psi-1)+h)}{2},
\label{eq:magneticurn}
\end{equation}
where $J\geq0$ is the interaction strength and $h$ the external field,
with temperature absorbed into these parameters.
The social logistic response in Eq.~(\ref{eq:app-social-logistic})
is this same urn function with $J=b/2$ and
$h=-\log[(1+d)/(1-d)]/2$. For matching initial histories,
the binary history laws therefore coincide.
Both probabilities are strictly between zero and one for finite $J$ and $h$. The papers
\cite{Jack=000020LD,KGGW} study the fluctuations and rare growth
histories of this class of processes; \cite{Jack=000020LD-1,Whitelam2018} place
these results in the broader setting of dynamical large deviations.

The fixed--point equation is
\begin{equation}
2\psi-1=\tanh(J(2\psi-1)+h).
\label{eq:magneticfixed}
\end{equation}
At zero field the symmetric composition $\psi=1/2$ is fixed. For
$J<1$ it is attracting; for $J>1$ it is repelling and two attracting
compositions appear, related by interchange of the species. These are
the familiar mean--field fixed points, now reached through irreversible
accumulation. Their equation alone does not identify the distribution
of a growing cluster with a Curie--Weiss Gibbs measure.

The trajectory equation also gives the rate of approach. At a fixed
composition $\psi_*$, differentiation of Eq.~(\ref{eq:magneticurn}) gives
\begin{equation}
\partial_\psi\pi(\psi_*)=4J\psi_*(1-\psi_*).
\label{eq:magnetic_slope}
\end{equation}
For a deterministic trajectory starting near an attracting fixed point
at $\tau_0>0$, linearizing
Eq.~(\ref{eq:Caucy}) yields
\begin{equation}
\psi(\tau)-\psi_*\simeq(\psi_0-\psi_*)
(\tau_0/\tau)^{1-4J\psi_*(1-\psi_*)}.
\label{eq:magnetic_relaxation}
\end{equation}
In particular, at zero field and $J<1$ the exponent is $1-J$.
Relaxation slows down as $J$ approaches one. At $J=1$ and $h=0$, the linear
term vanishes; retaining the leading cubic term gives, near $1/2$,
\begin{equation}
\psi(\tau)-1/2\simeq
\frac{\psi_0-1/2}{[1+(8/3)(\psi_0-1/2)^2\log(\tau/\tau_0)]^{1/2}}.
\label{eq:magnetic_critical_trajectory}
\end{equation}
The decay is then logarithmic in the ratio of cluster sizes. These
are local approximations to the deterministic scaling trajectories;
the full trajectories are obtained from the transformed urn function
$\Pi$. For the three--sample DEK urn at $p=5/6$ in
Section~\ref{app:economics}, the cubic drift is exact, so the right--hand
side of Eq.~(\ref{eq:magnetic_critical_trajectory}) gives its exact displacement
from $1/2$.
Correlations retaining information about the first addition
are studied for the magnetic response in \cite{MoriHisakado2015}.
Finite-size scaling and the stochastic-differential analysis of
\cite{Kazuaki} describe the associated critical slowing of memory.
For deterministic trajectories in the regime with two attractors,
the initial composition selects a basin. A prescribed macroscopic seed and a fixed finite
seed impose different constraints on the admissible histories, as
illustrated by the rare-growth calculations in
\cite{Jack=000020LD,KGGW}.
Whitelam and Jacobson \cite{WhitelamJacobson2021} examine branch
selection and transient histories in this same growth model.
Their finite-time estimates distinguish observed steady growth
from the transients contributing to the full asymptotic rate.

Substitution in the Lagrangian of Section~\ref{sec:lft} gives
\begin{equation}
\mathcal L(\sigma_{n+1},\psi_n)
=(2\sigma_{n+1}-1)[J(2\psi_n-1)+h]
-\log[2\cosh(J(2\psi_n-1)+h)].
\label{eq:magnetic_lagrangian}
\end{equation}
Consequently, for a fixed finite initial cluster, the scaled action is
\begin{equation}
\Phi(\varphi)=\int_0^1d\tau\,
\{(2\partial_\tau\varphi-1)[J(2\psi-1)+h]
-\log[2\cosh(J(2\psi-1)+h)]\}.
\label{eq:magnetic_action}
\end{equation}
Here all time-dependent quantities are evaluated at $\tau$. For example, the probability
that the final share of positive spins exceeds a prescribed threshold
is obtained from Eq.~(\ref{eq:VARPRINZ}) by restricting $Q(E^*)$ to
that terminal share. The maximizing action difference
$\Phi-\Phi_0^*$ gives the logarithmic probability and the corresponding
most probable growth histories; its negative is the nonnegative cost.
For $h=0$ and $1/2<J<1$, Jack and Harris
\cite{Jack=000020Harris} analyze rare final compositions generated
by an initial excursion followed by relaxation. This mechanism gives
a bound on the rate near the balanced composition, rather than an
exact expression for the whole rate function.
At $J=0$ the response is constant,
$a=(1+\tanh h)/2$, and the entropy reduces to the explicit Bernoulli
formula in Eq.~(\ref{eq:app-bernoulli-entropy}).
A field chosen once for each growing sample is considered in
Section~\ref{app:quenched-magnetic}.

\subsection{Random-walk range, trapping and polymers}\label{app:polymers}

The connection with random walks and polymers illustrates a different use
of the urn description. The quantity represented by the black balls is
an accumulated geometric observable. Its distribution determines the
number of configurations at a given energy, and hence the thermodynamic
quantities of several polymer models. The same distribution also enters
the survival probability of a walker among random traps
\cite{FranchiniBalzanRANGE2018,Franchini=000020Range,FranchiniPhD2015,Huges}.

Consider a simple random walk on the cubic lattice $\mathbb Z^d$, with
each of the $2d$ nearest-neighbor steps chosen with equal probability.
Write its history as $\omega_n=(x_0,\ldots,x_n)$, with $x_0=0$, and
use the uniform walk measure for walk probabilities and expectations below.
We count a black increment whenever the next step reaches a site
already visited:
\begin{equation}
\sigma_n=\mathbb I(x_n\in\{x_0,\ldots,x_{n-1}\}).
\label{eq:range_increment}
\end{equation}
The variables $\varphi_n$ and $\psi_n$ retain their definitions: they
count these increments, divided by $T$ and $n$, respectively. If $R_n$
is the number of distinct visited sites, including the origin, then
\begin{equation}
R_n=n+1-n\psi_n=n+1-T\varphi_n.
\label{eq:range_share}
\end{equation}
Each return to an old site contributes one increment, irrespective of
how often that site was visited before. This observable is therefore
different from a count of all pairs of intersecting monomers. The
first step visits a new site, so $\sigma_1=0$.

For a given geometric history, the probability of a black increment is
the fraction of neighbors of the endpoint that have already been
visited. Averaging this fraction over histories with the same urn share
defines the microcanonical atmosphere
\cite{FranchiniBalzanRANGE2018}:
\begin{equation}
\pi_n(\psi)=\mathbb P(\sigma_{n+1}=1\mid\psi_n=\psi).
\label{eq:range_atmosphere}
\end{equation}
Equivalently,
\begin{equation}
\pi_n(\psi)=\mathbb E[
\#\{e:|e|_1=1,\ x_n+e\in\{x_0,\ldots,x_n\}\}/(2d)
\mid\psi_n=\psi].
\label{eq:range_atmosphere_geometry}
\end{equation}
This is an exact conditional average, defined at attainable values of
$\psi$. It is a function of both the share and the walk length; it does
not require an affine fit.

The relation with urns follows directly by separating the two ways
to reach a given return count. At step $n+1$, a count $j$ is obtained
either from $j-1$ by a return or from $j$ by visiting a new site:
\begin{multline}
\mathbb P((n+1)\psi_{n+1}=j)=\pi_n((j-1)/n)
\mathbb P(n\psi_n=j-1)\\
+[1-\pi_n(j/n)]\mathbb P(n\psi_n=j).
\label{eq:range_master}
\end{multline}
Terms with an unattainable preceding count are zero. This is the
master equation of an urn with time-dependent function $\pi_n$.
Consequently, that urn and the random walk have exactly the same
distribution of the return count at every length. The scalar urn
compresses the geometry through the conditional average: two walks
with the same count can have different endpoint neighborhoods. Thus
the identity concerns the count distributions, while geometric
observables retain their own conditional distributions.

Under the uniform walk measure, every $n$-step history has probability
$(2d)^{-n}$. The urn therefore also determines the density of states:
\begin{equation}
\#\{\omega_n:n\psi_n=j\}=(2d)^n\mathbb P(n\psi_n=j).
\label{eq:range_density_states}
\end{equation}
The atmosphere provides a practical way to study this density: it
averages a local continuation probability over a constrained ensemble,
rather than counting all configurations separately. This is the
microcanonical approach developed in
\cite{FranchiniBalzanRANGE2018,Franchini=000020Range}.

Figure~\ref{fig:range-atmospheres}\label{call:fig:range-atmospheres}
compares numerical estimates of the atmosphere with the affine
mean--field approximation developed in \cite{FranchiniBalzanRANGE2018}.
The agreement at low return fractions in higher dimensions connects
this geometric problem with the explicitly solvable affine urns of
Section~\ref{app:classical}. At a critical return fraction, denoted
$m_c$ in \cite{Franchini=000020Range}, the conditioned walks undergo
a coil--globule transition: extended chains with the scaling of
self--avoiding walks give way to compact clusters. The affine approximation neglects
intersections between different chain segments; these become important
in compact configurations \cite[Section~V]{FranchiniBalzanRANGE2018}.
This change in geometry accounts for the departure from the affine
branch at high return fractions. The analysis and simulations in
\cite{Franchini=000020Range} predict the critical fraction at the
random--walk value for $d=3,4$ and above it for $d\geq5$.
Thus the vertical dotted guides locate the predicted transition only
in the first two panels. At the finite lengths shown, the change
spreads over a band of return fractions.

\afterpage{%
\clearpage
\begin{figure}[p]
\centering
\makebox[\linewidth][l]{\hspace*{-8.894462pt}\includegraphics[width=356.654812pt]{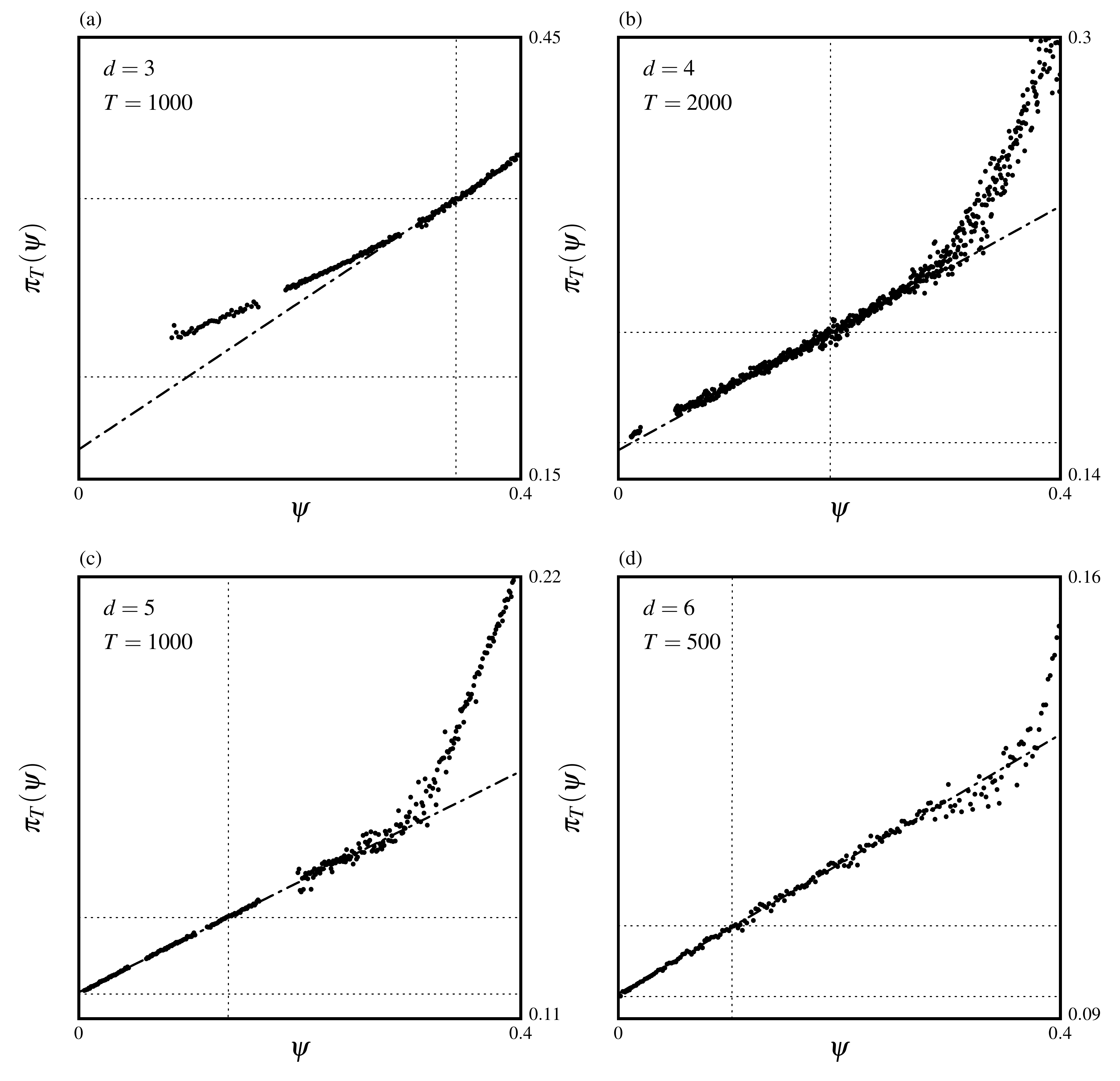}}
\caption{Finite--length estimates of the return atmosphere $\pi_T(\psi)$
for simple random walks on $\mathbb Z^d$, adapted from Figs.~1--2 of
\cite{FranchiniBalzanRANGE2018}. Panels (a--d) show $d=3,4,5,6$, with
chain lengths $T=1000,2000,1000,500$, respectively, as reported in the
source. Points were obtained with the pruned--enriched Rosenbluth
method. Dash--dotted lines show the affine mean--field approximation
of \cite[Eqs.~(28)--(29)]{FranchiniBalzanRANGE2018}, without fitting
the plotted data. The vertical and upper horizontal dotted lines
intersect at the exact limiting random--walk point, where the return
probability equals the return fraction. The lower horizontal dotted
line gives the limiting self--avoiding value at zero return fraction, obtained
from literature numerical estimates of the connective constant through
Eq.~(\ref{eq:saw_connective}). The high--density upturn in panels
(b--d) is consistent with the coil--globule crossover studied in
\cite{Franchini=000020Range}. The critical fraction is predicted at
the random--walk value for $d=3,4$ and above it for $d\geq5$;
the vertical guide therefore does not mark the transition in (c,d).
Panel (a) also shows substantial deviations from the affine
approximation within the extended regime.}
\label{fig:range-atmospheres}
\end{figure}
\clearpage
}

Passing to an autonomous HLS function is a further asymptotic step.
The approximation result in \cite{Franchini_URNS} requires uniform
convergence on $[0,1]$, the stated regularity, and urn functions
$\pi_n$ and $\pi$ strictly between zero and one. Convergence on an
interior density interval alone does not control the boundary regimes
selected by every exponential tilt. The numerical atmosphere in
\cite{FranchiniBalzanRANGE2018} supports a limiting description in
the dimensions and density ranges studied; extending that description
to the full range requires the corresponding convergence and boundary
control. This asymptotic question is separate from the exact
finite-length identity \eqref{eq:range_master}.

In this subsection, $\zeta$ and $\phi$ are defined by
Eqs.~(\ref{eq:cgf}) and~(\ref{eq:terminal-entropy}) under the uniform
walk measure. The matching count distributions give the
same functions for the nonhomogeneous urn. Their identification with
the quantities of an autonomous urn requires the additional step
described above.

The Stanley model gives a direct thermodynamic interpretation of this
construction. Its energy is proportional to the number of visits to
previously occupied sites, or equivalently to the deficit of the range
\cite{FranchiniBalzanRANGE2018,Franchini=000020Range,FranchiniMS2011}.
Let $\kappa$ be the dimensionless coupling multiplying this count.
The partition function at length $T$, summed over all nearest-neighbor
histories $\omega_T$ starting at the origin, is
\begin{equation}
Z_T(\kappa)=\sum\limits_{\omega_T}\exp(-\kappa T\psi_T).
\label{eq:stanley_partition}
\end{equation}
A positive $\kappa$ penalizes returns and favors extended
configurations; a negative value favors repeated visits. Grouping the
walks by their return count gives
\begin{equation}
Z_T(\kappa)=(2d)^T\sum\limits_{0\leq j\leq T}
\mathbb P(T\psi_T=j)\exp(-\kappa j).
\label{eq:stanley_density_partition}
\end{equation}
Hence the atmosphere determines the finite-length partition function.
With the convention of Eq.~(\ref{eq:cgf}), its limiting value is
\begin{equation}
\lim_{T\to\infty}(1/T)\log Z_T(\kappa)
=\log(2d)+\zeta(-\kappa).
\label{eq:stanley_cgf}
\end{equation}
The argument is $-\kappa$ because the physical weight penalizes a
positive energy, whereas the generating function uses
$\exp(\beta T\psi_T)$.

When the passage from $\pi_n$ to $\pi$ is justified on the full share
range, including the relevant boundary regimes, and the sample--path
assumptions of Section~\ref{sec:ldp} hold, the
variational principle gives the same
quantity directly in terms of the histories of the associated urn:
\begin{equation}
\lim_{T\to\infty}(1/T)\log Z_T(\kappa)=\log(2d)+
\sup_{\varphi\in Q}\{\Phi(\varphi)-\Phi_0^*(\varphi)-\kappa\varphi(1)\}.
\label{eq:stanley_variational}
\end{equation}
The terminal energy competes with the probability cost of the urn
history. The equality of the terminal distributions transfers this
free energy to the walk; it does not identify the joint law of its
successive returns with that of the urn. At a coupling where $\zeta$ is differentiable,
the limiting mean return and range fractions in this Gibbs ensemble
are, respectively,
\begin{equation}
\partial_\beta\zeta(\beta),\qquad 1-\partial_\beta\zeta(\beta),\qquad \beta=-\kappa.
\label{eq:stanley_densities}
\end{equation}
Thus the urn generating function determines the thermodynamic response
to the coupling, not just the correspondence of observables. The
entropy density $\phi$ gives the exponential abundance at a prescribed
return density, after removing the common factor $(2d)^T$.
Spatial shape can then be studied within the corresponding
microcanonical ensembles \cite{Franchini=000020Range,FranchiniBalzanRANGE2018}.

This use of HLS concerns a Gibbs ensemble of complete walks at a fixed
length. It does not require generating those walks by a local urn
rule. In particular, the conditional probability of the next step in
the Gibbs ensemble can depend on the possible future continuations.
The terminal weight in \eqref{eq:stanley_density_partition} already
contains everything needed for the partition function, even when a
scalar description of that conditioned spatial dynamics is unavailable.

The Rosenstock trapping problem uses the same range distribution for
a different physical question \cite{Huges,FranchiniMS2011,FranchiniPhD2015}.
Suppose that each lattice site independently contains an immobile trap
with probability $\rho$, where $0<\rho<1$. A walker survives a given history if every
visited site is free. The survival probability averages over both
the traps and the walk. Averaging first over the traps leaves the
following expectation over walk histories:
\begin{equation}
\mathbb P(\text{survival to }T)
=\mathbb E[(1-\rho)^{R_T}].
\label{eq:rosenstock_range}
\end{equation}
Here the origin is included among the sites that must be free. Using
\eqref{eq:range_share}, the same probability becomes
\begin{equation}
\mathbb P(\text{survival to }T)
=(1-\rho)^{T+1}\mathbb E[
\exp(-\log(1-\rho)T\psi_T)].
\label{eq:rosenstock_urn}
\end{equation}
If the origin is known to be free initially, the prefactor is instead
$(1-\rho)^T$. These are exact identities for independent traps. They
show why the trapping problem belongs to the urn application: changing
the trap density selects different exponential moments of the same
return count. The survival probability weights compact histories more
strongly, since fewer distinct sites need to be free. The relevant
information is the distribution of the range, including its tails,
rather than only its typical value. For $0<\rho<1$, the generating
function gives
\begin{equation}
\lim_{T\to\infty}(1/T)\log\mathbb P(\text{survival to }T)
=\log(1-\rho)+\zeta(-\log(1-\rho)).
\label{eq:rosenstock_entropy}
\end{equation}
For the simple walk on the infinite lattice this limit is zero:
annealed survival decays more slowly than exponentially in $T$
\cite{Huges}. Consequently, $\zeta(\beta)=\beta$ for $\beta>0$, and the attractive
Stanley branch in Eq.~(\ref{eq:stanley_cgf}) equals $\log(2d)-\kappa$
for $\kappa<0$. The subexponential survival cost and the spatial geometry
are not resolved by this leading free energy.

Self-avoiding walks occupy the other extreme of the same description.
A history is self-avoiding precisely when it has no black increments.
If $c_T$ denotes the number of $T$-step self-avoiding walks starting
at the origin, then
\begin{equation}
c_T=(2d)^T\mathbb P(T\psi_T=0).
\label{eq:saw_count}
\end{equation}
Setting the return count to zero in \eqref{eq:range_master} yields
\begin{equation}
c_{n+1}=2d[1-\pi_n(0)]c_n.
\label{eq:saw_atmosphere_count}
\end{equation}
Thus the atmosphere at zero directly measures the average fraction
of forbidden extensions in the uniform ensemble of self-avoiding
walks. Whenever $\pi_n(0)$ has a limit, its relation with the
connective constant $\mu_d$ is
\begin{equation}
\pi(0)=1-\mu_d/(2d),\qquad
\mu_d=\lim_{T\to\infty}c_T^{1/T}.
\label{eq:saw_connective}
\end{equation}
Here $\pi(0)$ denotes the limit of $\pi_n(0)$ alone; no global
autonomous approximation is needed for this boundary relation.
It is part of the atmosphere approach
\cite{FranchiniBalzanRANGE2018}. It connects an urn probability with a
central enumeration problem in polymer theory. In the notation of
Section~\ref{sec:thermodynamic}, Eq.~(\ref{eq:saw_count}) gives
\begin{equation}
\phi(0)=\log[\mu_d/(2d)].
\label{eq:saw_entropy}
\end{equation}
If the boundary atmosphere converges, this is also
$\phi(0)=\log[1-\pi(0)]$. The zero-return endpoint is attainable
here because the first increment is white. Thus the boundary entropy
is the logarithmic probability of self-avoidance per step, and determines the
connective constant through $\mu_d=2d\exp\phi(0)$.
At finite $T$, the limit $\kappa\to+\infty$ of the normalized Stanley measure selects
the uniform self-avoiding ensemble, because only histories with zero
returns retain a nonzero relative weight.

Attractive contacts require one additional geometric observable. In
the model of \cite{Franchini=000020Range=000020Line}, the range is
studied together with the number of nearest-neighbor links between
visited sites. Denote that number here by $\ell(\omega_T)$, counting
each pair of distinct neighboring sites once. On a self-avoiding
walk, $T$ of these links are the consecutive bonds of the chain;
the remaining $\ell(\omega_T)-T$ are nonconsecutive contacts.
Introducing a dimensionless attractive coupling $\theta>0$ gives
\begin{equation}
Z_T(\kappa,\theta)=\sum\limits_{\omega_T}
\exp[-\kappa T\psi_T+\theta\ell(\omega_T)].
\label{eq:isaw_joint_partition}
\end{equation}
Conditioning on the return count separates the two contributions:
\begin{equation}
Z_T(\kappa,\theta)=(2d)^T\sum\limits_{0\leq j\leq T}
\mathbb P(T\psi_T=j)\exp(-\kappa j)
\times\mathbb E[\exp(\theta\ell(\omega_T))\mid T\psi_T=j].
\label{eq:isaw_conditional_partition}
\end{equation}
Terms with zero probability are omitted. The first factor is supplied
by the range atmosphere. The second
contains the contact statistics within each constrained ensemble.
Both are needed: an old-site count does not specify how many
neighboring contacts a spatial configuration contains.

Taking $\kappa\to+\infty$ at fixed length restricts the model to
self-avoiding walks with attractive contacts, giving the interacting
self-avoiding walk. Including the consecutive links only multiplies
its partition function by $\exp(\theta T)$. In three dimensions, the competition between self-avoidance and
contact attraction produces the coil-to-globule transition, from an
extended chain to a compact configuration; the collapse threshold
is called the $\theta$ point. The microcanonical study in
\cite{Franchini=000020Range=000020Line} varies both the range and
the link count to investigate this transition. This extends the
use of the atmosphere description: the range component remains
binary, while the joint problem retains the additional contact
information. It is not necessary to force both observables into a
single binary increment.

\subsection{Biological growth and chemical reinforcement}
\label{app:biology}

Biological growth provides a direct interpretation of the two colors:
each addition can be a unit of length, a new cell, or a reaction
product. The index $n$ counts these additions, rather than elapsed
physical time. A binary HLS description applies when the probability of each type
depends on the current share alone. The following examples separate
published urn models from simple effective rules whose biological
implementation still has to be established.

\subsubsection{Competition between neurites}

Khanin and Khanin \cite{Khanin} model the early development of neuron
polarity by allowing competing neurites, the projections growing from
a neuron, to grow one unit at a time.
The probability of selecting a neurite is proportional to a power
of its current length. For two neurites, let $\sigma_n=1$ record an
addition to the first, and let $\psi_n$ be its fraction of the total
length. With lengths measured in growth units, the initial lengths
specify the initial urn history. The exact urn function is
\begin{equation}
\pi(\psi)=\psi^{q}/[\psi^{q}+(1-\psi)^{q}],
\qquad q>0.
\label{eq:bio-neurites}
\end{equation}
For $q=1$ this is the Polya rule. Superlinear reinforcement,
$q>1$, favors the longer neurite sufficiently strongly that,
almost surely, after a finite random time only one continues to
grow. This gives a qualitative mechanism for selecting an axon from
two initially similar neurites.

For $q\ne1$, the fixed shares are $0$, $1/2$ and $1$.
The central point attracts when $q<1$ and repels when $q>1$;
the endpoint stability is reversed. Indeed,
\begin{equation}
\partial_\psi\pi(1/2)=q.
\end{equation}
For sublinear reinforcement, a small initial length imbalance therefore
relaxes along a zero--cost trajectory as
\begin{equation}
\psi(\tau)-1/2\simeq(\psi_0-1/2)(\tau_0/\tau)^{1-q},
\qquad 0<\tau_0\leq\tau.
\end{equation}
At $q=1$ every share is fixed and the deterministic trajectory
is constant. For superlinear reinforcement, convergence of a share
to zero alone would not imply that growth stops; the stronger
monopoly statement uses the particular power rule.

The same function gives the local action density, whose integral
is the scaled log likelihood:
\begin{equation}
L(\partial_\tau\varphi,\pi(\psi))=
q\partial_\tau\varphi\log\psi+
q(1-\partial_\tau\varphi)\log(1-\psi)
-\log[\psi^q+(1-\psi)^q].
\label{eq:bio-neurite-action}
\end{equation}
Thus the entropy of exceeding a terminal length share $x\in(0,1)$
is obtained from Eq.~(\ref{eq:VARPRINZ}) with
\begin{equation}
Q(E^*)=\{\varphi\in Q:\varphi(1)\geq x\}.
\label{eq:bio-threshold-event}
\end{equation}
Both initial lengths are positive and fixed as $T$ grows. Because the urn function
attains zero and one, the boundary accessibility conditions of
\cite{Franchini_URNS} apply; the logarithmic formula is interpreted
with the corresponding endpoint conventions. As in the general
variational principle, equality requires matching large--deviation
bounds for the chosen event.

The equal-ability response without salaries in Figure~\ref{fig:review-capital-wages} also illustrates Eq.~(\ref{eq:bio-neurites}), with $q=2$.

\subsubsection{Clonal growth and developmental history}

For two cell types $A$ and $B$, suppose that the only events are
births, with per--cell birth rates $r_A(\psi)$ and $r_B(\psi)$
positive when both types are present. At population size $n$, the rates of
adding the two types are respectively $n\psi r_A(\psi)$ and
$n(1-\psi)r_B(\psi)$. Observing the process at successive births
therefore gives the exact rule
\begin{equation}
\pi(\psi)=\psi r_A(\psi)/[\psi r_A(\psi)+(1-\psi)r_B(\psi)].
\label{eq:bio-clonal}
\end{equation}
Here $\sigma_n$ records the newborn's type, $n\psi_n$ counts type
$A$, and $\varphi_n=n\psi_n/T$ records this count relative to the
chosen terminal population $T$. Waiting times between births are
additional information. This construction is part of the broader
connection between urns and evolving populations
\cite{BioBenaim2004}. It also covers two--type evolutionary games
implemented through pure births, when positive fitnesses depend
only on the current share.

Equal constant rates give the Polya urn. An interior fixed share
requires $r_A(\psi)=r_B(\psi)$ and attracts if their difference
crosses from positive to negative there. The zero--cost share
trajectory is
\begin{equation}
\tau\partial_{\tau}\psi=
\psi(1-\psi)[r_A(\psi)-r_B(\psi)]/
[\psi r_A(\psi)+(1-\psi)r_B(\psi)].
\end{equation}
For constant positive unequal rates, integration gives the implicit
trajectory from $0<\psi_0<1$ at $\tau_0>0$:
\begin{equation}
\log[\psi(\tau)/\psi_0]
-(r_A/r_B)\log\{[1-\psi(\tau)]/(1-\psi_0)\}
=(r_A/r_B-1)\log(\tau/\tau_0).
\label{eq:bio-clonal-trajectory}
\end{equation}
The share increases when $r_A>r_B$ and decreases when $r_A<r_B$;
setting $\tau=1$ gives the terminal composition for the specified
initial population.

Figure~\ref{fig:review-clonal-growth}\label{call:fig:review-clonal-growth} illustrates how a per--cell birth advantage changes the composition of a growing population.

\afterpage{%
\clearpage
\begin{figure}[p]
\centering
\includegraphics[width=\linewidth]{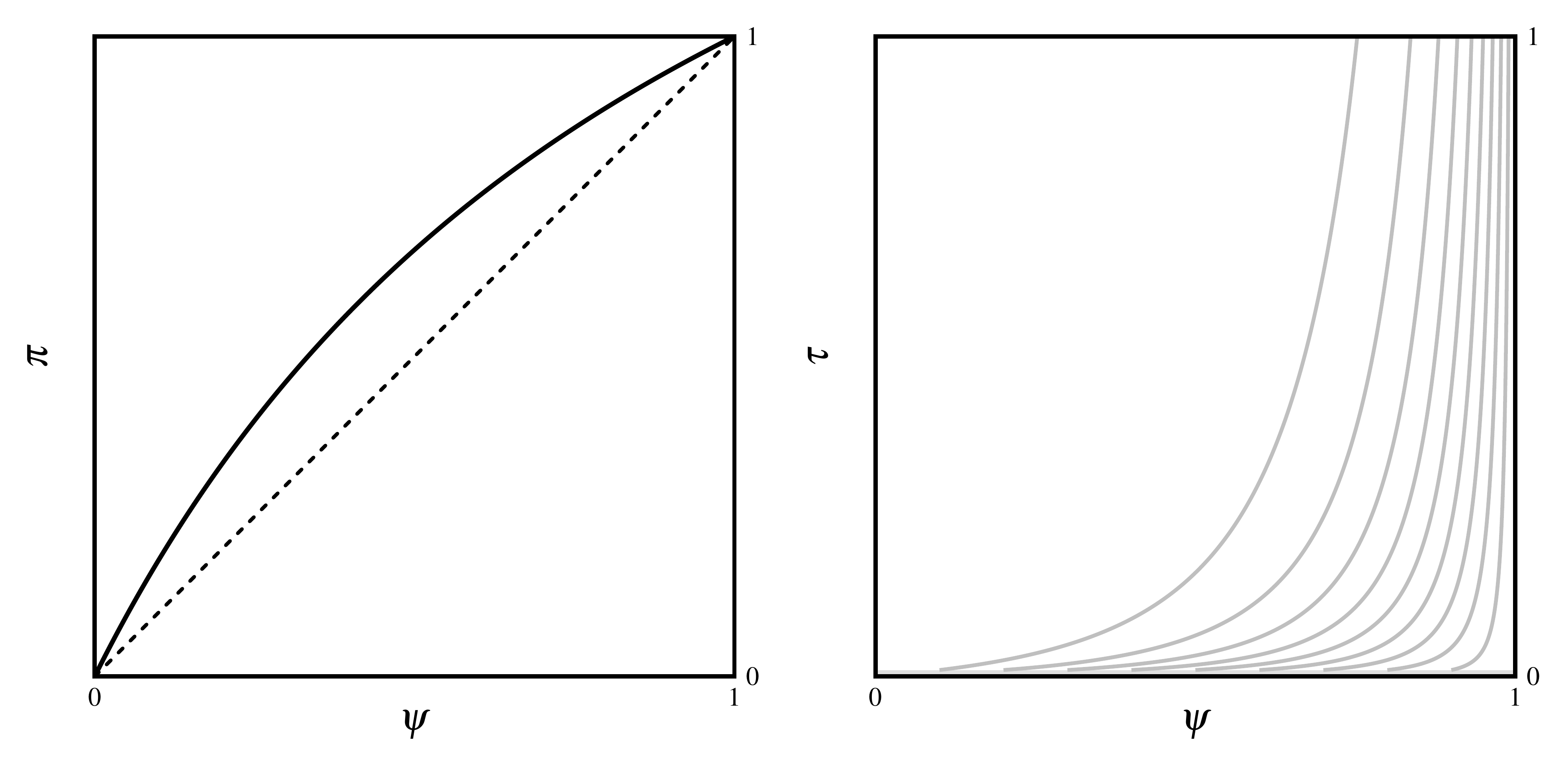}
\caption{Urn function (left) and zero--cost share trajectories (right) for two--clone growth with $r_A/r_B=2$. The response $\pi(\psi)=2\psi/(1+\psi)$ drives every interior share towards clone $A$ dominance. Gray trajectories start from $\psi_0=0.1,0.2,\ldots,0.9$ at $\tau_0=0.01$ and satisfy Eq.~(\ref{eq:bio-clonal-trajectory}). All axes span $[0,1]$. Saturation increases upwards and counts births relative to the terminal population, rather than physical time. The shaded strip marks $0\leq\tau<\tau_0$. The diagonal is dashed.}
\label{fig:review-clonal-growth}
\end{figure}
\clearpage
}

The birth-only growth rule in Kreger et al. \cite{BioKreger2024},
restricted to two clones with constant fitnesses, gives
Eq.~(\ref{eq:bio-clonal}). A shared developmental history specifies
the initial composition; subsequent independent growth is described
by two urns starting from that inherited composition.
Ma et al. \cite{MaRothschild2024} use the equal-rate P\'olya case
to interpret diversity and composition differences between confined
bacterial communities. This published urn model reproduces features
of experiments and mechanical simulations; the correspondence concerns
the statistical growth model, rather than an exact reduction of the
underlying cell mechanics.

\subsubsection{Control of cell--type proportions}

The opposite feedback is useful when growth should preserve a
prescribed mixture. An illustrative controller with target share $p$ is
\begin{equation}
\pi(\psi)=\{1+[(1-p)/p]\exp[\kappa(\psi-p)]\}^{-1},
\qquad 0<p<1,\qquad \kappa\geq0.
\label{eq:bio-controller}
\end{equation}
This is a proposed effective growth rule. It gives $\pi(p)=p$ and a
unique attracting fixed point.
At $\kappa=0$, cell--type assignments (commitments) are independent with probability $p$;
positive $\kappa$ suppresses the overrepresented type. The slope is
\begin{equation}
\partial_\psi\pi(p)=-\kappa p(1-p).
\end{equation}
A small displacement along a deterministic share trajectory then satisfies
\begin{equation}
\psi(\tau)-p\simeq(\psi_0-p)
(\tau_0/\tau)^{1+\kappa p(1-p)},\qquad 0<\tau_0\leq\tau.
\end{equation}
Here $n$ counts irreversible commitments; the relaxation concerns
further additions, rather than conversion of existing cells.
For an initial history of length fixed as $T$ grows, an error of at least
$0<\epsilon<\min(p,1-p)$ in the final ratio is the event
\begin{equation}
Q(E^*)=\{\varphi\in Q:|\varphi(1)-p|\geq\epsilon\}.
\label{eq:bio-control-event}
\end{equation}
Its exponential cost is $-\phi(E^*)$, where
\begin{equation}
\phi(E^*)=\sup_{\varphi\in Q(E^*)}
\{\Phi(\varphi)-\Phi_0^*(\varphi)\}.
\label{eq:bio-control-cost}
\end{equation}
In this optimization $\Phi$ uses Eq.~(\ref{eq:bio-controller});
the entropy contribution remains $-\Phi_0^*$.
The rule is smooth and strictly between zero and one, so the
regular sample--path theory applies. The variational equality has
the event regularity condition of Eq.~(\ref{eq:VARPRINZ}).
At $\kappa=0$ it reduces to the explicit Bernoulli cost
\begin{equation}
\phi(E^*)=\max_{x\in\{p-\epsilon,p+\epsilon\}}
\{-x\log(x/p)-(1-x)\log[(1-x)/(1-p)]\}.
\end{equation}
Section~\ref{app:quenched-biological} allows the target and feedback
strength to vary between replicas while remaining fixed within each one.

\subsubsection{Asymmetric autocatalysis}

An autocatalytic product can promote production of its own type,
amplifying an initial imbalance between enantiomers, the two
non-superimposable mirror-image forms of a molecule. A direct urn
precedent is the work of Barab\'as et al. \cite{BioBarabas2009}, which
uses Polya simulations to analyze experimental distributions of
enantiomeric excess, the excess of one form over the other relative
to their total amount.

To illustrate what a nonlinear urn can represent, let $\sigma_n=1$
record a newly formed product of one handedness and $\psi_n$ its
share. A phenomenological family combines spontaneous production
$a$, cooperative order $q$, and bias $h$:
\begin{equation}
\pi(\psi)=
\exp(h)(a+\psi^q)/\{\exp(h)(a+\psi^q)+\exp(-h)[a+(1-\psi)^q]\}.
\label{eq:bio-chiral}
\end{equation}
Here $a>0$, $q>1$ and $h$ is real; positive $h$ favors the counted enantiomer.
The formula is exact for the illustrative growth process with the
two stated composition--dependent propensities.

At zero bias the racemic share $1/2$ is fixed, with
\begin{equation}
\partial_\psi\pi(1/2)=q/(1+2^q a).
\end{equation}
It attracts above the threshold
\begin{equation}
a=(q-1)/2^q.
\end{equation}
and repels below it. In the attracting regime the local trajectory is
\begin{equation}
\psi(\tau)-1/2\simeq(\psi_0-1/2)
(\tau_0/\tau)^{1-q/(1+2^q a)}.
\end{equation}
Below threshold, a small displacement grows towards an ordered
branch on the same side of $1/2$. Fixed points that cross the
diagonal upwards separate attraction basins; when further branches
coexist, their positions are found from $\pi(\psi)=\psi$.
Since $a>0$, neither pure composition is fixed.

The action density includes the normalization of the reaction
propensities:
\begin{multline}
L(\partial_\tau\varphi,\pi(\psi))=
h(2\partial_\tau\varphi-1)+
\partial_\tau\varphi\log(a+\psi^q)
+(1-\partial_\tau\varphi)\log[a+(1-\psi)^q]\\
-\log\{\exp(h)(a+\psi^q)+\exp(-h)[a+(1-\psi)^q]\}.
\label{eq:bio-chiral-action}
\end{multline}
Integrating this expression gives $\Phi$ for the variational cost
of a prescribed terminal handedness, for example the threshold
event in Eq.~(\ref{eq:bio-threshold-event}), with $\Phi_0^*$ unchanged.
Here $a>0$ makes the rule strictly between zero and one, so the
regular sample--path theory applies for a fixed initial history.

The effect of a bias fixed within each reactor and varying between
reactors is developed in Section~\ref{app:quenched-chiral}.

Figure~\ref{fig:chiral-multistability}\label{call:fig:chiral-multistability} displays the five--fixed--point regime. A stable balanced composition coexists with two stable asymmetric compositions; the two intervening repellers delimit three distinct basins.

\afterpage{%
\clearpage
\begin{figure}[p]
\centering
\includegraphics[width=\linewidth]{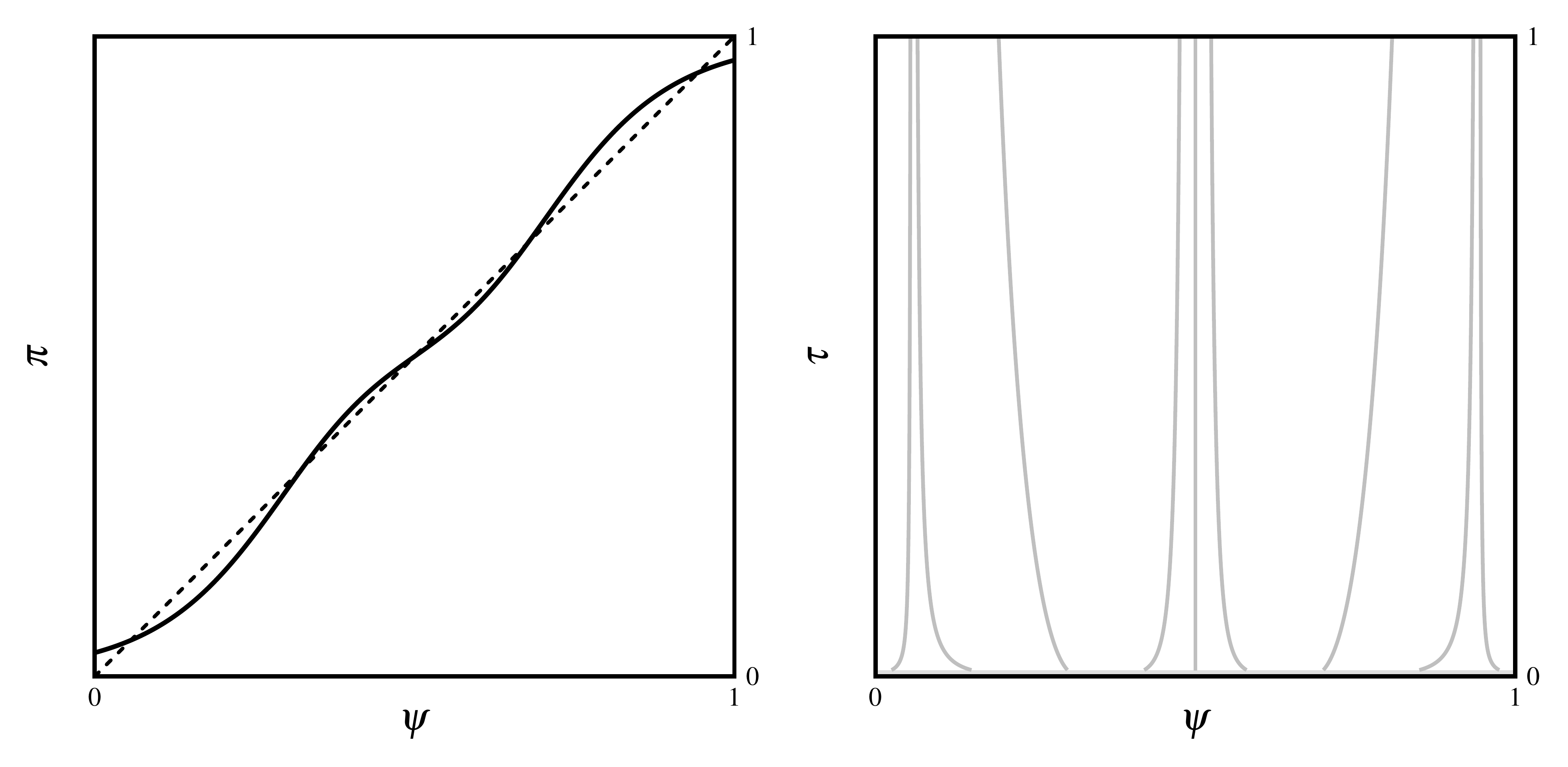}
\caption{Urn function (left) and zero--cost share trajectories (right) for autocatalytic growth with $q=8$, $a=1/25$ and $h=0$. The five intersections with the dashed diagonal give three attractors separated by two repellers. Gray trajectories start from $\psi_0=0.025,0.15,0.30,0.42,0.50,0.58,0.70,0.85,0.975$ at $\tau_0=0.01$ and remain within their respective basins. All axes span $[0,1]$, with saturation increasing upwards. The shaded strip marks $0\leq\tau<\tau_0$. The trajectory endpoints at $\tau=1$ need not coincide with the limiting compositions.}
\label{fig:chiral-multistability}
\end{figure}
\clearpage
}

\subsection{Sequential allocation}\label{app:allocation}

Sequential experiments provide a direct use of the urn as an allocation
rule; Zhang \cite{Zhang2024Urns} reviews its connections with stochastic
approximation and adaptive designs. Set $\sigma_n=1$ when participant $n$ receives treatment $A$ and
$\sigma_n=0$ for treatment $B$. Then $\psi_n$ is the fraction assigned to
$A$, and $\pi(\psi_n)$ is the probability of assigning the next participant
to it. Each assignment adds one observation, so this correspondence
preserves the HLS dynamics exactly whenever the current fraction is the
only input to the rule. The same identification describes allocating
successive users between two versions of a product. Here the purpose of
feedback is usually to correct an imbalance: the less represented group
is favoured while assignments remain random.

Efron's biased coin \cite{Efron1971} is an early example. With a bias
parameter $1/2<p<1$, its urn function is
\begin{equation}
\pi(\psi)=\begin{cases}
p,&\psi<1/2,\\
1/2,&\psi=1/2,\\
1-p,&\psi>1/2.
\end{cases}
\label{eq:app-efron}
\end{equation}
The underrepresented treatment receives the next participant with
probability $p$; a tie is resolved by a fair coin. Larger $p$ makes the
correction stronger and the next assignment more predictable. The usual
formulation uses the difference between the group sizes, but only its
sign enters the rule. Since that sign is also the sign of $\psi_n-1/2$,
the next-assignment probability depends only on $\psi_n$. The urn function is discontinuous: its exact
urn representation does not place it within the continuity assumptions
of the preceding large-deviation analysis. Its allocation share still
converges almost surely to $1/2$, by the discontinuous-rule convergence
result in \cite{BaldiAntogniniZagoraiou2015}.

Wei's adaptive biased coin \cite{Wei1978} replaces the jump with a
continuous response to the relative imbalance:
\begin{equation}
\pi(\psi)=f(2\psi-1),
\label{eq:app-wei}
\end{equation}
where $f:[-1,1]\to[0,1]$ is continuous and decreasing, with
\begin{equation}
f(-u)=1-f(u),\qquad -1\leq u\leq1.
\end{equation}
This symmetry gives $f(0)=1/2$. A directly solvable member of this
family is, for $-1<b\leq0$,
\begin{equation}
f(u)=(1+bu)/2,\qquad
\pi(\psi)=(1-b)/2+b\psi.
\label{eq:app-wei-affine}
\end{equation}
Its only fixed share is $1/2$. Applying Eq.~(\ref{eq:Caucy}) gives
the complete zero--cost trajectory from $\tau_0>0$:
\begin{equation}
\psi(\tau)=1/2+(\psi_0-1/2)(\tau_0/\tau)^{1-b}.
\label{eq:app-wei-trajectory}
\end{equation}
Thus a prescribed initial imbalance is reduced by the factor
$\tau_0^{1-b}$ at capacity. Making the slope more negative
accelerates this deterministic relaxation without changing the target.
For $b=0$ and a fair first assignment, the increments are
independent fair coins, and the entropy and cumulant functions are
\begin{equation}
\phi(x)=-x\log(2x)-(1-x)\log[2(1-x)],\qquad 0\leq x\leq1,
\end{equation}
\begin{equation}
\zeta(\beta)=\log[(1+\exp(\beta))/2].
\end{equation}
The endpoint values use $0\log0=0$.

For the affine rule with an initial history of length fixed as $T$
grows, take the event that
the final allocation misses balance by at least $\epsilon\in(0,1/2)$:
\begin{equation}
Q(E^*)=\{\varphi\in Q:|\varphi(1)-1/2|\geq\epsilon\}.
\label{eq:app-allocation-event}
\end{equation}
The function is smooth and takes values strictly between zero and one. Hence its terminal logarithmic probability is obtained
by the variational principle, with
\begin{equation}
\phi(E^*)=\sup_{\varphi\in Q(E^*)}\int_0^1d\tau\,
\{L(\partial_\tau\varphi,(1-b)/2+b\psi)
-L(\partial_\tau\varphi,\partial_\tau\varphi)\},
\label{eq:app-allocation-cost}
\end{equation}
under the event regularity condition of Eq.~(\ref{eq:VARPRINZ}).
For the independent case this optimization gives
\begin{equation}
\phi(E^*)=-(1/2+\epsilon)\log(1+2\epsilon)
-(1/2-\epsilon)\log(1-2\epsilon).
\end{equation}
This is a concrete cost of allocation imbalance, rather than an error
probability for estimating a treatment effect.

McNulty and Raban \cite[Section~2.3]{McNultyRaban2026} consider
self-balancing sampling with time-dependent bias parameters. For two
choices, take $\beta_i(n)=b_i/n$, with $b_i>0$, and positive initial
sampling weights $p_i(0)$, for $i=1,2$. After the first draw,
$\sigma_n=1$ records choice one. Each initial sampling weight is reduced
exponentially in the number of previous selections of that choice,
with rate $\beta_i(n)$. The two counts are $n\psi_n$ and $n(1-\psi_n)$,
so choosing $\beta_i(n)=b_i/n$ makes the weights depend only on the
current share. Normalizing these weights gives the autonomous logistic response
\begin{equation}
\pi(\psi)=\{1+[p_2(0)/p_1(0)]
\exp[b_1\psi-b_2(1-\psi)]\}^{-1}.
\label{eq:app-balancing-sampler}
\end{equation}
This is the logistic response of the cell--type controller in
Eq.~(\ref{eq:bio-controller}), with $\kappa=b_1+b_2$ and $p$ its
unique attracting fixed share. It decreases strictly, with
local relaxation given by Eq.~(\ref{eq:app-local-relaxation}). For the
stated parameters it is smooth and uniformly separated from zero and
one on $[0,1]$. More
generally, $n\beta_i(n)\to b_i$ as $n\to\infty$ gives uniform
convergence to this response, so Section~\ref{sec:asymptotic-urns} transfers the path LDP
and CGF under its stated initial-history conditions.

The inverse relation of
Section~\ref{sec:entropy} also connects an admissible entropy profile
to its feedback function. This concerns the asymptotic assignment
law; finite-sample inference about measured outcomes remains a
separate statistical problem.
Response-adaptive designs that update treatment estimates from observed
outcomes have additional statistical state: robustness and efficiency
for such designs are studied in \cite{YeHuMa2024,Zhang2026Designs}.

\subsection{Memory and admission feedback}\label{app:memory}

A growing memory supplies a natural interpretation of the two colors.
Let $\sigma_n=1$ mean that the $n$th stored item is erroneous, and
$\sigma_n=0$ that it is correct according to a specified criterion.
Then $\psi_n$ is the contaminated fraction. If the conditional
probability of an erroneous addition depends on the stored history
only through this fraction, the memory follows the HLS rule
\begin{equation}
\pi(\psi)=\mathbb P(\sigma_{n+1}=1\mid\psi_n=\psi).
\label{eq:memory_urn}
\end{equation}
Here $n$ counts stored items, rather than queries or attempted writes.
The reduction is exact under the stated dependence and one-item
accumulation rule. Whether these conditions describe a particular
memory system is a separate modeling question.

In the binary accumulation model, an affine admission response has
the familiar form
\begin{equation}
\pi(\psi)=a+b\psi,\qquad 0\leq a\leq1,\qquad 0\leq a+b\leq1.
\label{eq:memory_affine}
\end{equation}

For $b<1$, the fixed share is $a/(1-b)$. Starting with share
$\psi_0$ at a positive saturation $\tau_0$, the Cauchy problem
in Section~\ref{sec:trajectories} gives the exact scaling trajectory
\begin{equation}
\psi(\tau)=a/(1-b)+[\psi_0-a/(1-b)](\tau_0/\tau)^{1-b},\qquad
0<\tau_0\leq\tau\leq1.
\label{eq:memory_affine_trajectory}
\end{equation}
The terminal contamination is therefore
\begin{equation}
\psi(1)=a/(1-b)+[\psi_0-a/(1-b)]\tau_0^{1-b}.
\label{eq:memory_affine_terminal}
\end{equation}
The deviation from equilibrium decays as a power of the ratio of
memory sizes. When $b$ approaches one from below, this power becomes
small: an initially contaminated bank relaxes slowly even though the
fixed point is attracting. These formulas concern the deterministic
scaling trajectory, not each random history. The contaminated share
also does not determine prediction accuracy without knowing which
items are stored and used.

For $b=0$ and $0<a<1$, admissions are independent Bernoulli trials.
From a fixed initial memory, the entropy density of
Section~\ref{sec:entropy} becomes
\begin{equation}
\phi(x)=x\log(a/x)+(1-x)\log[(1-a)/(1-x)],\qquad 0<x<1,
\label{eq:memory_bernoulli_entropy}
\end{equation}
For a terminal contamination threshold $0<x<1$, the exponential
cost is zero when $x\leq a$; when $x>a$,
its logarithmic probability per stored item is $\phi(x)<0$.

A simple filter modifies this admission rule. Suppose that,
at share $\psi$, a candidate is erroneous with probability $p(\psi)$.
A filter accepts an erroneous candidate with probability $\alpha$
and a correct candidate with probability $t$, with $0<\alpha,t\leq1$.
Conditioning on an actual admission gives
\begin{equation}
\pi(\psi)=\frac{\alpha p(\psi)}{\alpha p(\psi)+t[1-p(\psi)]}.
\label{eq:memory_filter}
\end{equation}
This is the exact urn function for the specified growing-memory
model.
Its derivative, when $p$ is differentiable, is
\begin{equation}
\partial_\psi\pi(\psi)=\alpha t\,\partial_\psi p(\psi)/[t+(\alpha-t)p(\psi)]^2.
\label{eq:memory_filter_derivative}
\end{equation}
For affine candidate contamination $p(\psi)=a+b\psi$, with the
probability constraints in Eq.~(\ref{eq:memory_affine}), the filtered
fixed shares solve
\begin{equation}
(\alpha-t)b\psi^2+[t+(\alpha-t)a-\alpha b]\psi-\alpha a=0.
\label{eq:memory_filter_fixed}
\end{equation}
Only roots in the accessible share interval represent equilibria.
An interior root with $\partial_\psi\pi<1$ is attracting, while one with
$\partial_\psi\pi>1$ is repelling; equality needs the nonlinear drift.
Thus the filter changes both the equilibrium contamination and
its local relaxation. In particular, $\alpha<t$ lowers the erroneous
admission probability wherever $0<p<1$. The equal-acceptance case
$\alpha=t$ leaves the response unchanged. Inserting the filtered
response into Eq.~(\ref{eq:Caucy}) gives its trajectories and terminal
shares.

To quantify a high-contamination outcome, start from a fixed finite
memory containing both kinds of items and choose $0<x<1$. The event
and its scaling limit are
\begin{equation}
E=\{\psi_T\geq x\},\qquad E^*=\{\psi(1)\geq x\}.
\label{eq:memory_threshold}
\end{equation}
If $p$ is H\"older-continuous and strictly between zero and one,
Section~\ref{sec:ldp} gives the following specialization of the path
variational problem:
\begin{equation}
\phi(E^*)=\sup_{\varphi\in Q(E^*)}\int_0^1d\tau\,
[L(\partial_\tau\varphi,\alpha p(\psi)/[t+(\alpha-t)p(\psi)])
-L(\partial_\tau\varphi,\partial_\tau\varphi)].
\label{eq:memory_threshold_entropy}
\end{equation}
Here all time-dependent quantities are evaluated at $\tau$, with $\psi=\varphi/\tau$
for $\tau>0$. This variational problem determines the most probable histories leading to
the threshold as well as their exponential cost. A negative value
means exponential suppression. Zero means no suppression on this
scale, not a positive limiting probability.

The discrete action in Section~\ref{sec:lft} is equally explicit:
each admitted item contributes
\begin{multline}
\mathcal{L}(\sigma_{n+1},\psi_n)=\sigma_{n+1}\log[\alpha p(\psi_n)]\\
+(1-\sigma_{n+1})\log[t(1-p(\psi_n))]
-\log[t+(\alpha-t)p(\psi_n)].
\label{eq:memory_filter_lagrangian}
\end{multline}
The final term normalizes the accepted candidates. Summing these
contributions over a history gives its log weight, conditional on
the initial memory. Thus the same admission mechanism determines
the trajectory, the terminal risk and the field action; applying
these conclusions to a particular memory requires the probability that
the next admitted item is erroneous to depend on the stored history only
through the current contaminated share.

Figure~\ref{fig:memory_feedback}\label{call:fig:memory_feedback} shows how selective admission changes both the contamination response and the equilibrium approached by the growing memory.

\subsection{Machine learning and language models}\label{app:ml}

The next symbol in an HLS history is generated from the empirical
frequency of the preceding symbols. This gives a direct connection
with sequential prediction: a binary token can denote a word in a
two-word vocabulary, membership in a specified class of words, or an
event extracted from a longer text. The reduction depends on what the
bit records. We first consider predictive laws and document models,
then feedback from algorithmic decisions, and finally the appearance
of new words. These examples connect the urn to different parts of a
learning system: its data generator, its predictive distribution, or
the data accumulated through its decisions.

\afterpage{%
\clearpage
\begin{figure}[p]
\centering
\includegraphics[width=\linewidth]{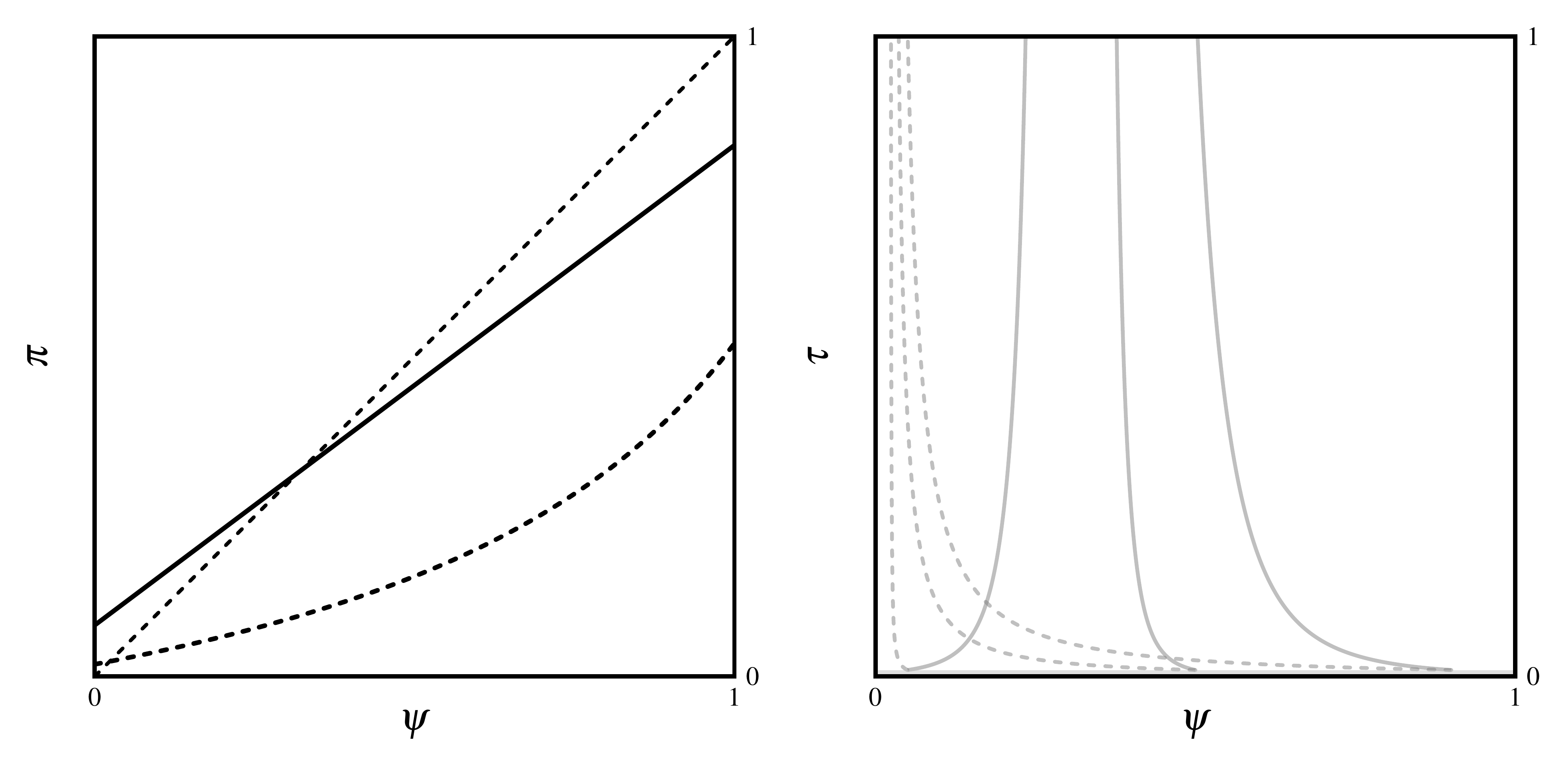}
\caption{Urn functions (left) and zero--cost share trajectories (right) for the memory admission rule. The candidate response $p(\psi)=0.08+0.75\psi$ is shown without filtering (solid, $\alpha=t=1$) and with selective admission (dashed, $\alpha=0.2$, $t=0.9$), as in Eq.~(\ref{eq:memory_filter}). The filter lowers the attracting contamination share. Gray trajectories retain these line styles and start from $\psi_0=0.05,0.5,0.9$ at $\tau_0=0.01$. All axes span $[0,1]$, with saturation increasing upwards. The shaded strip marks $0\leq\tau<\tau_0$. The diagonal is dashed.}
\label{fig:memory_feedback}
\end{figure}
\clearpage
}

\subsubsection{Sequential prediction and in-context learning}
\label{app:ml-prediction}

For a binary vocabulary, Eq.~(\ref{eq:definition}) is an
autoregressive generation rule. Its complete dependence on the prefix
is through $\psi_n$: two prefixes with the same empirical frequency
give the same next-token probability. Using the function $L$ from
Eq.~(\ref{eq:scaler}), the log probability of a continuation of an
assigned prefix of length $m$ is
\begin{equation}
\log\mathbb P(\sigma_{m+1},\ldots,\sigma_T\mid\sigma_1,\ldots,\sigma_m)
=\sum\limits_{m\leq n<T}L(\sigma_{n+1},\pi(\psi_n)).
\label{eq:ml-sequence-loss}
\end{equation}
Its negative is precisely the cumulative binary logarithmic loss.
The urn action therefore measures the predictive score of a whole
sequence. Dependence of the next prediction on a count does not imply
exchangeability: a nonlinear response generally assigns different
probabilities to different orderings of the same symbols, through
the intermediate frequencies in this sum.

An exact and widely used example is the Beta--Bernoulli generator.
A propensity is sampled once from a Beta distribution and, conditional
on it, tokens are independent Bernoulli draws. Such generators provide
controlled tasks for meta-learning and prompt-conditioned prediction
\cite{Mikulik2020,Genewein2025,Wenliang2026}. To use the clock of
Section~\ref{sec:definition}, represent the two positive integer
prior parameters by the black and white counts in an initial history
of length $m$. At time $n$, the posterior parameters are then
$n\psi_n$ and $n(1-\psi_n)$. Integrating the propensity out gives
\begin{equation}
\mathbb P(\sigma_{n+1}=1\mid\sigma_1,\ldots,\sigma_n)=\psi_n,
\qquad \pi(\psi)=\psi.
\label{eq:ml-polya-predictive}
\end{equation}
\begingroup\clubpenalty=10000
Thus the marginal token sequence is exactly the P\'olya urn of
Section~\ref{app:classical}; $n$ includes the prior counts and the
observed tokens. Arbitrary positive real prior parameters give the
same urn with initial masses, as in the document models below.
A prompt supplies further observations and updates these two counts.
It changes the subsequent predictive distribution without changing
the rule used to process the next observation.\par\endgroup

This mechanism also appears in the solvable topic model of
Cagnetta and Christensen \cite{CagnettaChristensen2026}: for a new
topic with a uniform Dirichlet prior and a binary vocabulary, the
in-context predictor adds one to each observed count. It is
Eq.~(\ref{eq:ml-polya-predictive}) with one initial ball of each color.
The same Beta--Bernoulli reference is used to study the beliefs and
uncertainty of Bayes-filtered transformers
\cite{AswadiMaWei2026,FortiniEtAl2026}. P\'olya sequences also serve
as a controlled prediction task in the study of next-token prediction
and regret minimization, where regret is the excess cumulative
prediction loss relative to the best predictor in a specified comparison
class \cite{MohriEtAl2026}. In these applications
the exact identification concerns the generator or the ideal
predictor; how closely a trained network reproduces it is a separate
question addressed by the learning experiments.

The theory describes the sequences generated by this predictor as
well as its one-step probabilities. With a fixed finite positive
seed, the limiting share has the corresponding Beta distribution;
its random value is the slope of the limiting count trajectory.
A fixed finite prompt changes that distribution, although the
limiting CGF remains
\begin{equation}
\zeta(\beta)=\max\{0,\beta\}.
\label{eq:ml-polya-cgf}
\end{equation}
The terminal entropy remains zero in the interior, as in
Eq.~(\ref{eq:app-polya-entropy}). The CGF follows because
neighborhoods of both endpoint shares have subexponential probabilities,
meaning that their logarithms divided by $T$ tend to zero, while the
share lies between zero and one.
Positive real prior masses give the same exponential limits: the
Beta--binomial formula shows that the log ratio of any string
probability to its value for unit prior masses is $O(\log T)$,
uniformly in the count. A prompt whose length is proportional to $T$
instead imposes a macroscopic initial condition; its subsequent
zero-cost share trajectory is constant, as explained in
Section~\ref{app:classical}.

For a fixed finite seed, the vanishing terminal cost does not make all
paths equally likely.
For $x=\varphi(1)$, the P\'olya response yields
\begin{equation}
\Phi(\varphi)=L(x,x).
\label{eq:ml-polya-action}
\end{equation}
Indeed, its integrand is the derivative of
$\tau L(\psi(\tau),\psi(\tau))$, whose boundary term at zero
vanishes. The path cost is consequently
\begin{equation}
\Phi_0^*(\varphi)-\Phi(\varphi)
=\int_0^1d\tau\,L(\partial_\tau\varphi,\partial_\tau\varphi)-L(x,x).
\label{eq:ml-polya-path-cost}
\end{equation}
Convexity makes it nonnegative, with equality precisely for
$\varphi(\tau)=x\tau$. The same predictive law thus permits all
limiting token proportions but penalizes departures from the straight
count trajectory at a prescribed proportion. Finally, integrating
the latent propensity in the token probabilities, as above, differs
from averaging the conditional CGF of a propensity kept fixed.
The latter is the quenched operation of Section~\ref{sec:quenched}.

\subsubsection{Reinforced document models and linguistic inference}
\label{app:ml-documents}

Repeated words provide a direct linguistic use of the predictive
mechanism in Section~\ref{app:ml-prediction}. A term that is rare across
a collection may occur many times in an individual document. The
Dirichlet compound multinomial model accounts for this word burstiness
by drawing a vector of word probabilities once per document and then
drawing tokens independently conditional on that vector
\cite{MadsenKauchakElkan2005}. After the vector is integrated out,
the joint token law is that of a multicolor P\'olya urn. This gives
a generative interpretation of the observed dependence between
repeated occurrences.

For the binary representation, fix a word or a group of words and let
$\sigma_n=1$ when the token belongs to that group. Dirichlet parameters
add under grouping: the initial mass of each color is the sum of the
parameters of the words assigned to it. Each occurrence increases its
group's mass by one, independently of which word within the group was
chosen. Thus the grouped process closes exactly. With integer masses
represented by the assigned initial history, its response is
Eq.~(\ref{eq:ml-polya-predictive}); real masses have the interpretation
given there. The beta--binomial count law and the P\'olya trajectories
therefore apply to word frequency within a document. The model assigns
equal probabilities to token sequences with the same counts, so its
dependence describes repetition without encoding word order.

Cummins, Paik and Lv \cite{CumminsPaikLv2015} use this reinforced
document model in information retrieval. Their Smoothed P\'olya Urn
Document model (SPUD) combines document information with estimates
from the collection. In its Dirichlet smoothing construction, the
initial parameter vectors are combined into a single urn, which
retains the grouping property above. Query likelihood is evaluated
using a nonreinforced multinomial model: the HLS representation is
the document generator underlying the retrieval score.

Different words can also have different propensities to repeat.
Cummins \cite{Cummins2017Burstiness} models this by a positive diagonal
replacement matrix, whose entries specify the reinforcement added
after each word is drawn. A binary vocabulary has an exact HLS
representation; the same holds for two word classes when reinforcement
is constant within each class. Let these two reinforcement amounts be
$a_1,a_2>0$, and choose initial class weights $a_1m\psi_m$ and
$a_2m(1-\psi_m)$ for an assigned binary seed of length $m$. At time
$n$, the class weights are then $a_1n\psi_n$ and $a_2n(1-\psi_n)$,
because $n\psi_n$ counts class-one tokens including the seed. Sampling
in proportion to weight gives
\begin{equation}
\pi(\psi)=\frac{a_1\psi}{a_1\psi+a_2(1-\psi)}.
\label{eq:ml-document-reinforcement}
\end{equation}
The binary clock counts token occurrences, while the physical urn mass
increases by the reinforcement of the selected class. Equation~(\ref{eq:ml-document-reinforcement})
is the $q=1$ case of the capital-allocation rule
(\ref{eq:app-capital}). Its zero--cost trajectories are given by
Eq.~(\ref{eq:bio-clonal-trajectory}), with the two birth rates replaced
by $a_1,a_2$. Equal reinforcement gives the neutral P\'olya rule;
with both classes initially present, unequal reinforcement drives
the share towards the more strongly reinforced class. This reduction requires equal reinforcement within
each class, since otherwise its total weight also depends on the
unrecorded frequencies of its individual words.

The ordinary P\'olya count law also enters statistical explanations
of term weighting. Ahmed et al. \cite{AhmedEtAl2026} use a
beta--binomial model for the number of occurrences of a target word
against all other words in a document. This is precisely the grouped
P\'olya count law. Their derivation of terms related to TF--IDF
(term frequency--inverse document frequency), which combines a word's
frequency within a document with its rarity across documents, then
uses an approximation and a penalty in a likelihood-ratio score. The
exact urn identification concerns the underlying generative count
model; the subsequent score is a statistic constructed from it.

A different linguistic observable is a judgment of similarity between
two word usages. Aida et al. \cite{AidaEtAl2026} model such judgments
on four score levels using a categorical distribution with a uniform
Dirichlet prior, with one initial unit in each category. After the
probabilities are integrated out, the judgments have a four-color
P\'olya law. Selecting $j$ score levels as color one gives the binary
seed $m=4$, $m\psi_m=j$, for $j=1,2,3$, and again the response
Eq.~(\ref{eq:ml-polya-predictive}). Their semantic-change inference
compares a common distribution of judgments with separate
distributions for blocks of within-period and between-period
comparisons. Each block admits the same urn representation; the full
test uses the four-category likelihood. Here $\sigma_n$ records
membership of the selected group of similarity scores, and the data
determine the evidence for a change in word meaning.

The same aggregation is useful when the reinforced sequence is hidden
and must be inferred from measurements. Brady et al.
\cite{BradyEtAl2026} evaluate a differentiable particle filter on a
model with eight latent regimes whose selection probabilities are
proportional to their previous counts plus one. Membership of a fixed
group of $j$ regimes is therefore P\'olya with seed $m=8$,
$m\psi_m=j$, for $j=1,\ldots,7$. The selected regime controls the
continuous dynamics and the noisy observations. The binary HLS
process describes the latent regime history, while the learned
filter estimates the hidden state from the observations. Before
conditioning on those observations, its grouped occupations have
the beta--binomial law, with the
limiting spread described in Section~\ref{app:ml-prediction}.

\subsubsection{Recommendation and algorithmic feedback}
\label{app:ml-feedback}

Recommendation and data collection can reinforce the observations on
which later decisions depend. Several published models reduce to the
weighted power response already encountered in Eq.~(\ref{eq:app-capital}):
\begin{equation}
\pi(\psi)=\frac{a_1\psi^q}{a_1\psi^q+a_2(1-\psi)^q},
\qquad a_1,a_2>0,\qquad q>0.
\label{eq:ml-feedback-family}
\end{equation}
What changes between applications is the quantity counted and the
meaning of the coefficients.

In the recommendation model of Bellina et al.\ \cite{Bellina2023},
restricting to one user and two items makes the next click probability
proportional to a power of the item's accumulated rating. Each click
increases that rating by one. Taking positive integer initial ratings as
initial counts gives Eq.~(\ref{eq:ml-feedback-family}) with $a_1=a_2$;
$q$ is the popularity bias. Thus Eq.~(\ref{eq:bio-neurites}) describes
the same feedback in a different setting. Sublinear reinforcement
balances the two click shares, whereas superlinear reinforcement
eventually confines clicks to one item. The full collaborative model
also uses similarities between users and does not close on one share.

An observation clock is essential in the discovery-only predictive
policing model of Ensign et al.\ \cite{Ensign2018}. The next patrol
visits the first area with probability $\psi$, its share of recorded
incidents. Let $a_1,a_2\in(0,1]$ be the probabilities of discovering an
incident during a visit to the respective areas. A patrol adds one
record of the first type with probability $a_1\psi$, one of the
second type with probability $a_2(1-\psi)$, and otherwise adds nothing.
Unsuccessful visits leave the record composition unchanged.
Conditioning on the next addition therefore gives
Eq.~(\ref{eq:ml-feedback-family}) with $q=1$: $n$ counts records,
rather than patrols. This is also Eq.~(\ref{eq:memory_filter}) with
$p(\psi)=\psi$, $\alpha=a_1$ and $t=a_2$.

Reinforced bandit models describe repeated choices with random rewards
\cite{YangHuLiBai2024,YangYanJiang2026}. The same calculation applies
to the following case of the latter: two options, called arms, one selection per round,
and independent Bernoulli rewards with fixed positive success
probabilities $a_1,a_2$. Arms are selected in proportion to
their accumulated rewards. With positive initial reward counts,
observing successive successes gives the $q=1$ urn above. Here $\psi$
is the first arm's share of accumulated successes, including the
initial counts. It is not its share of all selections; the success
clock does not directly describe regret per round, the expected
reward shortfall relative to always selecting the arm with the highest
success probability.

Laufer's model of repeated risk assessments \cite{Laufer2020} instead
reinforces each sampled binary decision by one unit. Its basic
initialization has one unit of each type, and its response is the
P\'olya case $q=1$, $a_1=a_2$. A finite positive seed then produces a
random limiting share, even though every share is a fixed point of
the urn function.

Yu et al.\ \cite[Section~6]{Yu2026Capture} model repeated citations
of two competing sources by sampling in proportion to their current
masses and adding a positive, source-dependent constant after each
citation. Denote these reinforcement amounts by $a_1,a_2$ and divide
each mass by its own reinforcement amount. The resulting counts
increase by one at each citation, and their share obeys
Eq.~(\ref{eq:ml-feedback-family}) with $q=1$. Integer rescaled initial
masses give an initial HLS history; positive noninteger masses are
fixed pseudocounts. The urn concerns this stochastic citation model;
the deterministic recurrence used in their computational illustration
is a separate model.

The fixed points and stability of these responses are those already
discussed for Eq.~(\ref{eq:app-capital}). In the unequal-coefficient
case with $q=1$, Eq.~(\ref{eq:bio-clonal-trajectory}) gives the
zero--cost trajectories after setting $r_A/r_B=a_1/a_2$. If $a_1>a_2$,
the first share tends to one from a finite positive seed, although
additions of the other type continue indefinitely. This distinguishes
dominance under unequal linear reinforcement from the finite-time
fixation of the superlinear power rule. Equation~(\ref{eq:VARPRINZ})
applies with the boundary conventions for responses attaining zero
and one, and with the relevant count used as the clock.

For the whole family, a fixed finite seed containing both types gives
a particularly simple CGF. Let its size be $m$, with
$0<\psi_m<1$, independently of $T$. The probability that every later
addition is of the first type is
\begin{equation}
\prod\limits_{0\leq j<T-m}
\{1+(a_2/a_1)[m(1-\psi_m)/(m\psi_m+j)]^q\}^{-1}.
\label{eq:ml-feedback-endpoint}
\end{equation}
Here $j$ counts first--type additions after the initial history, so
$m\psi_m+j$ is the current count of the first type along this history.
The logarithm of this probability is $-o(T)$, since the logarithm of each factor tends
to zero. The same holds for a history containing only additions of
the second type. Their terminal shares tend respectively to one and
zero. These two histories give lower bounds on the exponential
moment that match the upper bound from $0\leq\psi_T\leq1$, yielding
\begin{equation}
\zeta(\beta)=\max\{0,\beta\},\qquad \beta\in\mathbb R.
\label{eq:ml-feedback-cgf}
\end{equation}
Thus the CGF at speed $T$ does not distinguish these reinforcement
strengths, despite their different typical behavior. This argument
uses a fixed finite seed; it does not describe conditioning on a
macroscopic initial history at $\tau_0>0$.

Adaptive memory gives a related application of
Section~\ref{app:memory}. Balachandran \cite{Balachandran2026}
models contamination in growing detector memories by assuming that
the conditional mean erroneous fraction among new admissions depends
only on the current contamination. For one admission at a time,
this mean is the probability of an erroneous addition, so the
assumption gives an exact binary HLS rule. Its affine case has the
trajectory in Eq.~(\ref{eq:memory_affine_trajectory}). For batches,
the conditional mean alone does not specify the admission law and
does not by itself give the binary path probabilities used here.

\subsubsection{Lexical innovation}\label{app:ml-innovation}

An open vocabulary leads to a different binary observable: set
$\sigma_n=1$ when the $n$th word has not appeared previously, and
$\sigma_n=0$ otherwise. Thus $n\psi_n$ counts distinct words and
$\sigma_1=1$. Tani Raffaelli, Lalli and Tria use an innovation model
related to the two--parameter Poisson--Dirichlet, or Pitman--Yor,
process \cite{PitmanYor1997} for authorship attribution
\cite{TaniRaffaelli2024}.
We denote the discount by $p$, written $\alpha$ in
\cite{PitmanYor1997}. For $0<p<1$ and concentration $\theta\geq0$,
the novelty probability and limiting urn function are
\begin{equation}
\pi_n(\psi)=\frac{\theta+p n\psi}{\theta+n},
\qquad \pi(\psi)=p\psi.
\label{eq:ml-pitman-rule}
\end{equation}
The base distribution determines the identity selected by the
innovation mechanism. The selected word is previously unseen when
that distribution is nonatomic, assigning zero probability to each
individual label, or when it is renormalized over unseen words as in
\cite[Eq.~(5)]{TaniRaffaelli2024}, with positive unseen mass after every
finite history. The HLS representation
concerns the novelty indicators of this generator. The implemented
classifier also uses word identities and an author--dependent base
distribution, whose probabilities are not determined by $\psi_n$.

For $\theta=0$, Eq.~(\ref{eq:ml-pitman-rule}) is exactly the
Bernoulli elephant urn already discussed above.
For $\theta>0$, the comparison below gives the same path LDP and
scaled CGF. The only fixed share is zero, which is attracting;
Eq.~(\ref{eq:app-affine-trajectory}) describes its macroscopic
approach. Taking expectations in the exact update gives
\begin{equation}
\mathbb E(n\psi_n)=O(n^p).
\label{eq:ml-pitman-mean}
\end{equation}
The vanishing share therefore describes sublinear vocabulary growth,
not the disappearance of new words. Indeed, if the count were to
remain fixed, its future novelty probabilities would have a divergent
sum, so another novelty occurs almost surely.

The boundary affine case also has an explicit CGF
\cite{BercuFavaro2026}:
\begin{equation}
\zeta(\beta)=
\begin{cases}
0,&\beta\leq0,\\
-\log[1-(1-\exp(-\beta))^{1/p}],&\beta>0.
\end{cases}
\label{eq:ml-pitman-cgf}
\end{equation}
For negative $\beta$, Jensen's inequality and
Eq.~(\ref{eq:ml-pitman-mean}) give the zero branch.
The path LDP ensures existence of the CGF, which is convex and
$1$--Lipschitz, hence locally absolutely continuous. For positive
$\beta$, Eq.~(\ref{eq:cgf-equation}) therefore holds almost everywhere
and is a linear equation in $\exp(-\zeta)$, with
$\pi(\psi)=p\psi$. Its integrating factor determines the solution
on the whole positive semiaxis. For $\theta=0$, an
all--novel continuation has probability $p^{T-1}$, and every
success probability is at most $p$. These bounds give the
positive--infinity condition in Eq.~(\ref{eq:app-affine-cgf-boundaries})
with $a=0$ and $b=p$, selecting the displayed solution. The comparison
below makes it independent of $\theta$, although the finer
vocabulary fluctuations retain a dependence on this parameter.

The same observable applies beyond language, whenever repeated
observations are interspersed with discoveries. In the urn model
with triggering of Tria et al.\ \cite{Tria2014}, the urn initially
contains $N_0\geq1$ distinct colors. Each sampled ball is replaced
with $\rho>0$ additional copies; a color sampled for the first time
also introduces $\nu+1$ previously unavailable colors, with $\nu>0$.
The parameters are integers for literal balls. After $n$ draws,
the unseen colors have total mass $N_0+\nu n\psi_n$, giving
\begin{equation}
\pi_n(\psi)=\frac{N_0+\nu n\psi}{N_0+\rho n+(\nu+1)n\psi},
\qquad
\pi(\psi)=\frac{\nu\psi}{\rho+(\nu+1)\psi}.
\label{eq:ml-umt-rule}
\end{equation}
The limiting trajectory equation is
\begin{equation}
\partial_\tau\psi(\tau)
=\frac{\psi(\tau)[\nu-\rho-(\nu+1)\psi(\tau)]}
{\tau[\rho+(\nu+1)\psi(\tau)]}.
\label{eq:ml-umt-trajectory}
\end{equation}
For $\nu<\rho$, zero is attracting; at $\nu=\rho$, the linear
restoring term vanishes but the drift remains negative for positive
shares. For $\nu>\rho$, zero is repelling and the attracting share is
\begin{equation}
\psi=\frac{\nu-\rho}{\nu+1}.
\label{eq:ml-umt-fixed}
\end{equation}
As in the preceding model, discoveries continue almost surely even
when their fraction tends to zero: if their count stayed fixed,
the probability of another discovery would have a divergent
harmonic sum. The binary history also records runs of discoveries
and repetitions, studied for texts and other innovation data in
\cite{Santoro2026}. Such observables retain information about the
ordering of novelties that the terminal vocabulary size alone loses.

Time--dependent reinforcement can change this behavior. Bellina,
De Marzo and Loreto \cite{Bellina2025} first modify the preceding
model by adding $\rho-1$ copies when a color is new: its existing
ball then becomes an observed color, so the observed mass increases
by $\rho$ after either outcome. They allow reinforcement to depend
on the number of draws and triggering on the number of discoveries.
Consider reinforcement $c_\rho n$ at draw $n$ and net addition of
$c_\nu j$ unseen colors at the $j$th discovery, with positive
integers $c_\rho,c_\nu$. The number of new colors introduced then
is $c_\nu j+1$. Summing these additions gives the exact response
\begin{equation}
\pi_n(\psi)=
\frac{2N_0+c_\nu n\psi(n\psi+1)}
{2N_0+c_\rho n(n+1)+c_\nu n\psi(n\psi+1)}.
\label{eq:ml-tumt-rule}
\end{equation}
Its limiting urn function is
\begin{equation}
\pi(\psi)=\frac{c_\nu\psi^2}{c_\rho+c_\nu\psi^2}.
\label{eq:ml-tumt-limit}
\end{equation}
The limiting response gives the following fixed-point structure.
Zero is always attracting. If $c_\nu<4c_\rho$, there is no
other fixed share; if $c_\nu>4c_\rho$, there are two more:
\begin{equation}
\psi=\frac{1\pm\sqrt{1-4c_\rho/c_\nu}}{2}.
\label{eq:ml-tumt-fixed}
\end{equation}
The smaller is repelling and separates the basins of zero and the
larger, attracting share. At $c_\nu=4c_\rho$, the two meet at
$1/2$, which attracts from above and repels from below. These are
the directions of the macroscopic trajectories in
Eq.~(\ref{eq:Caucy}); a fixed microscopic seed leaves the choice
of the limiting share to the stochastic growth history.

There is also a qualitative difference at small counts. Along a
continuation without discoveries, $n\psi_n$ stays constant and
Eq.~(\ref{eq:ml-tumt-rule}) makes the novelty probabilities
$O(n^{-2})$. Their sum is finite, so the conditional probability
of no further discovery is positive. Here innovation can stop,
whereas a vanishing novelty share in the preceding models still
allows indefinitely many discoveries.

Figure~\ref{fig:ml-innovation}\label{call:fig:ml-innovation}
compares the two triggering mechanisms.%
\afterpage{%
\clearpage
\begin{figure}[p]
\centering
\includegraphics[width=\linewidth]{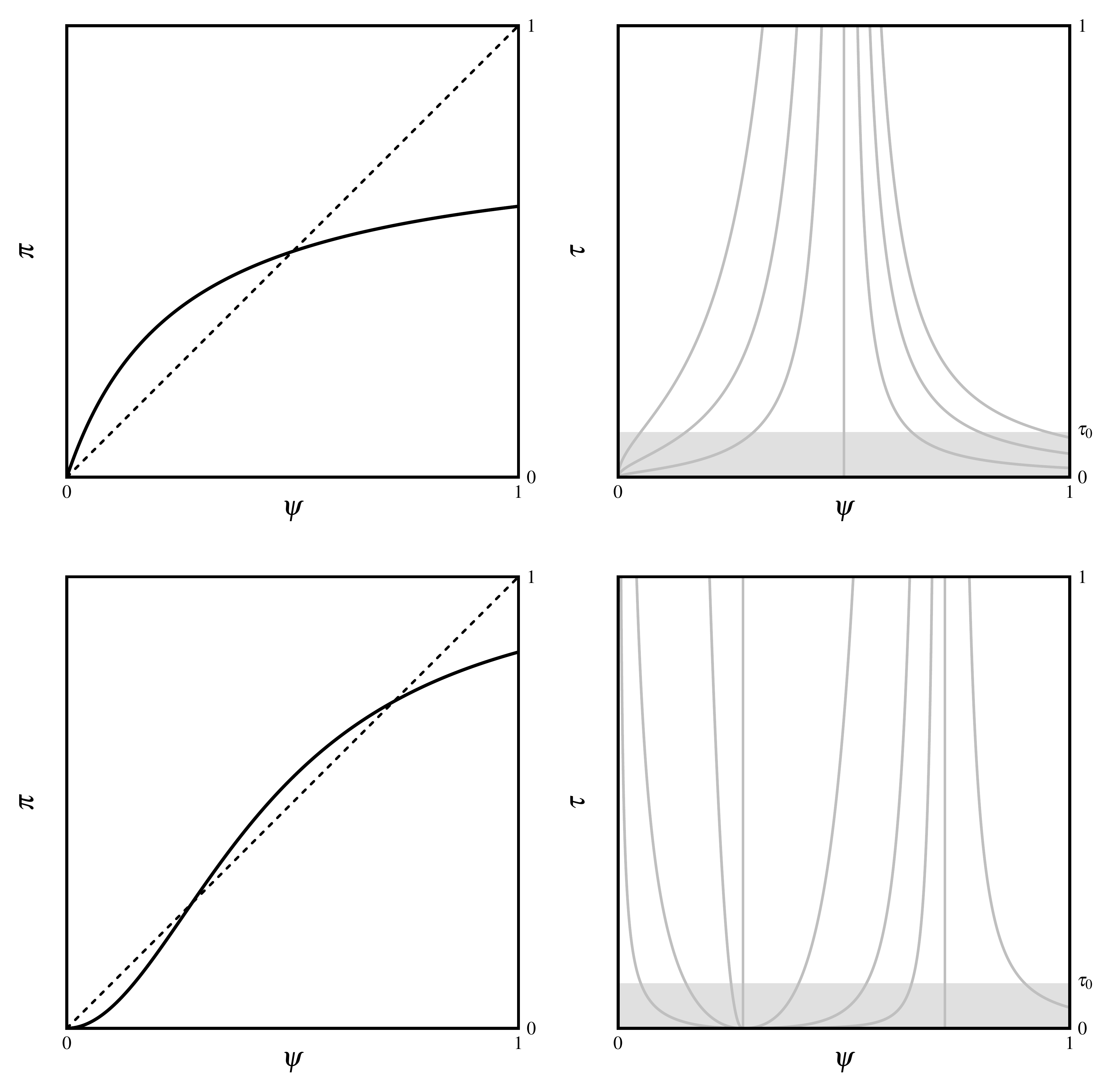}
\caption{Urn functions (left) and zero--cost share trajectories (right)
for innovation models. Upper row: the response of the urn model with
triggering in Eq.~(\ref{eq:ml-umt-rule}) with $\rho=1$ and $\nu=3$;
the attracting share is $1/2$. Lower row: the response with reinforcement
and triggering that depend on counts in Eq.~(\ref{eq:ml-tumt-limit}),
with $c_\rho=1$ and $c_\nu=5$.
Its unstable fixed share, approximately $0.276$, separates the basins
of zero and the attracting share, approximately $0.724$.
At $\tau_0=0.1$, the upper trajectories have initial shares
$0.05,0.15,0.3,0.5,0.65,0.8,0.95$; the lower trajectories have
initial shares $0.05,0.15,0.25,0.4,0.55,0.65,0.9$, together with the
two positive fixed shares. The vertical trajectories are stationary
solutions. All axes span $[0,1]$, saturation runs upwards, and the
diagonal is dashed. Below $\tau_0$, backward continuations of the
same solutions are visible through the gray strip, up to the
admissible boundary or to $\tau=0$; this strip lies before the
prescribed macroscopic initial condition.}
\label{fig:ml-innovation}
\end{figure}
\clearpage
}
In the upper row, a positive
fraction of discoveries attracts the share from both sides. In the
lower row, an unstable fixed share separates two basins: the novelty
fraction decreases towards zero on one side and approaches a
positive attracting share on the other. These are trajectories of
the limiting responses. Their geometry describes the evolution of
the fraction of novelties; the distinction between infinitely many
discoveries and their eventual arrest also uses the finite-time
probabilities examined above.

For all three models the limiting response vanishes at zero.
The comparison in Eq.~(\ref{eq:asymptotic-comparison}) therefore
requires direct control of the relative transition probabilities.
Starting from $\sigma_1=1$, every subsequent binary history has
positive probability under both the original process and its
limiting HLS urn. In particular, $n\psi_n\geq1$ throughout.
Direct division of the kernels gives, uniformly on these histories,
\begin{equation}
\log\frac{\pi_n(\psi_n)}{\pi(\psi_n)}
=O\left(\frac{1}{n\psi_n}\right),
\qquad
\log\frac{1-\pi_n(\psi_n)}{1-\pi(\psi_n)}
=O\left(\frac1n\right).
\label{eq:ml-relative-kernels}
\end{equation}
Indeed, the correction to the total urn mass is relative
$O(1/n)$. The novelty numerator has relative correction
$O(1/(n\psi_n))$; in the last model it also contains a term of
order $(n\psi_n)^{-2}$, bounded by the same estimate. The repeat
numerator is unchanged in the first two models and gains only
the factor $1+1/n$ in the third. All constants depend only on
the fixed model parameters.

A novelty increases the count by one, so its successive preceding
counts are distinct. Consequently,
\begin{equation}
\sum\limits_{1\leq n<T}\frac{\sigma_{n+1}}{n\psi_n}
=\sum\limits_{1\leq j<T\psi_T}\frac1j
\leq1+\log T.
\label{eq:ml-novelty-harmonic}
\end{equation}
Together with the harmonic bound for all steps,
Eq.~(\ref{eq:ml-relative-kernels}) therefore yields
\begin{equation}
\sup_{\sigma:\,\sigma_1=1}
\Bigl|\sum\limits_{1\leq n<T}
[L(\sigma_{n+1},\pi_n(\psi_n))-L(\sigma_{n+1},\pi(\psi_n))]\Bigr|
=O(\log T).
\label{eq:ml-innovation-comparison}
\end{equation}
This bounds the logarithmic ratio of the history probabilities
uniformly, including histories whose novelty share approaches zero.
The original models and their HLS limits consequently have the same
path LDP and CGF, with action $\Phi$, rate $\Phi_0^*-\Phi$ and the
boundary accessibility specified in Section~\ref{sec:ldp}
\cite[Corollary~2]{Franchini_URNS}. The comparison is on scale $T$;
sublinear vocabulary growth and its finer fluctuations require
their own analysis.

\subsection{Neural lattice field theory}\label{sec:lft-two-body}

The HLS urn provides a natural nonlinear benchmark for testing the accuracy
of general approximated LFT actions in the time domain, such as the two--body approximation
considered in the neural lattice field theory (NLFT) in
\cite{BardellaFranchiniShort2024,BardellaFranchini2024,FranchiniBardella2026}.
Its exact history law provides a reference for quantifying how truncations
of the action affect temporal trajectories and the free energy.
This LFT describes complete binary activity histories through an action
on a space--time lattice. Its fields and couplings determine mean
activity together with spatial and temporal correlations. The construction
also applies to binary systems outside neuroscience.

The reduced action follows the assumptions detailed in
\cite[Sections~5--6]{FranchiniBardella2026}: truncation to one--body
and two--body terms; causal ordering of temporal couplings; local
memory, excluding interactions between different neurons at different
times; and bi--stationarity, requiring time--independent couplings
between neurons and neuron--independent temporal couplings.
For a single neuron, the conditions involving distinct neurons
are automatic. The HLS example therefore tests the temporal two--body
approximation, with couplings that may depend on both observation times.

\subsubsection{Activity model and history law}

For $N=1$, let $\sigma_n$ indicate activity in time bin $n$,
$\psi_n$ the fraction of active bins among the first $n$, and $T$ the
number of observed bins. The neuron index is suppressed throughout.
The HLS response gives the conditional
probability of activity in the next bin. Consider the simple affine urn
\begin{equation}
\pi(\psi)=1/2+b(\psi-1/2),\qquad 0<b<1.
\label{eq:lft-affine}
\end{equation}
The probability of activation increases with the fraction of active
bins in the past. The slope $b$ controls this feedback, while $1/2$
remains the attracting activity fraction. Both outcomes have positive
probability for every history. We take an unbiased first bit, so the
law is invariant under simultaneous exchange of activity and inactivity.

The conditional probabilities determine the probability of a complete
recording by multiplication. Its logarithm is precisely the action in
Eq.~(\ref{eq:path-action}), with the Bernoulli Lagrangian
(\ref{eq:path-lagrangian}). Although the response is affine, these
logarithms generate terms of arbitrarily high degree in the past
activities. The two--body approximation retains only a quadratic action, discarding
these higher--order terms even for this elementary urn.

The formulas in Section~\ref{app:classical} give the reference
trajectories and limiting CGF without truncating the action. We compare
two observables: the typical continuation after an assigned initial
history and the exponential fluctuations of the total activity count.
The former uses a prescribed share at a positive saturation $\tau_0$;
the latter starts from the unbiased first bit. These different initial
conditions will also be imposed on the two--body model.

\subsubsection{Two--body approximation and temporal couplings}

For the affine response (\ref{eq:lft-affine}), the exact action
(\ref{eq:path-action}) is
\begin{multline}
\mathcal A(\sigma)=-T\log2+
\sum\limits_{1<n\leq T}\sigma_n\log\Bigl[1+\frac{b}{n-1}
\sum\limits_{1\leq m<n}(2\sigma_{m}-1)\Bigr]\\
+\sum\limits_{1<n\leq T}(1-\sigma_n)\log\Bigl[1-\frac{b}{n-1}
\sum\limits_{1\leq m<n}(2\sigma_{m}-1)\Bigr].
\label{eq:lft-exact-action}
\end{multline}
The time indices belong to $S$, as in the theory. The term $-T\log2$
includes the unbiased first bit and the factors $1/2$ in all subsequent
transition probabilities.

We apply the field expansion of
\cite[Sections~4.1.1--4.1.2]{BardellaFranchini2024} to the action in
Eq.~(\ref{eq:lft-exact-action}), expanding in the centered activities
$\sigma_n-1/2$, holding $b$ fixed. First, regard the Bernoulli
Lagrangian as a function of real activities near $1/2$. Then retain
terms of degree at most two and evaluate the resulting polynomial
on binary histories.
Centering at $1/2$ preserves the symmetry between activity and
inactivity; for $0<b<1$, the one--step logarithmic series converges uniformly
over all binary histories. Expanding the products and using
$\sigma_n^2=\sigma_n$ gives the following constant, linear
and quadratic contributions. Here $\psi_n$ is the fraction of active
bins among the first $n$ observations, defined in Eq.~(\ref{eq:urn-share}).
Each bit appears in these fractions from its own time onwards; each
pair appears in their squares once both bits are included. The sums
in brackets collect these contributions up to time $T-1$.
In pair terms, $s$ is the time in $\psi_s$; $n$ and $m$ are the
times of the two bits:
\begin{multline}
\mathcal A_2(\sigma)=-T\log2+\frac{b(2-b)(T-1)}2\\
-2b\sum\limits_{1\leq n\leq T}\Bigl[\mathbb I(n>1)
+\sum\limits_{n\leq m<T}\frac{(1-b)m+b}{m^2}\Bigr]\sigma_n\\
+4b\sum\limits_{1<n\leq T}\sum\limits_{1\leq m<n}\Bigl[\frac1{n-1}
-b\sum\limits_{n\leq s<T}\frac1{s^2}\Bigr]\sigma_n\sigma_{m}.
\label{eq:lft-quadratic-action}
\end{multline}
The terms can be traced directly to the quadratic part of the
one--step Lagrangian:
\begin{equation}
\mathcal L_2(\sigma_{n+1},\psi_n)=-\log2+b(2\sigma_{n+1}-1)(2\psi_n-1)
-\frac{b^2}{2}(2\psi_n-1)^2.
\label{eq:lft-quadratic-lagrangian}
\end{equation}
The product comes from the activity--dependent logarithm, while the
squared share comes from the symmetric part of the two logarithms.
Inserting the empirical share and summing over time yields
Eq.~(\ref{eq:lft-quadratic-action}). In particular, the square of a
share contains both pairs of distinct past bits and diagonal terms;
the latter become linear after restriction to binary activities.

In NLFT \cite[Eqs.~(6), (23)--(24)]{FranchiniBardella2026},
the off--diagonal entries of the matrix $A$ couple different neurons
at the same time and encode their spatial correlations. With a single
neuron, only the diagonal coefficient $A_{11}$ remains. Its contribution
involves the equal--time self--correlation $\sigma_n^2=\sigma_n$, so it
can be absorbed into the input field $I^n$ acting on activity in bin $n$.
In comparing $\mathcal A_2$ with NLFT, we therefore choose $A_{11}=0$
after this absorption. Correlations of the same neuron at distinct times
remain represented by $B^{nm}$, which couples its activities in bins
$n$ and $m$. Since our action is a logarithmic probability
weight, the coefficients below have opposite signs to those of the
linear and pair terms in $\mathcal A_2$. For each fixed $n$, both
sums defining $I^n$ run over $m$ from $n$ to $T-1$:
\begin{equation}
I^n=2\Bigl[b\,\mathbb I(n>1)
+(b-b^2)\sum\limits_{n\leq m<T}\frac1{m}
+b^2\sum\limits_{n\leq m<T}\frac1{m^2}\Bigr].
\label{eq:lft-linear-field}
\end{equation}
\begin{equation}
B^{nm}=-4\mathbb I(m<n)\Bigl[\frac{b}{n-1}
-b^2\sum\limits_{n\leq s<T}\frac1{s^2}\Bigr].
\label{eq:lft-temporal-coupling}
\end{equation}
The indicator $\mathbb I(n>1)$ in $I^n$ removes the direct transition
term for the unbiased first bit; the other field terms collect each bit's appearances in
later shares, including their diagonals. A temporal pair first appears when
its later bit is generated and then enters each subsequent squared
share. These contributions give the two terms in $B^{nm}$.
Every earlier bin $m$ has the same weight at fixed $n$ because
the share weights all past activities uniformly. Each pair is counted
once: the bracket is evaluated only for $m<n$, and the indicator sets
all remaining entries of $B$ to zero. The fields and pairs together preserve the exchange of activity and inactivity, despite
the nonzero linear coefficients.

On binary histories, the remainder of the logarithmic series gives
\begin{equation}
\sup_{\sigma\in\Omega^S}|\mathcal A(\sigma)-\mathcal A_2(\sigma)|
\leq\frac{(T-1)b^3}{3(1-b)}.
\label{eq:lft-remainder}
\end{equation}
Indeed, after extracting $-\log2$, each logarithm has argument
between $1-b$ and $1+b$;
the tail after its quadratic term is bounded by a geometric series.
Thus the coefficients and a uniform error bound follow from the
specified activity law, without fitting empirical moments.

Small feedback $b$ therefore controls the action error per step
uniformly over histories. At fixed $b$, the local remainder starts
at cubic order in $b(2\psi_n-1)$: empirical shares close to $1/2$
provide another regime of accuracy. The small quantity is the deviation of the past activity fraction
from $1/2$; each centered binary activity still has magnitude $1/2$.
Convergence for $0<b<1$ alone does not ensure that truncation after
two--body terms is accurate. A longer recording does not, by itself,
make it exact. The bound controls logarithmic weights. Bounding the displacement
of the maximizing trajectories requires a separate estimate of how
the maximizer changes when the action is perturbed.

\subsubsection{Solving the two--body model}

The two--body law is normalized over complete histories, as in
Eq.~(\ref{eq:path-gibbs}). Normalizing each transition separately
changes the law: the sum of the two extension weights depends on the
current share. Dividing by that sum therefore changes relative history
weights. The conditional probabilities of the globally normalized law
also depend on the weights of future continuations.
Separate normalization at each step would give precisely the zero--field
magnetic--growth urn of Eq.~(\ref{eq:magneticurn}), with $J=b$.

The finite computation groups histories by activity count. Let
$Z_{2,n}(j)$ be the total unnormalized two--body weight of histories
with $j$ activations in $n$ bins. Begin with the unbiased first bit,
$Z_{2,1}(0)=Z_{2,1}(1)=1/2$. Each attainable count then receives two
contributions:
\begin{equation}
Z_{2,n+1}(j)=Z_{2,n}(j)\exp[\mathcal L_2(0,j/n)]
+Z_{2,n}(j-1)\exp[\mathcal L_2(1,(j-1)/n)].
\label{eq:lft-count-recursion}
\end{equation}
Terms with an unattainable predecessor count are omitted. The sum of
the terminal weights gives the normalization; an exponential weight
on the terminal count gives the CGF in the convention of
Eq.~(\ref{eq:cgf}):
\begin{equation}
\zeta_2(\beta)=\lim_{T\to\infty}\frac1T\log
\frac{\sum\limits_{0\leq j\leq T}Z_{2,T}(j)\exp(\beta j)}
{\sum\limits_{0\leq j\leq T}Z_{2,T}(j)}.
\label{eq:lft-normalized-cgf}
\end{equation}
For a prescribed prefix, one instead initializes weight one at its
terminal count and applies the same recursion to the remaining bins.
This computes the conditional law, since the weight of the fixed
prefix cancels in its normalization.

For the unbiased first bit, the normalization is subexponential:
\begin{equation}
-\frac{b^2}{2}\sum\limits_{1\leq n<T}\frac1n
\leq\log\sum\limits_{0\leq j\leq T}Z_{2,T}(j)\leq0.
\label{eq:lft-normalization-bound}
\end{equation}
The lower bound follows from Jensen's inequality under independent
unbiased bits. The upper bound follows by summing the two extensions
and using $\cosh u\leq\exp(u^2/2)$.
For the unbiased first bit, this subexponential normalization allows
the bound in Eq.~(\ref{eq:lft-remainder}), divided by $T$, to control
the limiting CGF error uniformly in $\beta$.
A common mean activity fraction $1/2$ does not guarantee accurate
free energies for large $|\beta|$: the exponential weighting then
favors histories whose terminal shares are far from $1/2$.

Define $\Phi_2$ from
Eq.~(\ref{eq:scaledaction}) by using $\mathcal L_2(\partial_\tau\varphi,\psi)$
as the integrand. Applying the change of measure and path variational
principle of Section~\ref{sec:ldp} gives
\begin{equation}
\zeta_2(\beta)=\sup_{\varphi\in Q}
\{\beta\varphi(1)+\Phi_2(\varphi)-\Phi_0^*(\varphi)\}.
\label{eq:lft-variational-cgf}
\end{equation}
Here the paths start at zero and have derivatives between zero and
one, as in the theory. The binary mixing term remains essential:
a quadratic action on bits does not turn their path law into a
Gaussian model. Subtracting the quadratic part from the centered
Bernoulli rate function leaves a strictly convex function. Completing the
square in the remaining terms therefore shows that this variational
objective is strictly concave, with a unique maximizing path.
At $\beta=0$ its maximum is zero, attained by the balanced path.

The terminal share of the maximizer equals
$\partial_\beta\zeta_2(\beta)$. Rescaling this path in the variational
problem, or integrating its Euler--Lagrange equation, yields
\begin{equation}
\zeta_2(\beta)=\frac\beta2+
\log\cosh\Bigl[\frac\beta2+b(2\partial_\beta\zeta_2(\beta)-1)\Bigr]
-\frac{b^2}{2}(2\partial_\beta\zeta_2(\beta)-1)^2.
\label{eq:lft-quadratic-cgf}
\end{equation}
The CGF satisfies $\zeta_2(-\beta)=\zeta_2(\beta)-\beta$ and $\zeta_2(0)=0$.
The positive branch is fixed by
\begin{equation}
\lim_{\beta\to+\infty}[\zeta_2(\beta)-\beta]
=-\log2+b-b^2/2.
\label{eq:lft-quadratic-boundary}
\end{equation}
This is an exact differential characterization for the two--body
model, evaluated numerically below, rather than an elementary closed
formula.

To compare trajectories, fix the initial count through
$\varphi(\tau_0)=\tau_0\psi_0$ and maximize
$\Phi_2(\varphi)-\Phi_0^*(\varphi)$ over the remaining interval. Since $\psi=\varphi/\tau$, the
Euler--Lagrange equation becomes
\begin{equation}
\tau\partial_\tau^2\varphi
=-2b(1-b)(2\psi-1)\partial_\tau\varphi(1-\partial_\tau\varphi),
\qquad \tau_0<\tau<1.
\label{eq:lft-path-equation}
\end{equation}
With the terminal term $\beta\varphi(1)$, not fixing the terminal share
gives the boundary condition
\begin{equation}
\log\frac{\partial_\tau\varphi(1)}{1-\partial_\tau\varphi(1)}
=\beta+2b(2\psi(1)-1).
\label{eq:lft-path-boundary}
\end{equation}
Set $\beta=0$ for the typical unweighted continuation. Solving this
boundary problem gives the two--body trajectories; the HLS reference
instead follows Eq.~(\ref{eq:app-affine-trajectory}). For a
macroscopic prefix, the normalization need not have zero limiting
logarithm per step, but it is independent of the continuation and
therefore does not change its maximizing path.

\subsubsection{Comparison with the exact urn}

We now give the trajectories and free energies for the benchmark
$b=1/2$, with activation probability $\pi(\psi)=1/4+\psi/2$.
As above, subscript $2$ denotes the two--body model. For trajectories
we prescribe the same initial share $\psi_0$ at saturation $\tau_0>0$;
for free energies we use the unbiased first bit. The distinction
matters because a prescribed macroscopic history changes the
normalization of its subsequent continuation.

The HLS share follows directly from
Eq.~(\ref{eq:app-affine-trajectory}):
\begin{equation}
\psi(\tau)=1/2+(\psi_0-1/2)\sqrt{\tau_0/\tau},
\qquad \tau_0\leq\tau\leq1.
\label{eq:lft-comparison-hls-path}
\end{equation}
The normalized activity count is $\varphi(\tau)=\tau\psi(\tau)$.
For the two--body model, integrating the logarithm of the odds of
$\partial_\tau\varphi$ in Eqs.~(\ref{eq:lft-path-equation})--(\ref{eq:lft-path-boundary})
with $\beta=0$, and then integrating the activity rate, gives
\begin{equation}
\tau\psi_2(\tau)=\tau_0\psi_0+\frac12\int_{\tau_0}^{\tau}
\Bigl\{1+\tanh\Bigl[\psi_2(1)-1/2+
\frac14\int_s^1\frac{2\psi_2(u)-1}{u}\,du\Bigr]\Bigr\}\,ds.
\label{eq:lft-comparison-quadratic-path}
\end{equation}
This is an exact integral equation for the whole trajectory,
including its terminal share $\psi_2(1)$. It incorporates the same
initial condition and the free terminal condition; its solution is
unique by the concavity established above. The integration variables
$s$ and $u$ are saturation times. The count trajectory is
$\varphi_2(\tau)=\tau\psi_2(\tau)$.
The dependence on the terminal share reflects normalization over
complete continuations. Both models retain the constant solution
$\psi=1/2$, but their relaxation from other initial shares differs.
For $\tau_0=0.1$ and $\psi_0=0.9$, the final HLS and two--body
shares are approximately $0.62649$ and $0.62276$; the largest
difference along the trajectory is about $0.0052$.

For the free energies, we return to the unbiased microscopic seed.
Evaluating the integral in Eq.~(\ref{eq:app-affine-cgf-positive})
with $a=1/4$ and $b=1/2$ gives, for $\beta>0$,
\begin{equation}
\zeta(\beta)=\beta/2-\log\Bigl[\cosh(\beta/2)
+\sinh^2(\beta/2)\log\tanh(\beta/4)\Bigr].
\label{eq:lft-comparison-hls-cgf}
\end{equation}
For the two--body model, Eq.~(\ref{eq:lft-quadratic-cgf}) reduces to
\begin{equation}
\zeta_2(\beta)=\beta/2+
\log\cosh\Bigl[\beta/2+\partial_\beta\zeta_2(\beta)-1/2\Bigr]
-\frac12\Bigl[\partial_\beta\zeta_2(\beta)-1/2\Bigr]^2.
\label{eq:lft-comparison-quadratic-cgf}
\end{equation}
Its positive branch is selected by
\begin{equation}
\lim_{\beta\to+\infty}[\zeta_2(\beta)-\beta]=3/8-\log2.
\label{eq:lft-comparison-quadratic-boundary}
\end{equation}
Thus the HLS free energy has the explicit expression above, while
the two--body free energy is determined by a first--order differential
equation with this boundary condition. The latter is an exact
characterization of the two--body model, evaluated numerically;
quadratic dependence of the action does not remove the binary mixing
entropy. Both free energies have continuous value zero at $\beta=0$,
and exchanging activity and inactivity determines their negative branches:
\begin{equation}
\zeta(-\beta)=\zeta(\beta)-\beta,
\qquad \zeta_2(-\beta)=\zeta_2(\beta)-\beta.
\label{eq:lft-comparison-cgf-symmetry}
\end{equation}
At $\beta=1$, the values are approximately $0.79355$ for HLS
and $0.77852$ for the two--body model. These free energies describe
the exponential weights of atypical activity counts, even though
both models have mean activity fraction $1/2$.

An analytic comparison also holds for every slope $0<b<1$ in
Eq.~(\ref{eq:lft-affine}). Combining
Eqs.~(\ref{eq:app-affine-cgf-boundaries}) and
(\ref{eq:lft-quadratic-boundary}) gives
\begin{equation}
\lim_{\beta\to+\infty}[\zeta(\beta)-\zeta_2(\beta)]
=\log(1+b)-b+b^2/2>0.
\label{eq:lft-cgf-discrepancy}
\end{equation}
This is the action error per transition on a fully active history.
For $b=1/2$, the limiting gap is $\log(3/2)-3/8$, approximately
$0.03047$. The discrepancy therefore persists in the limiting free
energy, as well as in the conditional trajectories.

\penalty0
\afterpage{\clearpage
\begin{figure}[p]
\centering
\makebox[\linewidth][l]{\hspace*{-6.177927pt}\includegraphics[width=371.902464pt]{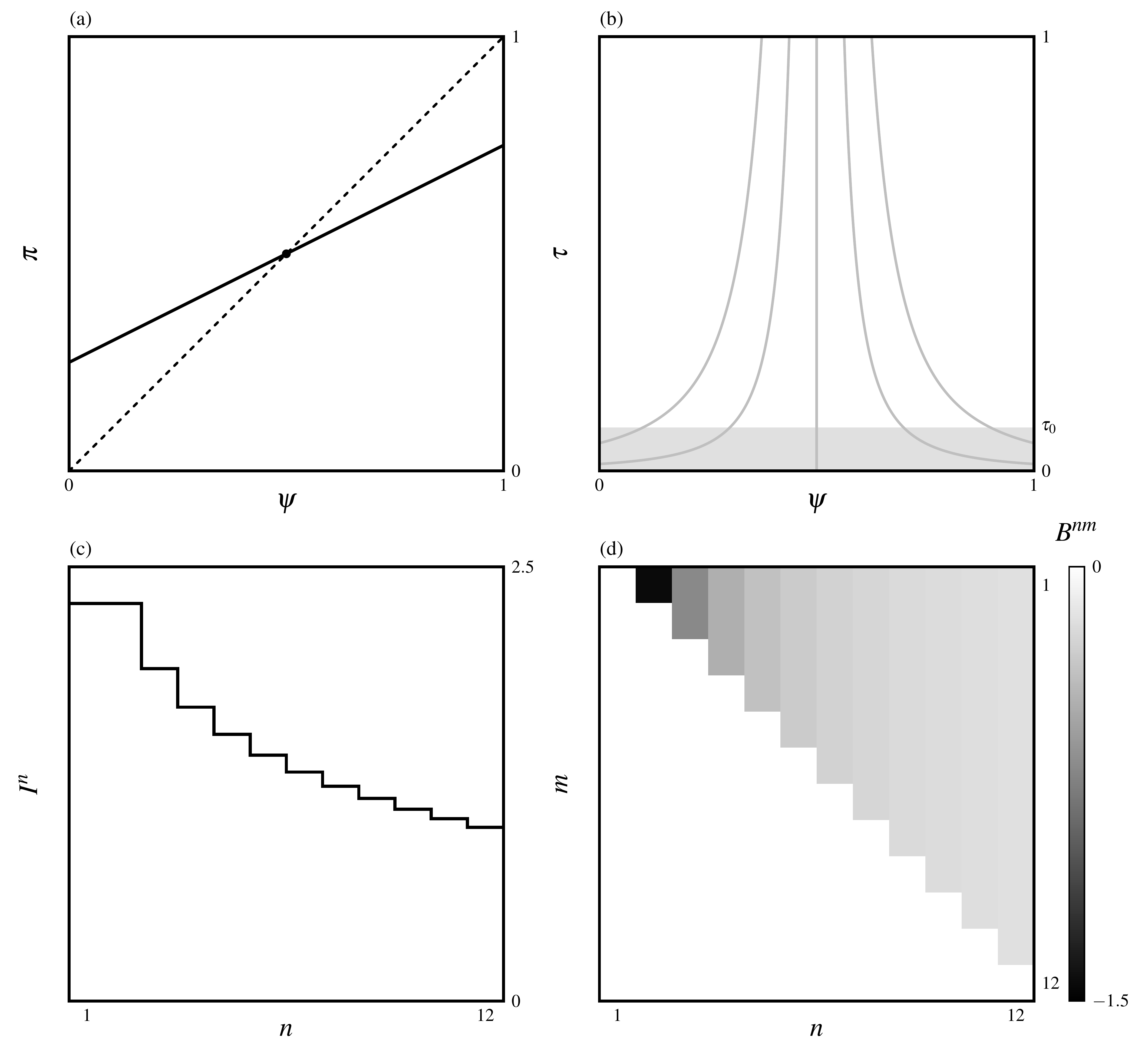}}
\caption{Affine HLS response and the two--body approximation of its action in neural LFT at $b=1/2$.
(a) Urn function $\pi(\psi)=1/4+\psi/2$ and the dashed diagonal;
the intersection gives the attracting share $1/2$.
(b) HLS share trajectories from Eq.~(\ref{eq:lft-comparison-hls-path}),
with $\tau_0=0.1$ and initial shares $0.1,0.3,0.5,0.7,0.9$.
The shaded initial strip contains backward continuations of the same
curves while their shares remain in $[0,1]$.
(c) Field coefficients $I^n$ from Eq.~(\ref{eq:lft-linear-field}),
shown as steps for $T=12$; the first two coincide, and the field decreases thereafter.
(d) Temporal coupling matrix $B^{nm}$ from
Eq.~(\ref{eq:lft-temporal-coupling}) for the same $T$.
Time $n$ increases to the right and $m$ downwards.
Each pair is counted once: cells with $m\geq n$ are zero (white),
while those with $m<n$ contain negative couplings.
The nonzero values within each column are constant because all past bins
have equal weight in the empirical share. The lateral scale gives the coupling values.}
\label{fig:lft-two-body}
\end{figure}
\clearpage
\begin{figure}[p]
\centering
\makebox[\linewidth][l]{\hspace*{-8.806833pt}\includegraphics[width=364.015748pt]{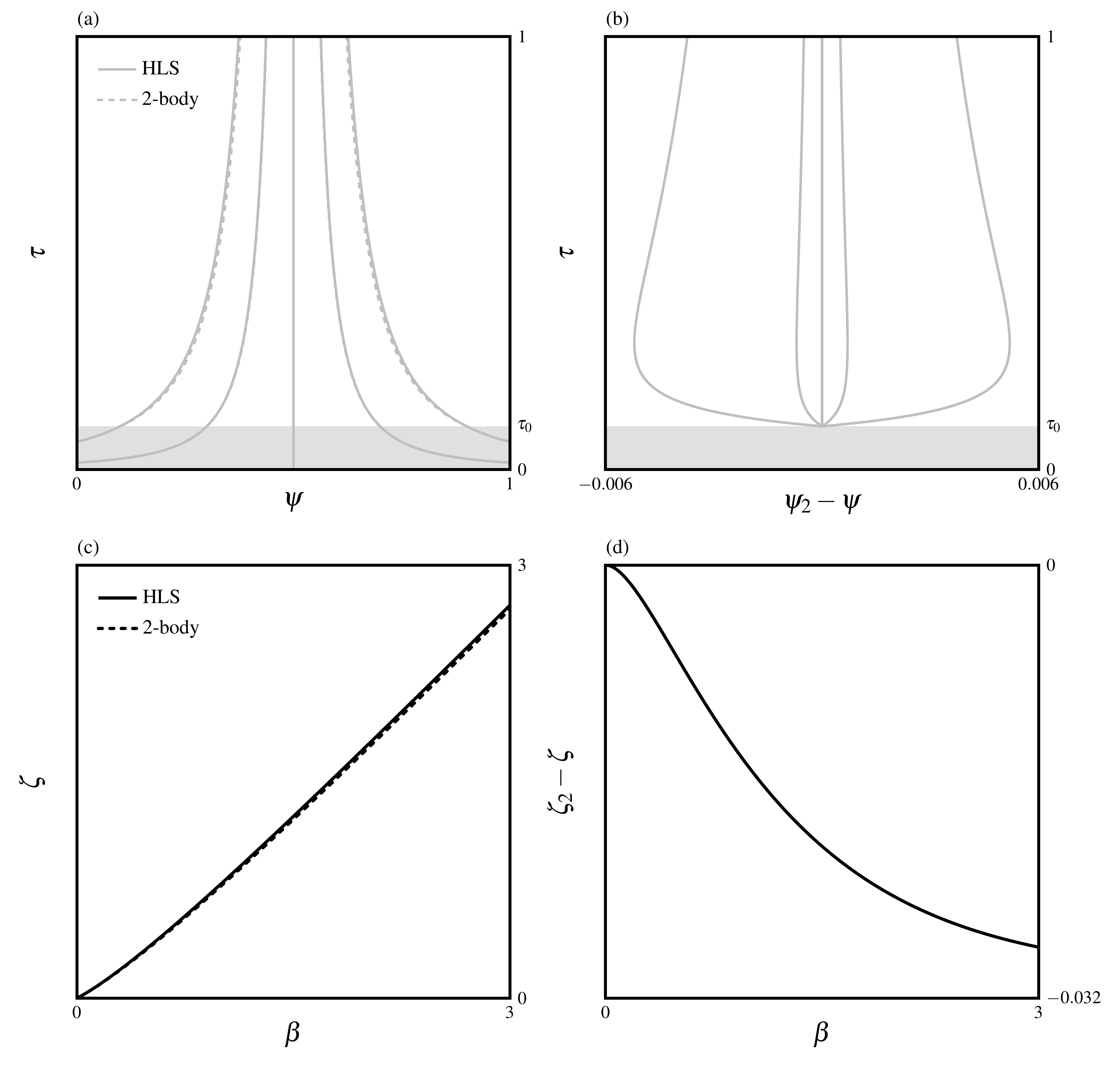}}
\caption{Trajectories, free energies and their differences for HLS and
the two--body approximation of its action in neural LFT at $b=1/2$.
(a) Typical continuations of histories with the same initial shares
$0.1,0.3,0.5,0.7,0.9$ at $\tau_0=0.1$, from
Eqs.~(\ref{eq:lft-comparison-hls-path}) and
(\ref{eq:lft-comparison-quadratic-path}).
Continuous curves denote HLS and dashed curves the two--body model,
as in panel (c).
(b) Differences $\psi_2-\psi$ for the same trajectories as in (a).
The balanced path has zero difference; for initial share $0.9$, the
largest magnitude of $\psi_2-\psi$ is about $0.0052$, and the final shares are
$0.62649$ and $0.62276$ for HLS and the two--body model.
In both upper panels, saturation spans $[0,1]$ and the shaded strip
marks the interval before the prescribed initial condition at $\tau_0$.
Only the solid HLS curves in (a) are continued backwards through this
strip, until their shares reach the admissible boundary or $\tau=0$.
(c) Limiting free energies for the unbiased first bit, from
Eqs.~(\ref{eq:lft-comparison-hls-cgf})--(\ref{eq:lft-comparison-quadratic-boundary}).
(d) Free-energy difference $\zeta_2-\zeta$ for the curves in (c).
At $\beta=1$, the two free energies are approximately $0.79355$
and $0.77852$. These differences persist as $T\to\infty$ and measure
the effect of truncating the action.}
\label{fig:lft-comparison}
\end{figure}
\clearpage}

Figure~\ref{fig:lft-two-body}\label{call:fig:lft-two-body}
shows the affine response and its HLS trajectories above the field
$I^n$ and temporal coupling matrix $B^{nm}$ of the two--body model.
The field varies with the observation time, while the pair couplings
are constant across earlier bins at each fixed later time.
Figure~\ref{fig:lft-comparison}\label{call:fig:lft-comparison}
compares the trajectories in the upper row and the free energies
in the lower row, with each difference shown in the adjacent panel.
The example isolates the effect of the two--body approximation while
preserving equilibrium activity. With $N=1$, it does not test the
additional restrictions on interactions between distinct neurons.

\subsection{Applications of quenched random urns}\label{app:quenched}

We now apply Section~\ref{sec:quenched} to four families. In each
case, the random parameters are drawn once, their law is independent
of $T$, and the resulting response is held fixed throughout growth.
Bold symbols identify random parameters, their derived functions and
the random share $\randombold{\psi}_n$, as in Section~\ref{sec:quenched}.
Fixed arguments and quantities evaluated at a fixed realization remain
in ordinary type.

\subsubsection{Quenched Bagchi--Pal and affine urns}\label{app:quenched-affine}

Bandyopadhyay and Thacker \cite[Section~6]{BandyopadhyayThacker2022}
formulate urns whose replacement kernel is chosen once, defining the
conditional law as quenched. For two colors, the extraction-count
reduction of Section~\ref{app:classical} gives affine feedback, with
an initial-mass correction of order $1/n$ for general seeds; the exact
finite-time representation retains the stated clock and seed conditions.

Choose a balanced Bagchi--Pal replacement matrix once for each urn
and retain it at every extraction. Its entries and common row sum
are random parameters, denoted by bold symbols:
\begin{equation}
\randombold{a}_{11}+\randombold{a}_{12}=\randombold{a}_{21}+\randombold{a}_{22}=\randombold{M},
\qquad \randombold{M}>0.
\label{eq:quenched-bagchi-pal-balance}
\end{equation}
For distinct rows and a compatible finite seed, the same exact
change of variables as in Section~\ref{app:classical} gives
\begin{equation}
\randombold{a}=\randombold{a}_{21}/\randombold{M},
\qquad \randombold{b}=(\randombold{a}_{11}-\randombold{a}_{21})/\randombold{M}.
\label{eq:quenched-bagchi-pal-coefficients}
\end{equation}
Figure~\ref{fig:quenched-affine}\label{call:fig:quenched-affine} shows the full parameter domain of the
conditional HLS response
\begin{equation}
\randombold{\pi}(\psi)=\randombold{a}+\randombold{b}\psi,
\qquad 0<\randombold{a}<1,\qquad 0<\randombold{a}+\randombold{b}<1.
\label{eq:quenched-affine-urn}
\end{equation}
This domain is used below for the uniform average and contains
the two negative--slope Bagchi--Pal examples considered below.

\afterpage{%
\clearpage
\begin{figure}[p]
\centering
\makebox[\linewidth][l]{\hspace*{79.655829pt}\includegraphics[width=181.219202pt]{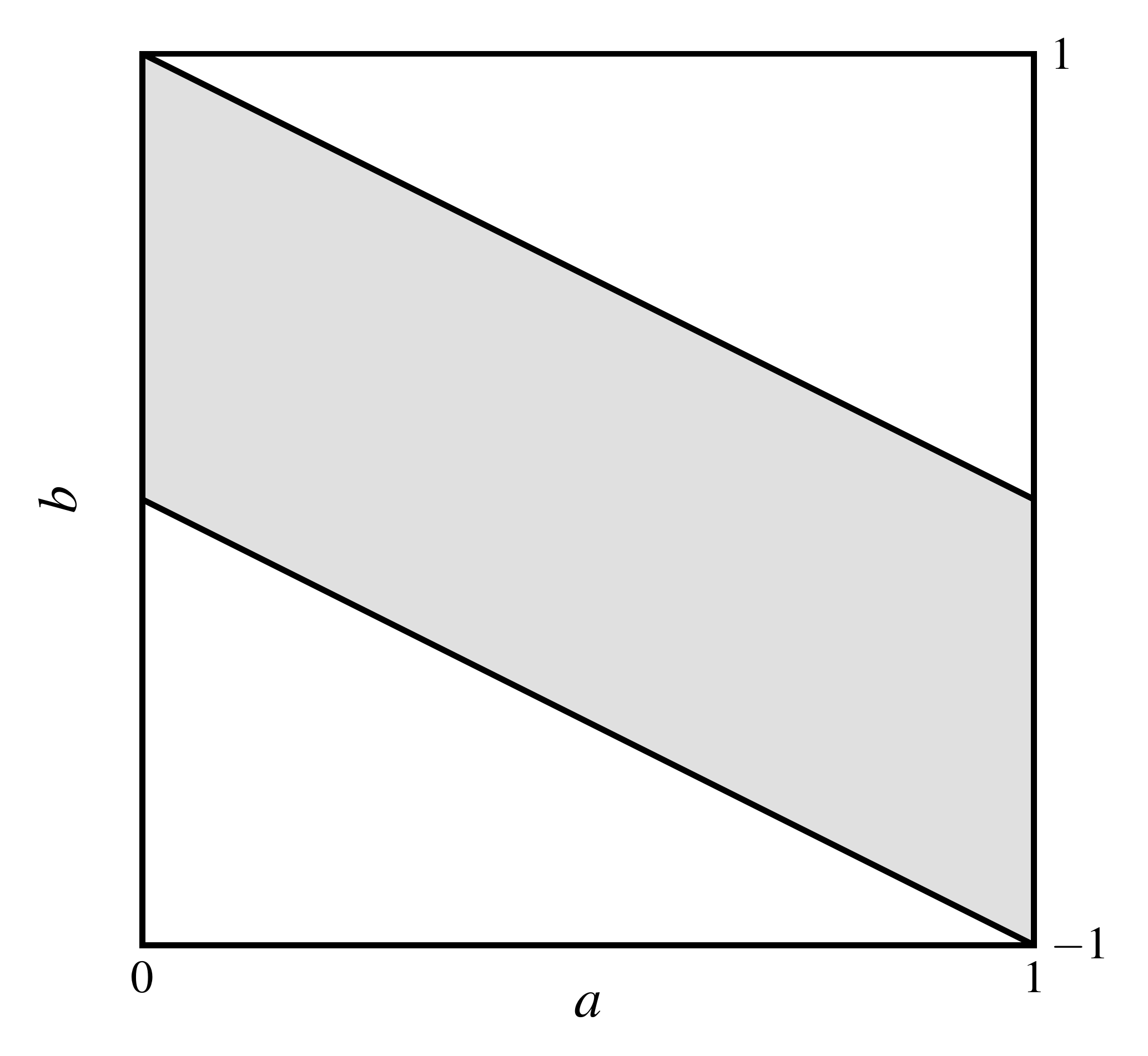}}
\caption{Full admissible parameter domain for the affine response $\pi(\psi)=a+b\psi$. The entire shaded parallelogram, $0<a<1$, $-a<b<1-a$, carries the uniform law of Eq.~(\ref{eq:quenched-affine-domain}), with density one. The box spans $0\leq a\leq1$ horizontally and $-1\leq b\leq1$ vertically. The outline delimits the open domain; its boundary has zero measure. Each sampled coefficient pair is retained throughout growth.}
\label{fig:quenched-affine}
\end{figure}
\clearpage
}

These inequalities keep the urn probabilities strictly between zero and one. Both signs
of $\randombold{b}$ are allowed; the realized function is invertible when
$\randombold{b}\ne0$. Its attracting composition is itself random:
\begin{equation}
\randombold{\psi}_*=\randombold{a}/(1-\randombold{b}).
\label{eq:quenched-affine-fixed}
\end{equation}
For a prescribed macroscopic starting share, the zero--cost trajectory
varies with the chosen coefficients:
\begin{equation}
\randombold{\psi}(\tau)=\randombold{\psi}_*+(\psi_0-\randombold{\psi}_*)(\tau_0/\tau)^{1-\randombold{b}},
\qquad 0<\tau_0\leq\tau.
\label{eq:quenched-affine-trajectory}
\end{equation}
With a fixed finite seed, the stochastic share converges almost
surely to $\randombold{\psi}_*$. Its variance over coefficients and
growth therefore tends to $\operatorname{Var}(\randombold{\psi}_*)$,
where the latter variance is over the coefficient law. Even when every urn has the
same fixed point, different slopes give different relaxation and
fluctuations. The conditional CGF $\randombold{\zeta}$ is obtained by
substituting $\randombold{a},\randombold{b}$ into
Eqs.~(\ref{eq:app-affine-cgf-positive})--(\ref{eq:app-affine-cgf-decreasing});
the quenched average is then taken over these fixed coefficients.

An integer--valued Bagchi--Pal example admits a closed expression
for that average. Choose the two symmetric replacement matrices
with equal probabilities:
\begin{equation}
\mathbb P_\pi\left((\randombold{a}_{ij})=\begin{pmatrix}1&3\\3&1\end{pmatrix}\right)=1/2,
\qquad
\mathbb P_\pi\left((\randombold{a}_{ij})=\begin{pmatrix}1&2\\2&1\end{pmatrix}\right)=1/2.
\label{eq:quenched-bagchi-pal-law}
\end{equation}
At each fixed matrix, $a,b,\zeta$ denote its coefficients and CGF.
Choose its row sum as the initial physical count of each color. This gives the same HLS seed in both cases:
$m=2$ and $m\psi_m=1$. The coefficient pairs are respectively
$(a,b)=(3/4,-1/2)$ and $(2/3,-1/3)$, so both fixed
points are $1/2$. Evaluating Eq.~(\ref{eq:app-affine-cgf-decreasing})
for the first matrix gives
\begin{equation}
\zeta(\beta)=\beta/2+2\log\cosh(\beta/4).
\label{eq:quenched-bagchi-pal-first-cgf}
\end{equation}
For the second matrix it gives
\begin{equation}
\zeta(\beta)
=\beta/2+\log\frac{4\sinh^3(\beta/2)}{3(\sinh\beta-\beta)}.
\label{eq:quenched-bagchi-pal-second-cgf}
\end{equation}
The ratio is positive for $\beta\ne0$ and tends to one at zero.
The quenched average is consequently explicit:
\begin{equation}
\mathbb E_\pi[\randombold{\zeta}(\beta)]
=\beta/2+\log\cosh(\beta/4)
+\frac12\log\frac{4\sinh^3(\beta/2)}{3(\sinh\beta-\beta)}.
\label{eq:quenched-bagchi-pal-mean-cgf}
\end{equation}
At zero it has the expansion
\begin{equation}
\mathbb E_\pi[\randombold{\zeta}(\beta)]
=\beta/2+11\beta^2/160-439\beta^4/537600+O(\beta^6).
\label{eq:quenched-bagchi-pal-mean-expansion}
\end{equation}
The derivative at zero records the common limiting share. The
curvature $11/80$ averages fluctuations at a fixed matrix, even
though both matrices have the same equilibrium. These are CGFs of
the extraction count. For the physical black count one instead
averages the right--hand side of
Eq.~(\ref{eq:app-bagchi-pal-physical-cgf}), using the random matrix
entries and $\randombold{\zeta}$; the normalization remains per HLS step.

A continuous law also allows the limiting share to vary. Choose
both coefficients uniformly over the full admissible domain:
\begin{equation}
0<\randombold{a}<1,\qquad -\randombold{a}<\randombold{b}<1-\randombold{a}.
\label{eq:quenched-affine-domain}
\end{equation}
The coefficients are dependent; the parallelogram has area one,
so their joint uniform density is one.
It includes both increasing and decreasing responses; the line
$\randombold{b}=0$ has probability zero. Equivalently, the endpoint
values $\randombold{a}$ and $\randombold{a}+\randombold{b}$ are independent
and uniform on $(0,1)$. This is an exact binary HLS model. Its
replacement matrix interpretation uses nonnegative real masses;
a continuous coefficient law is not a law on integer numbers of balls.

For each fixed coefficient pair $a,b$, let $\zeta(\beta)$ be the
CGF given by Eqs.~(\ref{eq:app-affine-cgf-positive})--(\ref{eq:app-affine-cgf-decreasing}),
with Eq.~(\ref{eq:app-bernoulli-cgf}) on $b=0$.
Integrating these fixed--coefficient CGFs gives the quenched mean:
\begin{equation}
\mathbb E_\pi[\randombold{\zeta}(\beta)]
=\int_0^1 da\int_{-a}^{1-a}db\,\zeta(\beta).
\label{eq:quenched-affine-uniform-mean}
\end{equation}
\afterpage{%
\clearpage
\begin{figure}[p]
\centering
\makebox[\linewidth][l]{\hspace*{-9.244983pt}\includegraphics[width=364.541528pt]{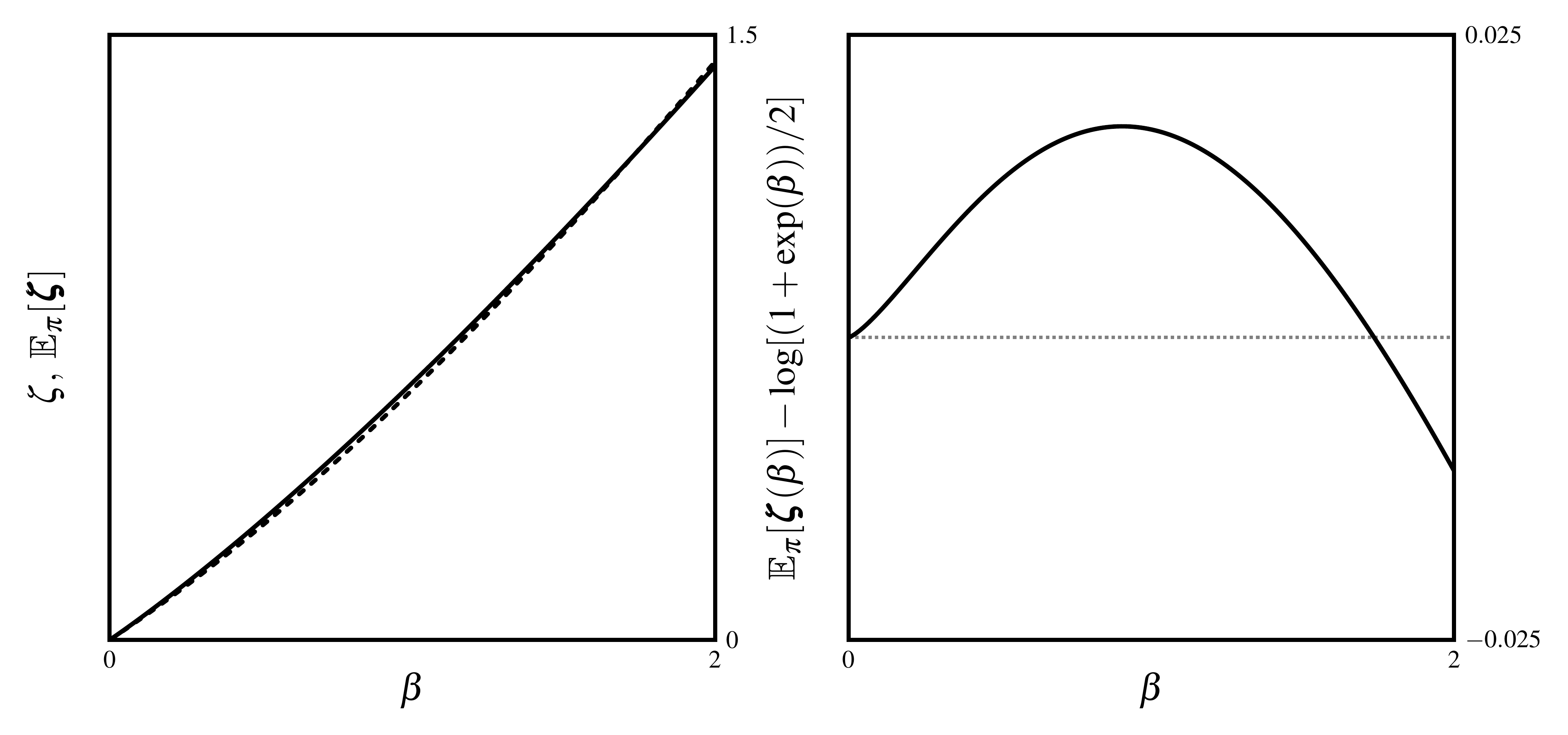}}
\caption{CGFs for affine urns under the uniform law on the full domain of Eq.~(\ref{eq:quenched-affine-domain}). Left: the quenched mean $\mathbb E_\pi[\randombold{\zeta}(\beta)]$ (solid) and the Bernoulli CGF $\log[(1+\exp(\beta))/2]$ of the urn with mean coefficients $a=1/2$, $b=0$ (dashed). Right: the former minus the latter, with a dotted zero line. Both panels span $0\leq\beta\leq2$; vertical ranges are $[0,1.5]$ and $[-0.025,0.025]$, respectively. The difference changes sign: averaging the solved CGFs and inserting the mean coefficients give distinct free energies, with neither curve uniformly above the other on this interval.}
\label{fig:quenched-affine-cgf}
\end{figure}
\clearpage
}
Figure~\ref{fig:quenched-affine-cgf}\label{call:fig:quenched-affine-cgf}
compares this quenched mean with the CGF of the urn with mean
coefficients and shows their difference.

Numerical quadrature over the full domain gives
\begin{equation}
\mathbb E_\pi[\randombold{\zeta}(1)]\simeq0.63729035,
\qquad \mathbb E_\pi[\randombold{\zeta}(2)]\simeq1.42274871.
\label{eq:quenched-affine-uniform-values}
\end{equation}

Color exchange preserves this coefficient law, hence
\begin{equation}
\mathbb E_\pi[\randombold{\zeta}(-\beta)]
=\mathbb E_\pi[\randombold{\zeta}(\beta)]-\beta.
\label{eq:quenched-affine-uniform-symmetry}
\end{equation}
The averaged CGF is differentiable at zero, with derivative $1/2$,
but has no finite quadratic coefficient:
\begin{equation}
\lim_{\beta\to0}
\frac{\mathbb E_\pi[\randombold{\zeta}(\beta)]-\beta/2}{\beta^2}=+\infty.
\label{eq:quenched-affine-uniform-curvature}
\end{equation}
Indeed, the conditional first derivative is
\begin{equation}
\partial_\beta\randombold{\zeta}(0)=\randombold{a}/(1-\randombold{b}),
\label{eq:quenched-affine-first-derivative}
\end{equation}
whose mean is $1/2$. The conditional CGFs have slopes in $[0,1]$,
so dominated convergence justifies this first derivative of the average.
For $\randombold{b}<1/2$, the conditional curvature is
\begin{equation}
\partial_\beta^2\randombold{\zeta}(0)
=\randombold{a}(1-\randombold{b}-\randombold{a})/[(1-\randombold{b})^2(1-2\randombold{b})].
\label{eq:quenched-affine-second-derivative}
\end{equation}
Its integral over $0<\randombold{b}<1/2$ already diverges logarithmically.
Convexity and Fatou's lemma then give
Eq.~(\ref{eq:quenched-affine-uniform-curvature}). In contrast, the mean response is
\begin{equation}
\mathbb E_\pi[\randombold{\pi}(\psi)]=1/2.
\label{eq:quenched-affine-mean-urn}
\end{equation}
Its Bernoulli CGF, given by Eq.~(\ref{eq:app-bernoulli-cgf}), has finite
curvature $1/4$. Thus replacing the random function by its mean
changes the free energy.

The equilibrium $\randombold{\psi}_*=\randombold{a}/(1-\randombold{b})$
has mean $1/2$ and second moment $1-\log2$ under this coefficient law.
The variance of the stochastic share over coefficients and growth therefore satisfies
\begin{equation}
\lim_{n\to\infty}\operatorname{Var}(\randombold{\psi}_n)=3/4-\log2.
\label{eq:quenched-affine-uniform-variance}
\end{equation}
This persistent variation across urns is distinct from the
conditional fluctuations measured by the curvature of the quenched
CGF. Increasing the population does not remove the spread of the
random equilibria.

\subsubsection{Asymmetric autocatalysis with a frozen bias}\label{app:quenched-chiral}

The illustrative autocatalytic rule in Eq.~(\ref{eq:bio-chiral}) can represent
persistent differences between reaction environments by drawing the
bias $\randombold{h}$ once for each reactor. Keeping $a>0$ and $q>1$ fixed gives
\begin{equation}
\randombold{\pi}(\psi)
=\exp(\randombold{h})(a+\psi^q)/\{\exp(\randombold{h})(a+\psi^q)+\exp(-\randombold{h})[a+(1-\psi)^q]\}.
\label{eq:quenched-chiral-urn}
\end{equation}
At fixed $h$, write $\pi$ for the realized urn function and $\psi_*(h)$
for a fixed--point branch. On its domain,
$\randombold{\psi}_*=\psi_*(\randombold{h})$ is the random equilibrium.
For every finite field the urn function is continuously differentiable,
strictly increasing and strictly between zero and one.

The fixed points and their basins are determined separately for each
field. Along a differentiable branch of fixed points with
$\partial_\psi\pi(\psi_*)\ne1$, implicit differentiation gives
\begin{equation}
\frac{d\psi_*}{dh}
=\frac{2\psi_*(1-\psi_*)}{1-\partial_\psi\pi(\psi_*)}.
\label{eq:quenched-chiral-sensitivity}
\end{equation}
Thus a stable branch becomes increasingly sensitive to the frozen
bias as the urn-function slope $\partial_\psi\pi(\psi_*)$ approaches one. The small--field response can be
computed even at the degenerate symmetric equilibrium. For $1<q<5$, on the
critical line $a=(q-1)/2^q$, the stable branch approaching $1/2$ satisfies
\begin{equation}
\psi_*(h)-1/2\sim\operatorname{sgn}(h)
\{3|h|/[4(q-1)(5-q)]\}^{1/3}.
\label{eq:quenched-chiral-critical}
\end{equation}
At the tricritical point, where the linear and cubic terms in the
expansion of $\pi(\psi)-\psi$ about $\psi=1/2$ vanish at $h=0$,
the response is stronger:
\begin{equation}
\psi_*(h)-1/2\sim\operatorname{sgn}(h)(5|h|/128)^{1/5},
\qquad q=5,\qquad a=1/8.
\label{eq:quenched-chiral-tricritical}
\end{equation}
Both limits are taken as $h\to0$. They follow from the cubic and
quintic terms of the zero--field drift, respectively, together with
$\partial_h\pi(1/2)=1/2$ at $h=0$. These responses show how a small spread of frozen biases produces a
much broader spread of equilibrium compositions in the proposed model. They are statements
about fixed points, without assuming a joint limit of field and
population size.

For
$h>0$ and a prescribed macroscopic seed, an example is the terminal
event
\begin{equation}
Q(E^*)=\{\varphi\in Q:\varphi(\tau_0)=\tau_0\psi_0,
\ \varphi(1)\leq1/2\},\qquad 0<\tau_0<1.
\label{eq:quenched-chiral-event}
\end{equation}
Equation~(\ref{eq:quenched-event}) gives its logarithmic probability
using Eq.~(\ref{eq:bio-chiral-action}) for
$\Phi$ at the fixed field, with both action integrals restricted to
$[\tau_0,1]$. The seed may be chosen in the favored basin to examine
histories ending with the opposite majority. Reachability and event
regularity must still hold; a positive cost requires exclusion of
zero--cost histories from the event. With a fixed finite seed, two
accessible attractors can instead both have zero exponential cost.

Figure~\ref{fig:quenched-chiral}\label{call:fig:quenched-chiral} compares equilibrium responses under the same frozen--bias law: the tricritical response amplifies the bias more strongly. For this illustration, $\randombold{h}$ is uniform on $[-1/50,1/50]$, with $(q,a)=(2,1/4)$ or $(5,1/8)$ fixed. The law is independent of $T$ and is sampled once for each reactor.

\afterpage{%
\clearpage
\begin{figure}[p]
\centering
\makebox[\linewidth][l]{\hspace*{79.480569pt}\includegraphics[width=174.997460pt]{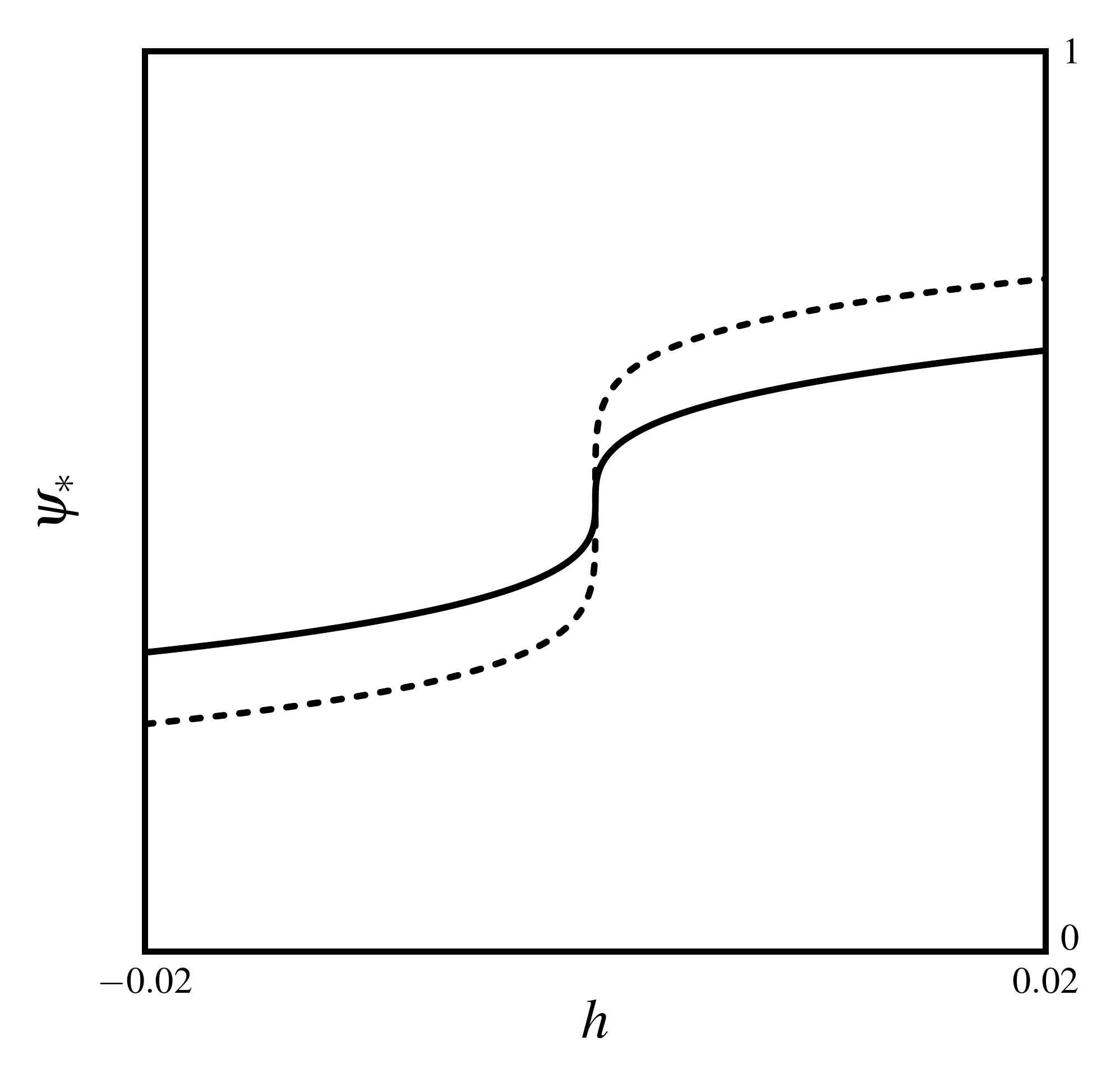}}
\caption{Equilibrium composition $\psi_*(h)$ for the critical autocatalytic response $(q,a)=(2,1/4)$ (solid) and the tricritical response $(q,a)=(5,1/8)$ (dashed). The horizontal field range is $[-1/50,1/50]$ and the vertical composition range is $[0,1]$. A uniform frozen bias on this interval produces a wider spread of equilibria in the tricritical case, according to the $1/3$ and $1/5$ laws in Eqs.~(\ref{eq:quenched-chiral-critical}) and (\ref{eq:quenched-chiral-tricritical}). Each reactor retains its sampled field.}
\label{fig:quenched-chiral}
\end{figure}
\clearpage
}

\subsubsection{Magnetic growth in a frozen field}\label{app:quenched-magnetic}

In the irreversible growth model of Section~\ref{app:magnetic}, draw
a finite field $\randombold{h}$ once for each sample and keep $J>0$ fixed:
\begin{equation}
\randombold{\pi}(\psi)
=\{1+\tanh[J(2\psi-1)+\randombold{h}]\}/2.
\label{eq:quenched-magnetic-urn}
\end{equation}
The same fixed--field notation applies, branch by branch.
The field is common to all additions in that sample. Temperature
remains absorbed into $J,h$, as in Eq.~(\ref{eq:magneticurn}).
This gives a frozen--field extension of the growth process studied
in \cite{KGW,Jack=000020LD,KGGW}.

The equilibrium sensitivity is again
Eq.~(\ref{eq:quenched-chiral-sensitivity}), now with
$\partial_\psi\pi(\psi_*)=4J\psi_*(1-\psi_*)$.
For $J<1$, the unique equilibrium has the small--field response
\begin{equation}
\psi_*(h)-1/2\sim h/[2(1-J)],\qquad h\to0.
\label{eq:quenched-magnetic-subcritical}
\end{equation}
At $J=1$, the leading response becomes
\begin{equation}
\psi_*(h)-1/2\sim
\operatorname{sgn}(h)(3|h|)^{1/3}/2,\qquad h\to0.
\label{eq:quenched-magnetic-critical}
\end{equation}
Consequently, weak variation of the field is strongly amplified
near the ordering threshold. For $J>1$, sufficiently small fields
leave two attractors; the frozen field changes their positions and
basins, while the growth history still affects which attractor is approached.
On each stable branch, Eq.~(\ref{eq:magnetic_relaxation}) describes
the local deterministic relaxation with the realized field and
fixed point.

Replacing $h$ by the random field $\randombold{h}$ in
Eqs.~(\ref{eq:magnetic_lagrangian})--(\ref{eq:magnetic_action}) gives
$\randombold{\Phi}$. Together with Eq.~(\ref{eq:quenched-event}),
this determines the cost and optimal histories of a prescribed
terminal magnetization, with the same seed distinction as in the
chiral example.

Figure~\ref{fig:quenched-magnetic}\label{call:fig:quenched-magnetic} shows how a frozen field selects different basins from the same starting composition. Here $\randombold{h}$ is uniform on $[-3/25,3/25]$, with $J=3/2$ fixed. The response is admissible for every finite field and positive $J$, and the entire chosen field interval lies within the regime with two attracting branches. The law is independent of $T$ and is sampled once per growing sample.

\afterpage{%
\clearpage
\begin{figure}[p]
\centering
\includegraphics[width=\linewidth]{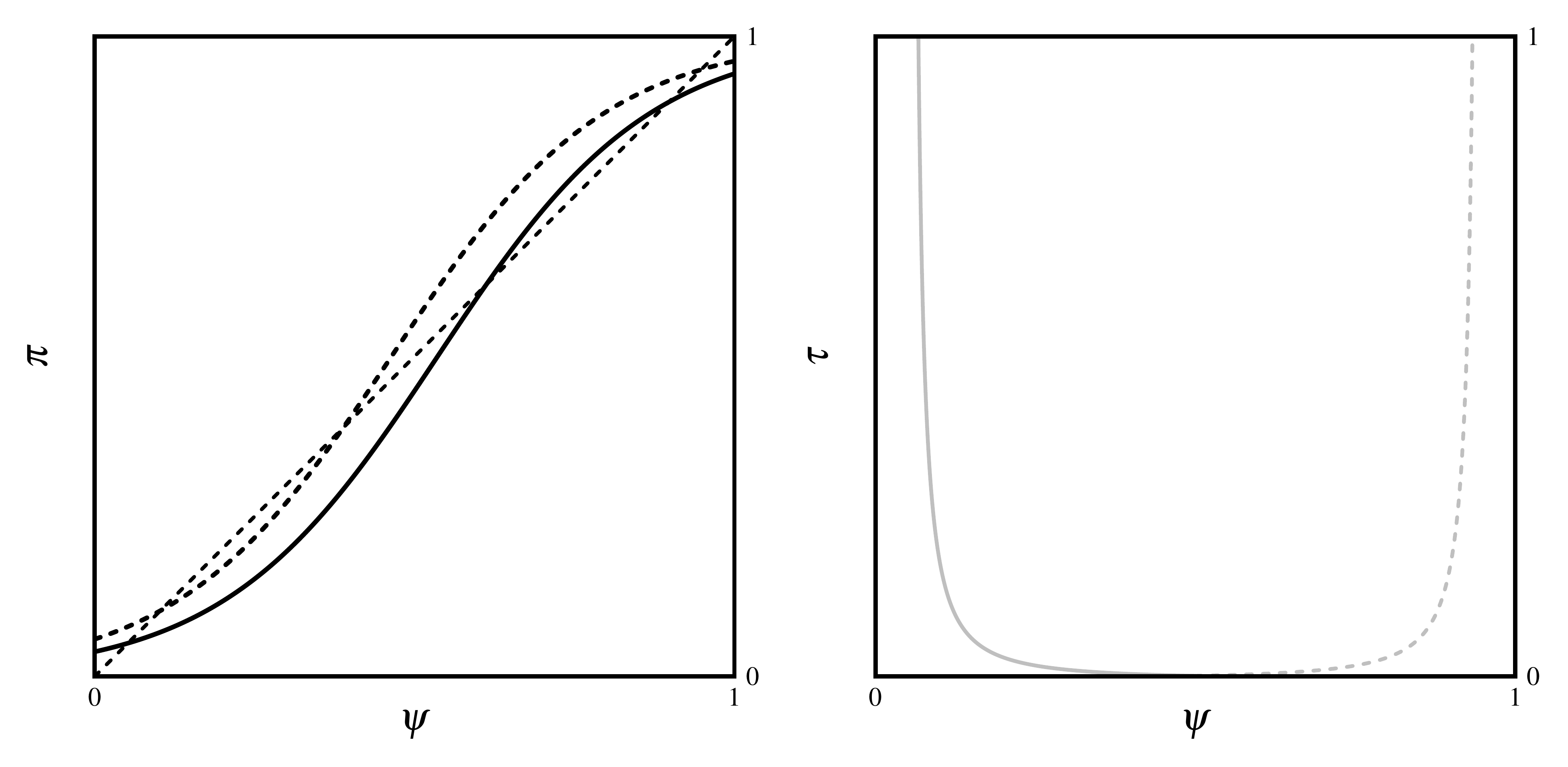}
\caption{Urn functions (left) and zero--cost share trajectories (right) for magnetic growth at $J=3/2$, with $h=-0.108$ (solid) and $h=0.108$ (dashed). These are the $5\%$ and $95\%$ quantiles of the uniform law on $[-3/25,3/25]$. Gray trajectories retain these styles and start from the common seed $\psi_0=1/2$ at $\tau_0=0.001$. Opposite frozen fields select opposite attracting branches. All axes span $[0,1]$, with saturation increasing upwards. The shaded strip marks $0\leq\tau<\tau_0$. The diagonal is dashed.}
\label{fig:quenched-magnetic}
\end{figure}
\clearpage
}

\subsubsection{Cell--type proportions across replicas}\label{app:quenched-biological}

The controller of Eq.~(\ref{eq:bio-controller}) can also represent
differences between growing biological replicas. Choose the target
$\randombold{p}$ and feedback strength $\randombold{\kappa}$ once for each replica, with
$0<\randombold{p}<1$ and finite $\randombold{\kappa}\geq0$ almost surely, and retain them
throughout its sequence of commitments. The random response is
\begin{equation}
\randombold{\pi}(\psi)
=\{1+[(1-\randombold{p})/\randombold{p}]\exp[\randombold{\kappa}(\psi-\randombold{p})]\}^{-1}.
\label{eq:quenched-biological-urn}
\end{equation}
The response and its variation between replicas are proposed effective
descriptions; their parameters must be related to the system being studied. No particular
distribution of $\randombold{p}$ or $\randombold{\kappa}$ is needed for the following identities.
At a fixed realization, $p,\kappa$ denote the selected coefficients
and $\pi$ the realized urn function.

Every realization has the unique attracting composition $p$, and
\begin{equation}
\partial_\psi\pi(p)=-\kappa p(1-p).
\label{eq:quenched-biological-slope}
\end{equation}
The response is strictly decreasing for positive $\kappa$ and
constant for $\kappa=0$. At fixed parameters, a small perturbation
of the share relaxes along the deterministic trajectory as
\begin{equation}
\psi(\tau)-p\simeq(\psi_0-p)
(\tau_0/\tau)^{1+\kappa p(1-p)},
\qquad 0<\tau_0\leq\tau.
\label{eq:quenched-biological-relaxation}
\end{equation}
The target controls the limiting proportion, while the strength of
the feedback also controls the recovery after a perturbation.
Growth still consists of new commitments, as in the original example.

\afterpage{%
\clearpage
\begin{figure}[p]
\centering
\includegraphics[width=\linewidth]{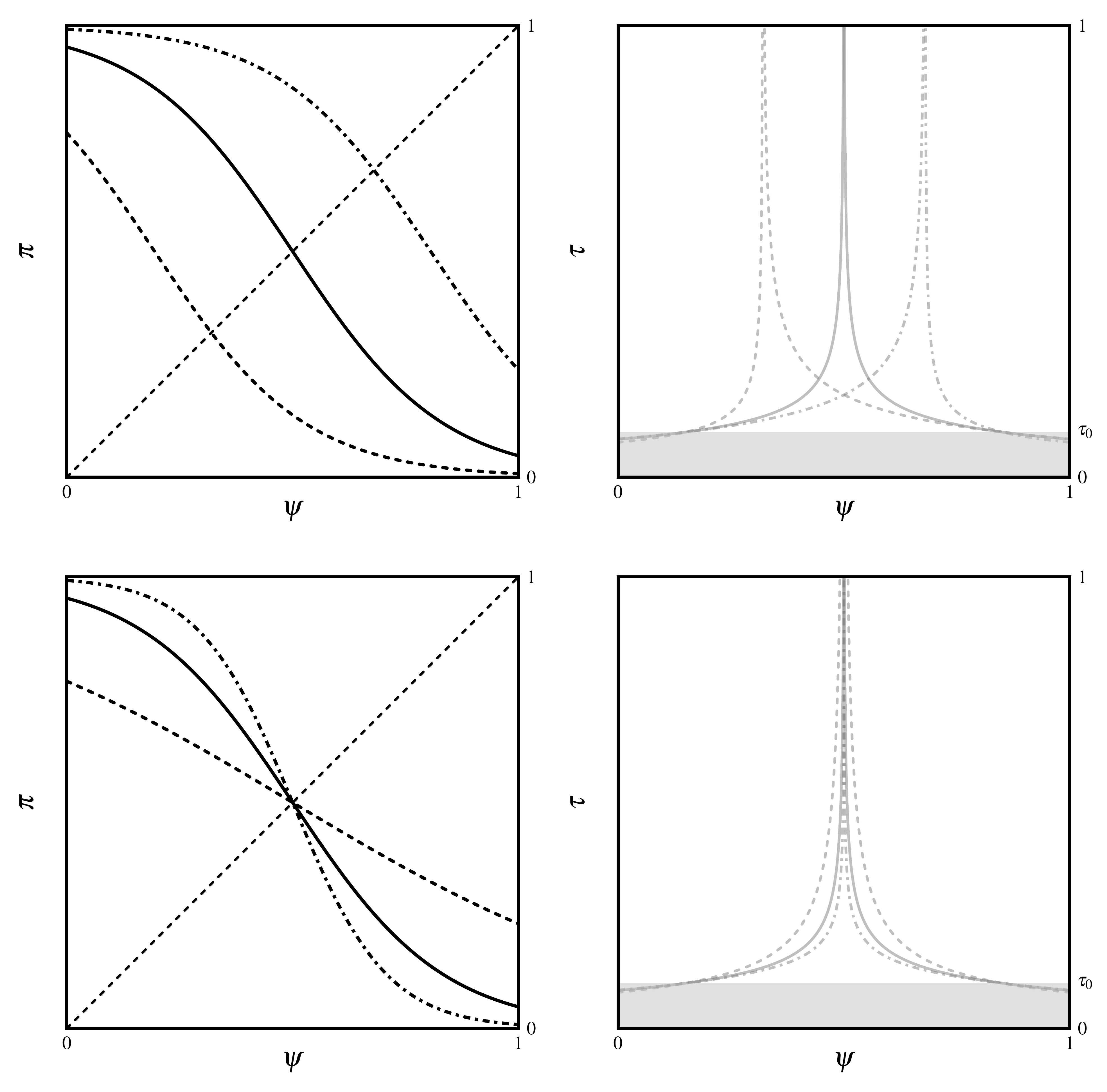}
\caption{Urn functions (left) and zero--cost share trajectories (right) for cell--type control. Top: $\kappa=6$, with $p=0.32,0.5,0.68$, the $5\%$, $50\%$ and $95\%$ quantiles of a uniform law on $[0.3,0.7]$. Bottom: $p=1/2$, with $\kappa=2.4,6,9.6$, the same quantiles of a uniform law on $[2,10]$. These quantiles use dashed, solid and dash--dotted curves, respectively. Gray trajectories start from $\psi_0=0.15,0.85$ at $\tau_0=0.1$. Axes cover $[0,1]$, with saturation increasing upwards. Random targets spread the limiting compositions; random feedback preserves their common limit. The shaded strip marks $0\leq\tau<\tau_0$ and shows backward continuations, stopped at the admissible boundary.}
\label{fig:quenched-cell}
\end{figure}
\clearpage
}

Figure~\ref{fig:quenched-cell}\label{call:fig:quenched-cell} contrasts two forms of heterogeneity. The first illustrative law takes $\randombold{p}$ uniform on $[0.3,0.7]$ and fixes $\kappa=6$; the second fixes $p=1/2$ and takes $\randombold{\kappa}$ uniform on $[2,10]$. Both lie within the stated admissible parameter range and are sampled once per replica, independently of $T$. The resulting limiting variances are $1/75$ and zero, respectively.

For the stochastic process, $\randombold{\psi}_n$ converges almost surely
to the random target $\randombold{p}$ of its replica.
Equation~(\ref{eq:quenched-limit-variance})
therefore specializes to
\begin{equation}
\lim_{n\to\infty}\operatorname{Var}(\randombold{\psi}_n)
=\operatorname{Var}(\randombold{p}).
\label{eq:quenched-biological-variance}
\end{equation}
The share variance is over parameters and growth; the target variance
is over the sampling law of $\randombold{p}$.
Thus increasing the number of commitments does not remove differences
in the proportions specified by different targets. If all targets
coincide and only $\randombold{\kappa}$ varies, the limiting proportion is shared,
but recovery and finite--size fluctuations can still differ.
At finite size, Eq.~(\ref{eq:quenched-total-variance}) remains the
exact decomposition. A conditional fluctuation estimate alone does
not determine the next term in the variance across replicas: the
conditional mean can also retain a finite--size bias.

For a fixed realization and a tolerance
$0<\epsilon<\min(p,1-p)$, the event
in Eq.~(\ref{eq:bio-control-event}) has its conditional cost given
by Eq.~(\ref{eq:quenched-event}), with
Eq.~(\ref{eq:quenched-biological-urn}) in the action. This quantifies
failure to reach the replica's own target; deviation from a common
externally assigned target is a different endpoint constraint.

More generally, the same distinction between fluctuations within
a replica and persistent differences between replicas is relevant
to circuits, embryos or reactors whose growth probabilities depend on the system's history only
through its current share, for each fixed choice of parameters. On the shares visited by a nonconstant scaling trajectory, the urn
function can be reconstructed using Section~\ref{sec:inverse}. A
stationary trajectory supplies only its fixed point. Comparing
replicas can thus reveal differences in equilibria and feedback
strengths; inference from noisy observations remains a statistical task.

\subsection{Replica Symmetry Breaking in the HLS model}\label{app:rsb}

We conclude by discussing the overlap properties of the HLS model
and the pure-state interpretation of its limiting trajectories.
At positive initial saturation, a prescribed initial share and urn
function determine, under the assumptions below, a trajectory with
definite local spin means and overlaps. The asymptotic factorization
of any fixed number of increments supports interpreting the
corresponding conditional component as a pure state specified by
local magnetizations \cite{MezardParisiVirasoro1987,KERNEL=000020THEO}.
In the spin-kernel representation and the Blanket construction
\cite{Franchini2024,RSBSPA3v2}, these magnetizations describe the
components over which the measure is disintegrated. Here they are obtained directly from temporal
sectors of the urn trajectory. We first derive these properties from
the urn dynamics, then compare a prescribed initial share, a uniform
share frozen in each system, and a uniform share integrated within
the path measure. Repeating these three preparations for a quenched
random urn function gives six comparisons of overlap laws and free
energies. The construction also applies to nonlinear responses under
the stated assumptions; affine urns provide examples with simple
expressions. The comparison with these representations follows at the end.

\subsubsection{Trajectory concentration and overlaps}\label{app:rsb-trajectories}

Fix the initial share $\psi_0$ and saturation $0<\tau_0<1$.
At finite capacity, the initial history occupies the following sector
of the index set $S$:
\begin{equation}
S_0=\{n\in S:n\leq\tau_0T\},
\qquad |S_0|=\lfloor\tau_0T\rfloor.
\label{eq:rsb-initial-sector}
\end{equation}
Prescribe its black count as the integer part of $\psi_0|S_0|$, so that
its share tends to $\psi_0$ as $T\to\infty$. Assume that $\pi$ is
Lipschitz on $[0,1]$, and let $\psi(\tau)$ be the unique solution of
Eq.~(\ref{eq:Caucy}) with this initial condition. The binary increments
remain $\sigma_n\in\{0,1\}$, and their spins are $2\sigma_n-1$.
The count trajectory and the local spin mean satisfy
\begin{equation}
\varphi(\tau)=\tau\psi(\tau),
\qquad
2\pi(\psi(\tau))-1
=2\partial_\tau\varphi(\tau)-1.
\label{eq:rsb-local-profile}
\end{equation}
The mean $2\pi(\psi(\tau))-1$ concerns newly
added spins; $2\psi(\tau)-1$ describes the accumulated spin mean.
The bounded martingale increments of the HLS process, followed by
Gronwall's inequality on $[\tau_0,1]$, give
\begin{equation}
\max_{|S_0|<n\leq T}
\bigl|\psi_n-\psi(n/T)\bigr|
\xrightarrow{\mathbb P}0.
\label{eq:rsb-trajectory-concentration}
\end{equation}
Continuity of $\pi$ and uniqueness of the prescribed trajectory also
suffice, as in Section~\ref{sec:trajectories}; the Lipschitz assumption
is a convenient sufficient condition. Uniqueness is not automatic for
every H\"older response.

Concentration also determines the joint law of any fixed number of
future increments. Condition on the history before each selected
increment, starting with the latest, and replace its conditional
probability by $\pi(\psi(n/T))$. Each replacement error tends to zero
by Eq.~(\ref{eq:rsb-trajectory-concentration}), since the transition
probabilities converge uniformly and are bounded. The joint law thus
approaches the product of the corresponding Bernoulli laws. This
factorization applies to a fixed number of increments; it does not
replace the probability law of the whole history by a product.

To compare the local spin means along the histories, divide the future
sites into $\ell_{\max}$ consecutive sectors. Their endpoints satisfy
\begin{equation}
\tau_0<\tau_1<\cdots<\tau_{\ell_{\max}}=1.
\label{eq:rsb-partition-endpoints}
\end{equation}
The sectors are
\begin{equation}
S_\ell=\{n\in S:T\tau_{\ell-1}<n\leq T\tau_\ell\},
\qquad 1\leq\ell\leq\ell_{\max}.
\label{eq:rsb-sectors}
\end{equation}
In each sector, the average spin converges to the time average of $2\pi(\psi(\tau))-1$. These sector magnetizations are
therefore determined by the increments of the count trajectory:
\begin{equation}
\frac{1}{|S_\ell|}\sum\limits_{n\in S_\ell}(2\sigma_n-1)
\xrightarrow{\mathbb P}
\frac{2[\varphi(\tau_\ell)-\varphi(\tau_{\ell-1})]}
{\tau_\ell-\tau_{\ell-1}}-1.
\label{eq:rsb-sector-magnetizations}
\end{equation}

Consider two independent continuations, denoted by $\sigma$ and
$\sigma'$, from initial shares $\psi_0$ and $\psi_0'$ under the same
urn function. A prime distinguishes the second copy: its trajectories are
$\varphi',\psi'$,
with $\psi'(\tau_0)=\psi_0'$.
Normalize the overlap by the number of future sites:
\begin{equation}
q_T:=\frac{1}{T-|S_0|}\sum\limits_{|S_0|<n\leq T}
(2\sigma_n-1)(2\sigma_n'-1).
\label{eq:rsb-overlap-definition}
\end{equation}
Given the two past histories, the next draws are independent, so the
conditional mean of their spin product is the product of the two
conditional spin means. Subtracting the conditional means of the spin
products leaves a martingale whose contribution to $q_T$ has variance
at most $1/(T-|S_0|)$. Equation~(\ref{eq:rsb-trajectory-concentration})
then allows the conditional means to be replaced by their deterministic
profiles. The resulting Riemann sum shows that $q_T$ converges in
probability to
\begin{equation}
q_\pi(\psi_0,\psi_0')=
\frac{1}{1-\tau_0}\int_{\tau_0}^{1}
[2\pi(\psi(\tau))-1][2\pi(\psi'(\tau))-1]\,d\tau.
\label{eq:rsb-overlap-limit}
\end{equation}
After taking $T\to\infty$, refine the partition to recover the
scalar product from the sector magnetizations:
\begin{multline}
q_\pi(\psi_0,\psi_0')=
\lim_{\max_\ell(\tau_\ell-\tau_{\ell-1})\to0}
\sum\limits_{1\leq\ell\leq\ell_{\max}}
\frac{\tau_\ell-\tau_{\ell-1}}{1-\tau_0}\\
\times\left[\frac{2[\varphi(\tau_\ell)-\varphi(\tau_{\ell-1})]}
{\tau_\ell-\tau_{\ell-1}}-1\right]
\left[\frac{2[\varphi'(\tau_\ell)-\varphi'(\tau_{\ell-1})]}
{\tau_\ell-\tau_{\ell-1}}-1\right].
\label{eq:rsb-sector-overlap}
\end{multline}
At a fixed partition, the product of two sector means need not equal the sector
mean of their product.

The initial sector must be included separately for an overlap of the
whole histories. If its normalized overlap converges to $q_0$, then
\begin{equation}
\frac{1}{T}\sum\limits_{1\leq n\leq T}
(2\sigma_n-1)(2\sigma_n'-1)
\xrightarrow{\mathbb P}
\tau_0q_0+(1-\tau_0)q_\pi(\psi_0,\psi_0').
\label{eq:rsb-full-overlap}
\end{equation}
A shared microscopic prefix has $\psi_0=\psi_0'$ and $q_0=1$. Independently sampled
uniform prefixes with prescribed shares $\psi_0,\psi_0'$ instead give
\begin{equation}
q_0=(2\psi_0-1)(2\psi_0'-1).
\label{eq:rsb-independent-prefix-overlap}
\end{equation}
For two assigned different prefixes, their actual
overlap is required; the shares alone do not determine it.
We use the future overlap of Eq.~(\ref{eq:rsb-overlap-definition}) below.
Its conditional concentration holds for every pair of prescribed shares;
a mixture can retain a nontrivial overlap distribution when it assigns
weight to profiles with different scalar products.

\subsubsection{Random initial conditions}
\label{app:rsb-initial}

Keep the urn function fixed. For a prescribed initial share $\psi_0$,
two independent continuations have limiting overlap
$q_\pi(\psi_0,\psi_0)$. Let $\delta_x$ denote unit mass at $x$ and
write $P(q)\,dq$ for this limiting law, leaving its dependence on
the prescribed urn function and initial share implicit. Here $P$ is a
generalized density, allowing Dirac masses: its associated measure is
the unit mass at $q_\pi(\psi_0,\psi_0)$. Averages of $P$ mean averages
of these measures; where the resulting law has an ordinary density,
the same notation denotes that density. The symbol $P$ denotes this
conditional law as a function of its inputs; each averaging subscript
specifies which input is integrated.

To compare different initial conditions, give the sampled share
$\randombold{\psi}_0$ the uniform law $\mathbb P_{\psi_0}$ on $[0,1]$,
with expectation $\mathbb E_{\psi_0}$. Boldface marks a sampled input,
or a quantity depending on it, before it is fixed; $\psi_0$ denotes
a prescribed value or the variable of integration.
At finite capacity, choosing the initial count uniformly among
$0,\ldots,|S_0|$ gives the same limiting share law. Giving equal
probabilities to all binary prefixes would instead concentrate the
share at $1/2$. If the sampled share is frozen, the overlap remains
concentrated within each system. Averaging over the sampled initial
conditions gives the limiting overlap law
\begin{equation}
\mathbb E_{\psi_0}[P(q)]\,dq=\int_0^1d\psi_0\,\delta_{q_\pi(\psi_0,\psi_0)}(dq).
\label{eq:rsb-diagonal-law}
\end{equation}
If the initial share is instead integrated within the path measure,
two independent configurations draw their shares independently from
the same law. Denote the second draw by $\randombold{\psi}_0'$;
$\mathbb E_{\psi_0'}$ averages this independent draw under the
same law $\mathbb P_{\psi_0}$. Write $P_{\mathrm{ann}}(q)$ for the
resulting generalized overlap density. The subscript refers only to
this initial-share mixture, whose limiting law is
\begin{equation}
P_{\mathrm{ann}}(q)\,dq=\int_0^1d\psi_0\int_0^1d\psi_0'\,
\delta_{q_\pi(\psi_0,\psi_0')}(dq).
\label{eq:rsb-mixture-law}
\end{equation}
These formulas apply to
nonlinear responses as well, whenever the concentration assumptions
above hold. A nonconstant
response and a uniform initial share give distinct local profiles on
sets of positive measure, so the mixed overlap law is non-concentrated.
Indeed, the squared distance between two profiles, normalized by the
duration of the future interval, is
\begin{equation}
q_\pi(\psi_0,\psi_0)+q_\pi(\psi_0',\psi_0')
-2q_\pi(\psi_0,\psi_0').
\label{eq:rsb-profile-distance}
\end{equation}
Continuity then prevents the overlap kernel from being constant almost
everywhere. For a constant response, the
future increment profile is independent of the initial share.

To illustrate both laws with simple expressions, take the affine
response $\pi(\psi)=a+b\psi$ on $[0,1]$, write $\psi_*=a/(1-b)$
when $b<1$ and use the trajectory
in Eq.~(\ref{eq:app-affine-trajectory}). Its local spin mean is
\begin{equation}
2\pi(\psi(\tau))-1
=2\psi_*-1+b(2\psi_0-2\psi_*)(\tau_0/\tau)^{1-b}.
\label{eq:rsb-affine-local-mean}
\end{equation}
Integrating the two affine profiles in
Eq.~(\ref{eq:rsb-overlap-limit}) gives the quadratic expression
\begin{multline}
q_\pi(\psi_0,\psi_0')=(2\psi_*-1)^2
+2(2\psi_*-1)(\psi_0+\psi_0'-2\psi_*)
\frac{\tau_0^{1-b}-\tau_0}{1-\tau_0}\\
+\frac{4b^2(\psi_0-\psi_*)(\psi_0'-\psi_*)}{1-\tau_0}
\int_{\tau_0}^1(\tau_0/\tau)^{2(1-b)}d\tau.
\label{eq:rsb-affine-kernel}
\end{multline}
Thus the diagonal law is obtained by setting $\psi_0'=\psi_0$ and
averaging a uniform initial share. The mixed law instead averages
Eq.~(\ref{eq:rsb-affine-kernel}) over two independent uniform shares.

For the symmetric family $a=(1-b)/2$, the overlap $q_\pi(1,1)$
of two continuations from an all-one initial population determines
the full kernel:
\begin{equation}
q_\pi(\psi_0,\psi_0')
=q_\pi(1,1)(2\psi_0-1)(2\psi_0'-1).
\label{eq:rsb-symmetric-kernel}
\end{equation}
Its value is
\begin{equation}
q_\pi(1,1)=\frac{b^2(\tau_0^{2-2b}-\tau_0)}
{(2b-1)(1-\tau_0)},\qquad b\ne1/2.
\label{eq:rsb-affine-C}
\end{equation}
At $b=1/2$, its continuous value is
\begin{equation}
q_\pi(1,1)=\frac{\tau_0\log(1/\tau_0)}{4(1-\tau_0)}.
\label{eq:rsb-affine-C-half}
\end{equation}
For $b\ne0$, both laws have densities with respect to $dq$.
For a uniformly sampled share frozen within each system, the density is
\begin{equation}
\mathbb E_{\psi_0}[P(q)]=\frac{1}{2\sqrt{q_\pi(1,1)q}},
\qquad 0<q<q_\pi(1,1).
\label{eq:rsb-symmetric-diagonal-density}
\end{equation}
For the initial share integrated within the path measure, it is
\begin{equation}
P_{\mathrm{ann}}(q)=\frac{\log[q_\pi(1,1)/|q|]}{2q_\pi(1,1)},
\qquad 0<|q|<q_\pi(1,1).
\label{eq:rsb-symmetric-mixed-density}
\end{equation}
They vanish outside their indicated supports. Their singularities at
zero are integrable and do not represent atoms. The symmetric affine
urn in Figure~\ref{fig:lft-two-body}(a), with $a=1/4$ and $b=1/2$,
provides an example of both distributions: what changes is how the
initial share enters the measure. For P\'olya's response,
$b=1$ and $a=0$, the share profile is constant and the same formulas
hold with $q_\pi(1,1)=1$. For $b=0$, both laws instead reduce to unit mass at
$(2a-1)^2$. Throughout, the initial population is macroscopic.
Taking $T$ to infinity first and then $\tau_0$ to zero removes the
dependence on the initial share for every fixed affine response with
$b<1$. In the symmetric family, both overlap laws then concentrate at
zero, whereas P\'olya's response retains $q_\pi(1,1)=1$.
This loss of initial--share dependence does not remove the variation
of equilibria under the random--coefficient law of
Section~\ref{app:quenched-affine}.

For the symmetric affine response, averaging the trajectory over a
uniform initial share gives
\begin{equation}
\mathbb E_{\psi_0}[\randombold{\psi}(\tau)]=\int_0^1\psi(\tau)\,d\psi_0=1/2,
\qquad \tau_0\leq\tau\leq1.
\label{eq:rsb-uniform-initial-mean-trajectory}
\end{equation}
This mean is the same whether the share is integrated within the path
measure or frozen and averaged afterwards. The overlap laws distinguish
these preparations: two histories choose their initial shares independently
in the first case and share the same frozen value in the second.
For the illustrated response, the frozen-share average has positive mean
overlap even though the mean local spin profile vanishes. Averaging
trajectories therefore does not replace averaging their overlaps.

\afterpage{%
\clearpage
\begin{figure}[p]
\centering
\makebox[\linewidth][l]{\hspace*{-4.600584pt}\includegraphics[width=344.999999pt]{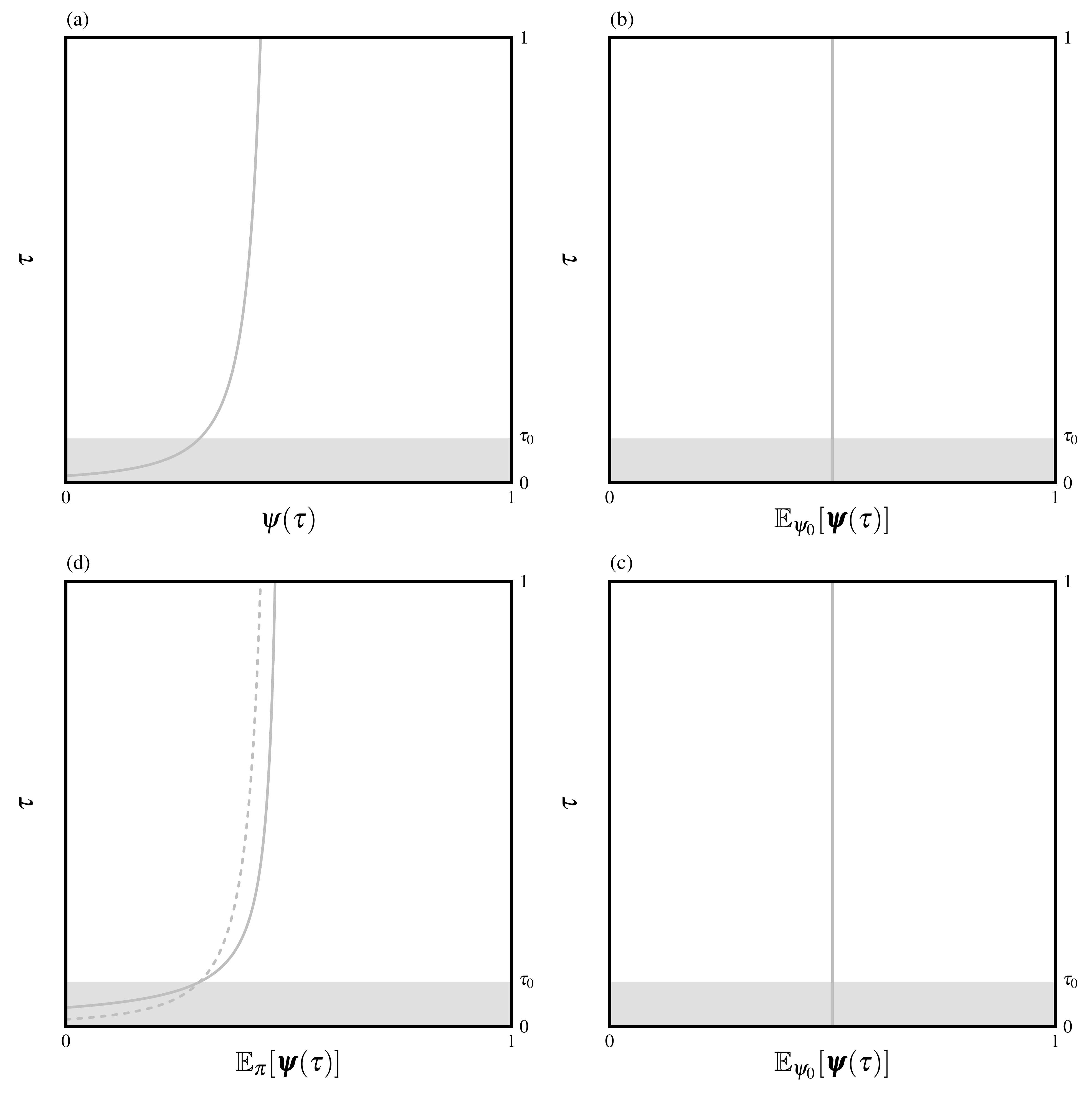}}
\caption{Share trajectories and their means at the fixed initial saturation \mbox{$\tau_0=0.1$}. Panels \mbox{(a--c)} use \mbox{$\pi(\psi)=1/4+\psi/2$} from Figure~\ref{fig:lft-two-body}. (a) Conditional trajectory from the prescribed share \mbox{$\psi_0=0.3$}. (b) Mean trajectory when the initial share is uniform on $[0,1]$ and integrated within the path measure. (c) Mean over uniformly sampled initial shares frozen within each system. Both means are $\mathbb E_{\psi_0}[\randombold{\psi}(\tau)]=1/2$, as in Eq.~(\ref{eq:rsb-uniform-initial-mean-trajectory}), although their overlap laws differ. (d) Mean trajectory (solid) from \mbox{$\psi_0=0.3$} under the full uniform quenched coefficient law of Figure~\ref{fig:quenched-affine}. The mean $\mathbb E_\pi[\randombold{\psi}(\tau)]$ is given by Eq.~(\ref{eq:rsb-quenched-mean-trajectory}). The dashed curve repeats the fixed-urn trajectory in (a). All averages use the full continuous laws. All boxes span $[0,1]^2$, with saturation increasing upwards. The shaded strip marks the prescribed initial sector $0\leq\tau<\tau_0$. Within the strip, the curves are formal backward continuations of the displayed expressions, not prescribed microscopic histories or averages of admissible continuations.}
\label{fig:rsb-trajectories}
\end{figure}
\clearpage
}

Figure~\ref{fig:rsb-trajectories}\label{call:fig:rsb-trajectories}
compares a prescribed trajectory with the three mean trajectories.
The annealed and quenched initial averages occupy the right column;
the prescribed share and the coefficient average occupy the left.
The same arrangement is used for the overlaps and free energies below.

\subsubsection{Quenched random urn functions}
\label{app:rsb-quenched}

Now choose the response once and keep it fixed throughout growth,
as in Section~\ref{sec:quenched}. Its sampling law is $\mathbb P_\pi$,
with expectation $\mathbb E_\pi$. Conditional on that choice, the
preceding construction applies unchanged. The same response is shared
by the configurations being compared. With an assigned initial share,
the overlap is concentrated for each urn. Averaging over urns gives
the following law at the prescribed initial share:
\begin{equation}
\mathbb E_\pi[P(q)]\,dq
=\mathbb E_\pi[\delta_{q_{\randombold{\pi}}(\psi_0,\psi_0)}(dq)].
\label{eq:rsb-random-fixed-initial}
\end{equation}
With a uniform frozen initial share, each system still has a single
conditional overlap. Average the law in Eq.~(\ref{eq:rsb-diagonal-law})
over the response to obtain its distribution between systems.
With a uniform share integrated within the path measure, each fixed
urn instead has the law in Eq.~(\ref{eq:rsb-mixture-law}); this law
can then be averaged over responses. These are the three initial
preparations of Section~\ref{app:rsb-initial}, followed by an external
average over the response. No parameter is resampled during a history.

For affine urns, take the uniform coefficient law of
Eq.~(\ref{eq:quenched-affine-domain}), independently of the initial share.
Both signs of $b$ contribute.
The densities at fixed coefficients follow by changing variables in
Eq.~(\ref{eq:rsb-affine-kernel}). For $b\ne0$, its diagonal is a
strictly convex quadratic polynomial in $\psi_0$. Differentiating
this diagonal polynomial varies both arguments together. Each root
in $(0,1)$ contributes its inverse absolute derivative. Thus, for a
uniformly sampled share frozen within each system, almost everywhere,
\begin{equation}
\mathbb E_{\psi_0}[P(q)]=
\sum\limits_{\substack{\psi_0\in(0,1):\\q_\pi(\psi_0,\psi_0)=q}}
\frac{1}{|\,d q_\pi(\psi_0,\psi_0)/d\psi_0\,|}.
\label{eq:rsb-general-diagonal-density}
\end{equation}
There are at most two roots; the density is zero where neither is
admissible. At an admissible minimum, a vanishing derivative produces
an integrable singularity.
For the share integrated within the path measure, first fix $\psi_0$ and change variables from the
independent uniform share $\psi_0'$ to the overlap. The dependence
on $\psi_0'$ is affine, so there is at most one admissible root.
Omitting the isolated value of $\psi_0$ where the slope vanishes,
we obtain, for almost every $q$,
\begin{equation}
P_{\mathrm{ann}}(q)=\int_0^1d\psi_0
\sum\limits_{\substack{\psi_0'\in(0,1):\\q_\pi(\psi_0,\psi_0')=q}}
\frac{1}{|\partial_{\psi_0'}q_\pi(\psi_0,\psi_0')|}.
\label{eq:rsb-general-mixed-density}
\end{equation}
Averaging the frozen-share density over the coefficient law, with
$\pi(\psi)=a+b\psi$ in the integrand, gives
\begin{equation}
\mathbb E_\pi\mathbb E_{\psi_0}[P(q)]=\int_0^1da\int_{-a}^{1-a}db
\sum\limits_{\substack{\psi_0\in(0,1):\\q_\pi(\psi_0,\psi_0)=q}}
\frac{1}{|\,d q_\pi(\psi_0,\psi_0)/d\psi_0\,|}.
\label{eq:rsb-average-diagonal-density}
\end{equation}
For the share integrated within each path measure, average
$P_{\mathrm{ann}}(q)$ with $\mathbb E_\pi$,
applying the same coefficient integral to
Eq.~(\ref{eq:rsb-general-mixed-density}). The two configurations
still share the same sampled response.
For a fixed initial share $\psi_0$, write
$q_\pi(\psi_0,\psi_0)$ as the quadratic polynomial in $a$ specified by
Eq.~(\ref{eq:rsb-affine-kernel}).
For almost every $q$, integrating over its admissible roots gives the
density of Eq.~(\ref{eq:rsb-random-fixed-initial}):
\begin{equation}
\mathbb E_\pi[P(q)]=\int_{-1}^1db
\sum_{\substack{a:\,q_\pi(\psi_0,\psi_0)=q\\
\max(0,-b)<a<\min(1,1-b)}}
\frac{1}{|\partial_a q_\pi(\psi_0,\psi_0)|}.
\label{eq:rsb-average-fixed-initial-density}
\end{equation}
The sum contains at most two roots. The inverse-derivative factors
can diverge where the derivative vanishes; their integral over $b$
need not diverge.
For the prescribed share $\psi_0=0.3$ in
Figure~\ref{fig:rsb-overlap-densities}(d), the overlap vanishes only at
$(a,b)=(1/2,0)$. Its quadratic growth in two parameter directions yields
a finite density at zero; varying a single initial share instead gives
the inverse--square--root singularity in
Figure~\ref{fig:rsb-overlap-densities}(c).
The line $b=0$ has zero probability under this coefficient law and
therefore introduces no atom into the averaged distributions.

The distinction between variation within an urn and variation between
urns is visible already in the first moment. At fixed coefficients
$a,b$, take $\randombold{\psi}_0$ and $\randombold{\psi}_0'$
independently and uniformly on $[0,1]$. The two expectations below
average these independent draws, and the variance refers to the same
product law. For $b<1$, the overlap has mean and variance
\begin{equation}
\mathbb E_{\psi_0}\mathbb E_{\psi_0'}[q_\pi(\randombold{\psi}_0,\randombold{\psi}_0')]
=\frac{(2\psi_*-1)^2}{1-\tau_0}
\int_{\tau_0}^1[1-b(\tau_0/\tau)^{1-b}]^2d\tau.
\label{eq:rsb-mixed-moments}
\end{equation}
\begin{multline}
\operatorname{Var}(q_\pi(\randombold{\psi}_0,\randombold{\psi}_0'))=
\frac{2b^2(2\psi_*-1)^2}{3(1-\tau_0)^2}
\left\{\int_{\tau_0}^1(\tau_0/\tau)^{1-b}
[1-b(\tau_0/\tau)^{1-b}]d\tau\right\}^2\\
+\frac{b^4}{9(1-\tau_0)^2}
\left\{\int_{\tau_0}^1(\tau_0/\tau)^{2(1-b)}d\tau\right\}^2.
\label{eq:rsb-mixed-variance}
\end{multline}
For the uniform coefficient law, the mean overlap is strictly positive.
By contrast, the mean response is $1/2$ and would give zero future
overlap. Thus averaging the response does not preserve the overlap law.
At $\tau_0=0.1$, an assigned share $\psi_0=1/2$ gives an averaged overlap
mean $0.185316$ and variance $0.047498$. Freezing a uniform share and
then averaging gives mean $0.197249$ and variance $0.051115$; integrating
the uniform share within each path measure gives mean $0.185316$ and
variance $0.050491$. These values follow by integrating the explicit
densities above. They illustrate the three cases without changing
the coefficient distribution.

For a prescribed initial share, the coefficient average determines the
mean trajectory. At fixed $b$, the conditional mean of $\randombold{a}/(1-b)$ is
$1/2$, and the marginal density of $\randombold{b}$ is $1-|b|$. Averaging
Eq.~(\ref{eq:quenched-affine-trajectory}) over the full coefficient law gives
\begin{equation}
\mathbb E_\pi[\randombold{\psi}(\tau)]
=1/2+(\psi_0-1/2)
\frac{(1-\tau_0/\tau)^2}{[\log(\tau/\tau_0)]^2},
\qquad \tau_0<\tau\leq1.
\label{eq:rsb-quenched-mean-trajectory}
\end{equation}
The expression extends continuously to $\tau=\tau_0$, where its value is
$\psi_0$. For $\psi_0\ne1/2$, each fixed--coefficient trajectory with $b<1$ approaches
its own fixed share as a power of $\tau_0/\tau$, whereas the mean
deviation from $1/2$ decays as $1/[\log(\tau/\tau_0)]^2$ for large
$\tau/\tau_0$. The slower relaxation comes from slopes $b$ arbitrarily
close to one, for which attraction to the fixed share is arbitrarily weak.
The overlap law averages the conditional laws, as above;
substituting the mean trajectory into the overlap integral would give
a different quantity.

\afterpage{%
\clearpage
\begin{figure}[p]
\centering
\makebox[\linewidth][l]{\hspace*{-8.894462pt}\includegraphics[width=352.185673pt]{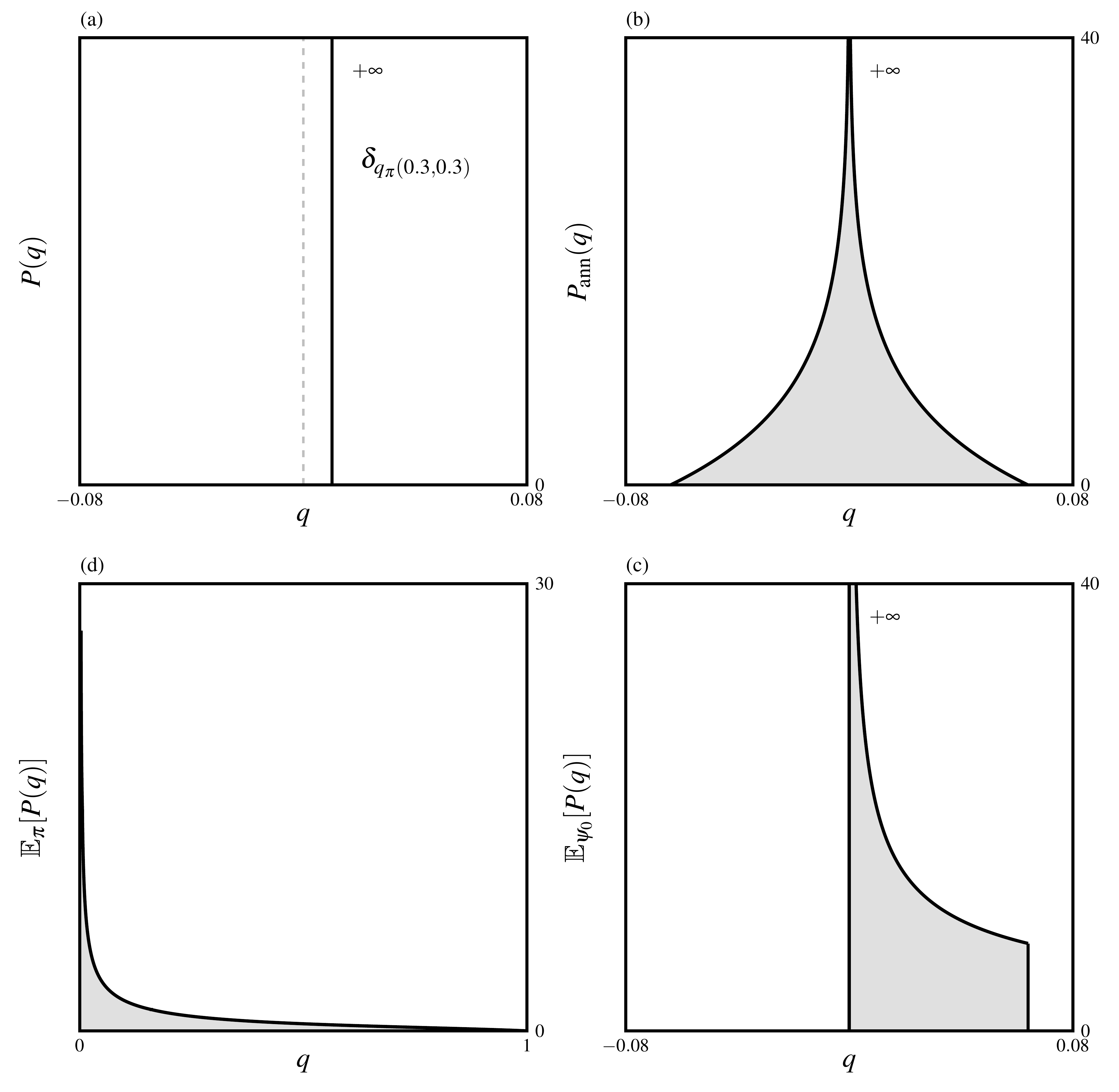}}
\caption{Limiting future overlap laws for the same four preparations. Panels \mbox{(a--c)} use \mbox{$\pi(\psi)=1/4+\psi/2$} from Figure~\ref{fig:lft-two-body}, with \mbox{$\tau_0=0.1$}. (a) The conditional law $P(q)$ at \mbox{$\psi_0=0.3$} has unit mass at $q_\pi(0.3,0.3)\simeq0.010234$, represented symbolically by the vertical line and the $+\infty$ label. (b) The uniform initial share integrated within the path measure gives $P_{\mathrm{ann}}(q)$ in Eq.~(\ref{eq:rsb-symmetric-mixed-density}). (c) The external average over uniform frozen shares gives $\mathbb E_{\psi_0}[P(q)]$ in Eq.~(\ref{eq:rsb-symmetric-diagonal-density}); the conditional law in each system remains concentrated. Both use $q_\pi(1,1)\simeq0.063961$. In (c), the density is zero for $q<0$ and jumps to zero at the upper support boundary, as marked by the vertical segment. (d) The mean $\mathbb E_\pi[P(q)]$ at \mbox{$\psi_0=0.3$}, under the full uniform quenched coefficient law in Figure~\ref{fig:quenched-affine}, follows Eq.~(\ref{eq:rsb-average-fixed-initial-density}). The densities in (b,c) diverge integrably at zero and carry no atom; their branches continue beyond the displayed vertical range. The solid line at zero in (c) marks this divergence. The density in (d) has a finite limit of approximately $21.13$ at zero. In all panels, the overlap compares two histories by multiplying the spins added at the same step after the initial saturation $\tau_0=0.1$ and averaging these products over those steps, as in Eq.~(\ref{eq:rsb-overlap-definition}). In particular, the common mean trajectory in (b,c) does not determine their different overlap laws.}
\label{fig:rsb-overlap-densities}
\end{figure}
\clearpage
}

Figure~\ref{fig:rsb-overlap-densities}\label{call:fig:rsb-overlap-densities}
shows the overlap laws for these four preparations, using the prescribed
share $\psi_0=0.3$ in the two left panels. Equal mean trajectories in
the right panels of Figure~\ref{fig:rsb-trajectories} coexist with
distinct overlap distributions because the two histories choose their
initial shares independently or retain the same frozen share.

\subsubsection{Free energy and the order of averaging}
\label{sec:rsb-free-energy}

The free energy of Section~\ref{sec:entropy} weights histories by
their final count, including the prescribed initial population.
We retain the CGF $\zeta(\beta)$ of Eq.~(\ref{eq:cgf}). The expectation
$\mathbb E$ is over continuations from the prescribed initial history,
with the urn function and initial share fixed. Their dependence
remains implicit, as in Section~\ref{sec:quenched}; integrals and
evaluations in $\pi$ and $\psi_0$ below act on these conditional
quantities. The operators $\mathbb E_\pi$ and $\mathbb E_{\psi_0}$
average the sampled inputs, while $\randombold{\zeta}$ denotes the
conditional CGF before those inputs are fixed. Thus
\begin{equation}
\zeta(\beta)=\lim_{T\to\infty}\frac{1}{T}
\log\mathbb E[\exp(\beta T\psi_T)].
\label{eq:rsb-conditional-cgf}
\end{equation}
The conditional law is normalized, so $\zeta(0)=0$. Here $\beta$
weights the final count; weighting an overlap would define a different
generating function. This CGF requires the probabilities of complete
histories, which are not determined by the limiting overlap or by
factorization of a fixed number of increments.

For the affine examples, add the terminal term $\beta\varphi(1)$
to the action difference in Eq.~(\ref{eq:VARPRINZ}), restrict its
integrals to $[\tau_0,1]$, and prescribe
$\varphi(\tau_0)=\tau_0\psi_0$. The equation for a fixed finite seed,
Eq.~(\ref{eq:cgf-equation}), does not impose this macroscopic initial
population. First take $b\ne0$, $\beta\ne0$, $0<a<1$ and
$0<a+b<1$. Along an optimizing path, write $\beta(\tau)$ for the
conjugate field, whose terminal value is the argument $\beta$ of
the CGF. Its initial value is denoted by $\beta_0$:
\begin{equation}
\beta(1)=\beta,\qquad \beta(\tau_0)=\beta_0.
\label{eq:rsb-characteristic-boundaries}
\end{equation}
The equations of the variational problem are
\begin{equation}
\partial_\tau\varphi(\tau)
=\frac{\pi(\varphi(\tau)/\tau)\exp(\beta(\tau))}
{1+\pi(\varphi(\tau)/\tau)[\exp(\beta(\tau))-1]}.
\label{eq:rsb-characteristic-path}
\end{equation}
\begin{equation}
\partial_\tau\beta(\tau)
=-\frac{b[\exp(\beta(\tau))-1]}
{\tau\{1+\pi(\varphi(\tau)/\tau)[\exp(\beta(\tau))-1]\}}.
\label{eq:rsb-characteristic-equations}
\end{equation}
The initial and terminal fields have the same sign. Moreover,
$|\beta_0|>|\beta|$ for $b>0$ and $|\beta_0|<|\beta|$ for $b<0$.
Integrating these equations determines $\beta_0$ through a single
integral, with $\beta'$ only the variable of integration:
\begin{multline}
\frac{1-\tau_0}
{\tau_0[1+(a+b\psi_0)(\exp(\beta_0)-1)]}
=\\
-\frac{1}{b}\int_{\beta_0}^{\beta}
\frac{\exp[(1-a)(\beta'-\beta_0)/b]}{\exp(\beta')-1}
\left|\frac{\exp(\beta_0)-1}{\exp(\beta')-1}\right|^{1/b}
\,d\beta'.
\label{eq:rsb-affine-field}
\end{multline}
Integration by parts in the variational expression then gives
\begin{multline}
\zeta(\beta)=\log\tau_0+\tau_0\psi_0\beta_0
+\frac{(1-a)(\beta-\beta_0)}{b}\\
+(1-\tau_0)\log[1+(a+b\psi_0)(\exp(\beta_0)-1)]
-\frac{1}{b}\log\left|
\frac{\exp(\beta)-1}{\exp(\beta_0)-1}\right|.
\label{eq:rsb-affine-macro-cgf}
\end{multline}
Joint convexity of Bernoulli relative entropy in its two arguments
makes this path a global optimizer. At nonzero field the optimizing
path, and hence $\beta_0$, is unique. At $\beta=0$, $\zeta(0)=0$.
For $b=0$, the independent future increments give directly
\begin{equation}
\zeta(\beta)=\beta\tau_0\psi_0
+(1-\tau_0)\log[1-a+a\exp(\beta)].
\label{eq:rsb-constant-cgf}
\end{equation}
An exactly all-zero prefix is absorbing when $a=0$, giving
$\zeta=0$; an exactly all-one prefix is absorbing when $a+b=1$,
giving $\zeta=\beta$. For $0<\psi_0<1$, the quadrature extends to
nonconstant affine responses on the parameter boundary. The accessible
shares satisfy
\begin{equation}
\tau_0\psi_0\leq\varphi(\tau)/\tau
\leq1-\tau_0(1-\psi_0),\qquad \tau_0\leq\tau\leq1.
\label{eq:rsb-accessible-shares}
\end{equation}
On this interval the transition probabilities stay away from zero and one. These parameter boundaries
have zero weight under the continuous coefficient law below.

The conditional CGF gives the first of the three preparations at a
fixed urn function. If the initial share is uniform and then frozen,
the logarithm is taken before its average:
\begin{equation}
\mathbb E_{\psi_0}[\randombold{\zeta}(\beta)]
=\lim_{T\to\infty}\int_0^1\frac{\log \mathbb E[\exp(\beta T\psi_T)]}{T}\,d\psi_0
=\int_0^1\zeta(\beta)\,d\psi_0.
\label{eq:rsb-fixed-response-quenched-seed}
\end{equation}
If the uniform share is instead integrated within the path measure,
the average precedes the logarithm. Denote this CGF by
$\zeta_{\mathrm{ann}}(\beta)$, with the subscript again referring only
to the initial-share mixture. The selected share is one for
$\beta>0$ and zero for $\beta<0$, so
\begin{equation}
\zeta_{\mathrm{ann}}(\beta)=\lim_{T\to\infty}\frac{1}{T}\log\int_0^1 \mathbb E[\exp(\beta T\psi_T)]\,d\psi_0
=\left.\zeta(\beta)\right|_{\psi_0=\mathbb I(\beta>0)}.
\label{eq:rsb-fixed-response-annealed-seed}
\end{equation}
At zero field, $\zeta_{\mathrm{ann}}(0)=0$. The same limit is obtained
by giving equal weight to the admissible initial counts $j=0,\ldots,|S_0|$. Their average partition lies between
its maximum divided by $|S_0|+1$ and that maximum. Coupling the integer
counts preserves their order, so the maximizing count is $|S_0|$ for
$\beta>0$ and zero for $\beta<0$. The logarithmic discrepancy divided
by $T$ is at most $\log(|S_0|+1)/T$. The initial share is averaged only
at initialization; it is not resampled during growth.

Next freeze the random coefficients $\randombold{a},\randombold{b}$
with the uniform law of Eq.~(\ref{eq:quenched-affine-domain}). In the
integrals, $\pi(\psi)=a+b\psi$. As in
Eq.~(\ref{eq:quenched-affine-uniform-mean}), boldface marks the CGF
before its sampled inputs are fixed. A prescribed share gives
\begin{equation}
\mathbb E_\pi[\randombold{\zeta}(\beta)]
=\int_0^1 da\int_{-a}^{1-a}db\,\zeta(\beta).
\label{eq:rsb-random-response-fixed-seed}
\end{equation}
If the initial share is also uniform and frozen, the average over
both sources of variation is
\begin{equation}
\mathbb E_\pi\mathbb E_{\psi_0}[\randombold{\zeta}(\beta)]
=\int_0^1 da\int_{-a}^{1-a}db\int_0^1d\psi_0\,\zeta(\beta).
\label{eq:rsb-random-response-quenched-seed}
\end{equation}
If the uniform initial share is instead inside each path measure,
only the coefficients remain quenched. The resulting
$\randombold{\zeta}_{\mathrm{ann}}$ is random only through these
coefficients:
\begin{equation}
\mathbb E_\pi[\randombold{\zeta}_{\mathrm{ann}}(\beta)]
=\int_0^1da\int_{-a}^{1-a}db\,\zeta_{\mathrm{ann}}(\beta)
=\left.\mathbb E_\pi[\randombold{\zeta}(\beta)]
\right|_{\psi_0=\mathbb I(\beta>0)}.
\label{eq:rsb-random-response-annealed-seed}
\end{equation}
These are averages of the solved conditional CGFs, rather than the
CGF evaluated at mean coefficients. The bound between
$\min(0,\beta)$ and $\max(0,\beta)$ permits dominated convergence.
Jensen's inequality puts the externally averaged free energy below
the one with the initial share integrated inside the measure, both
before and after the coefficient average. At nonzero field the latter
weights each initial share by the conditional exponential moment
of its final count; the frozen-share average retains the external
uniform law.

For comparison with the future overlap, the initial count is removed
before taking any internal initial-share average. Since that count
divided by $T$ tends to $\tau_0\psi_0$, its removal gives
\begin{equation}
\lim_{T\to\infty}\frac{1}{T}\log\mathbb E\!\left[
\exp\!\left(\beta\sum\limits_{|S_0|<n\leq T}\sigma_n\right)\right]
=\zeta(\beta)-\beta\tau_0\psi_0.
\label{eq:rsb-future-cgf}
\end{equation}
For affine urns, future additions increase with the initial count
when $b>0$ and decrease when $b<0$. Take the uniform initial-share
average of the exponential moment in Eq.~(\ref{eq:rsb-future-cgf})
before the logarithm and the limit. Its limiting free energy is
\begin{equation}
\left.[\zeta(\beta)-\beta\tau_0\psi_0]
\right|_{\psi_0=\mathbb I(\beta b>0)}.
\label{eq:rsb-future-annealed-selector}
\end{equation}
For $b=0$, the future partition is independent of the seed. A frozen
uniform share instead gives
\begin{equation}
\mathbb E_{\psi_0}[\randombold{\zeta}(\beta)-\beta\tau_0\randombold{\psi}_0]
=\int_0^1[\zeta(\beta)-\beta\tau_0\psi_0]d\psi_0
=\mathbb E_{\psi_0}[\randombold{\zeta}(\beta)]-\beta\tau_0/2.
\label{eq:rsb-future-quenched-seed}
\end{equation}
The corresponding three energies with quenched coefficients are
obtained by applying $\mathbb E_\pi$, given by the coefficient integral
in Eq.~(\ref{eq:rsb-random-response-fixed-seed}), respectively to
Eqs.~(\ref{eq:rsb-future-cgf}),
(\ref{eq:rsb-future-quenched-seed}) and
(\ref{eq:rsb-future-annealed-selector}). In the last case the selector
depends on $b$ and must remain inside that coefficient average. For the future
spin sum, replace $\beta$ by $2\beta$ in these expressions and subtract
$\beta(1-\tau_0)$.

Take the limit $T\to\infty$ before differentiating with respect
to $\beta$. The derivative at zero of the conditional CGF is the
limiting terminal share. At fixed coefficients and initial share,
\begin{equation}
\partial_\beta\zeta(0)
=\frac{a(1-\tau_0^{1-b})}{1-b}+\psi_0\tau_0^{1-b},
\qquad b<1.
\label{eq:rsb-conditional-slope}
\end{equation}
Its value at the P\'olya vertex is $\psi_0$. The frozen uniform-share
average uses $\psi_0=1/2$; the left and right derivatives after
integration within the path measure use $\psi_0=0$ and $\psi_0=1$.
Under the uniform coefficient law, the marginal density of
$\randombold{b}$ is $1-|b|$, giving
\begin{equation}
\mathbb E_\pi[\tau_0^{1-\randombold{b}}]
=\frac{(1-\tau_0)^2}{[\log(1/\tau_0)]^2}.
\label{eq:rsb-disorder-R}
\end{equation}
The conditional mean of $\randombold{a}/(1-b)$ at fixed $b$ is $1/2$.
Consequently, at a prescribed initial share,
\begin{equation}
\left.\partial_\beta\mathbb E_\pi[\randombold{\zeta}(\beta)]
\right|_{\beta=0}
=\frac12+\left(\psi_0-\frac12\right)
\frac{(1-\tau_0)^2}{[\log(1/\tau_0)]^2}.
\label{eq:rsb-averaged-conditional-slope}
\end{equation}
A further external uniform-share average gives derivative $1/2$.
When that share is instead integrated within the path measure, the
two one-sided derivatives are
\begin{equation}
\lim_{\beta\to0^\pm}\frac{1}{\beta}
\mathbb E_\pi[\randombold{\zeta}_{\mathrm{ann}}(\beta)]
=\frac12\pm\frac{(1-\tau_0)^2}{2[\log(1/\tau_0)]^2}.
\label{eq:rsb-total-annealed-slopes}
\end{equation}
Removing the initial population changes the coefficient of $\psi_0$
in the derivative from $\tau_0^{1-b}$ to
$\tau_0^{1-b}-\tau_0$. The difference between the extreme initial
conditions is governed by its absolute value, whose mean is
\begin{equation}
\mathbb E_\pi[|\tau_0^{1-\randombold{b}}-\tau_0|]
=\frac{1-\tau_0^2-2\tau_0\log(1/\tau_0)}
{[\log(1/\tau_0)]^2}.
\label{eq:rsb-disorder-S}
\end{equation}
The derivative after external averaging of the share and coefficients
is $(1-\tau_0)/2$. With the share integrated within each path measure,
the future-count derivatives are instead
\begin{multline}
\lim_{\beta\to0^\pm}\frac{1}{\beta}
\int_0^1da\int_{-a}^{1-a}db\,
\left.[\zeta(\beta)-\beta\tau_0\psi_0]
\right|_{\psi_0=\mathbb I(\beta b>0)}\\
=\frac{1-\tau_0}{2}\pm
\frac{1-\tau_0^2-2\tau_0\log(1/\tau_0)}{2[\log(1/\tau_0)]^2}.
\label{eq:rsb-future-annealed-slopes}
\end{multline}
At a fixed share, subtract $\tau_0\psi_0$ from the total-count
derivative. Bounded difference quotients justify these coefficient
averages.
A broad overlap law after external averaging can therefore coexist
with a differentiable averaged CGF. Integrating the initial share
within the measure produces the displayed cusps by selecting among
its initial sectors.

\afterpage{%
\clearpage
\begin{figure}[p]
\centering
\makebox[\linewidth][l]{\hspace*{-9.244983pt}\includegraphics[width=354.201167pt]{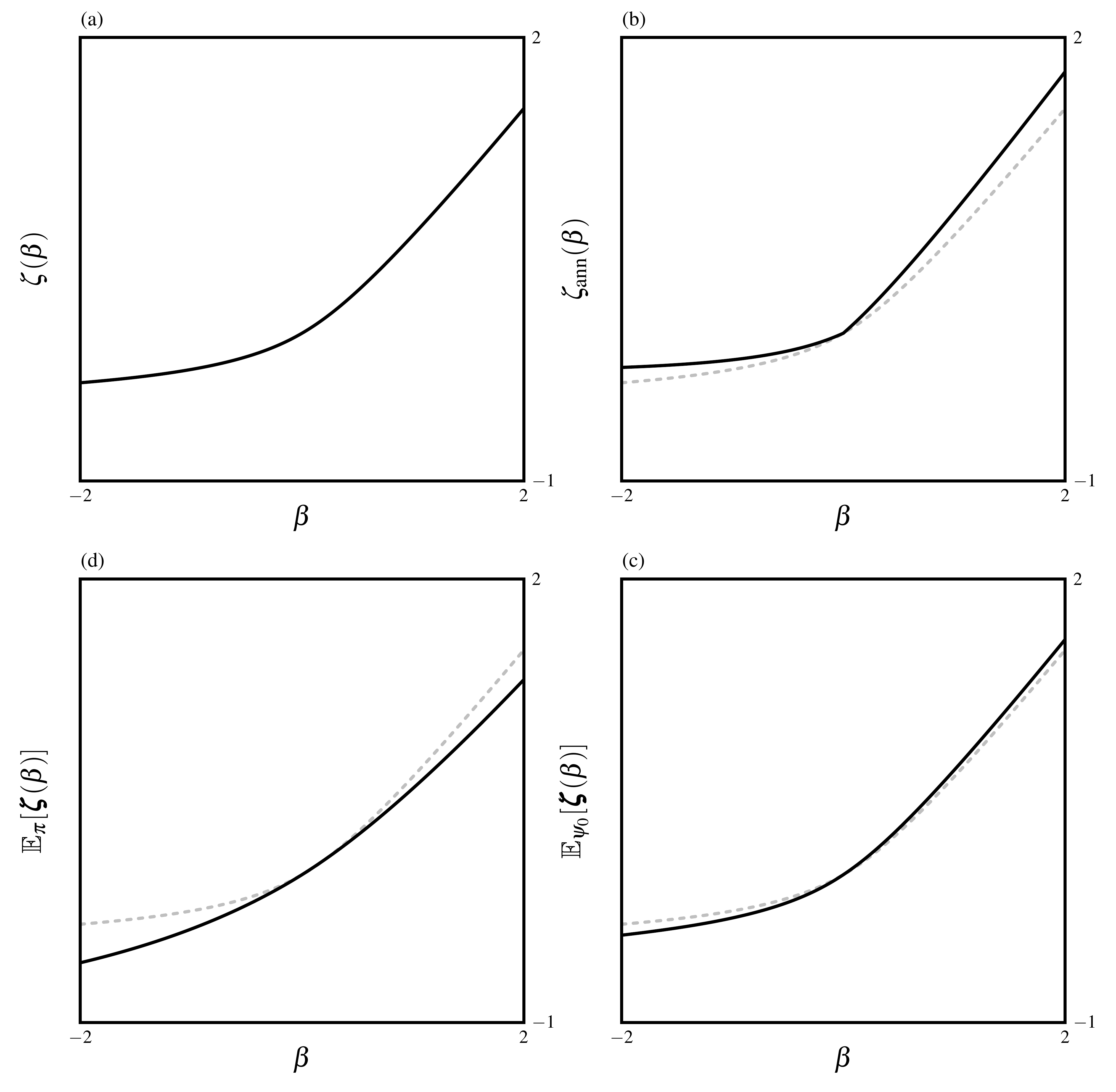}}
\caption{Total-count CGFs at \mbox{$\tau_0=0.1$} for the same four preparations. Panels \mbox{(a--c)} use the fixed response \mbox{$\pi(\psi)=1/4+\psi/2$} of Figure~\ref{fig:lft-two-body}. (a) Conditional CGF $\zeta(\beta)$ at \mbox{$\psi_0=0.3$}. (b) Annealed CGF $\zeta_{\mathrm{ann}}(\beta)$, with the uniform initial share integrated before the logarithm. (c) Quenched mean $\mathbb E_{\psi_0}[\randombold{\zeta}(\beta)]$, with the uniform initial share averaged after the logarithm. (d) Mean $\mathbb E_\pi[\randombold{\zeta}(\beta)]$ at \mbox{$\psi_0=0.3$}, under the full uniform coefficient law of Figure~\ref{fig:quenched-affine}. Light-gray dashed curves in \mbox{(b--d)} repeat the conditional CGF in (a). All CGFs include the initial population and are normalized by $T$. The conditional CGFs follow Eqs.~(\ref{eq:rsb-affine-field}) and (\ref{eq:rsb-affine-macro-cgf}); their averages follow Eqs.~\mbox{(\ref{eq:rsb-fixed-response-quenched-seed})--(\ref{eq:rsb-random-response-fixed-seed})}. Only (b) has a cusp at zero, with left and right derivatives approximately $0.341886$ and $0.658114$. The other curves are differentiable there; their derivatives equal the corresponding terminal mean shares, approximately $0.436754$ in (a), $1/2$ in (c), and $0.469445$ in (d).}
\label{fig:rsb-free-energy-four}
\end{figure}
\clearpage
}

Figure~\ref{fig:rsb-free-energy-four}\label{call:fig:rsb-free-energy-four}
compares the four total-count CGFs with the same prescribed share
$\psi_0=0.3$ in the left panels. The annealed initial mixture has a
cusp at zero, whereas the quenched averages remain differentiable.
Their derivatives at zero recover the corresponding terminal mean shares
in Figure~\ref{fig:rsb-trajectories}.

\subsubsection[Trajectory states and Kernel/Blanket]{Trajectory states and the Kernel/Blanket construction}\label{app:rsb-kernel}

The trajectory, overlap and free-energy calculations above provide an
explicit setting for the Kernel/Blanket description of spin
measures~\cite{KERNEL=000020THEO,Franchini2024}. A spin kernel records
configurations together with their probabilities. In the Blanket
construction, an ordered partition of the sites separates interactions
within each part from interfaces with preceding parts~\cite{RSBSPA1}. For the urn,
temporal order supplies a sequence of conditional laws, while the growth
dynamics determines the local magnetizations. We can therefore relate
the components and their overlaps to the weights entering the free-energy
calculation.

Fix the urn function and the initial preparation. In the microscopic
kernel, rows represent times and columns represent histories: each
column contains the spins $2\sigma_n-1$ of its history and has width
equal to that history's probability~\cite{KERNEL=000020THEO}. For the
future overlaps considered here, retain the rows after $S_0$, with equal
weights summing to one. The sectors $S_\ell$ group these rows. Histories
with the same spin counts in every sector can in turn be grouped,
retaining each group's total probability and its conditional history
law. This is the disintegration by sector magnetizations used in the
layer description~\cite{RSBSPA2,RSBSPA3v2}. For HLS, the conditional law
within such a group is generally nonuniform: the order of additions
changes intermediate shares and hence the transition probabilities.

The limiting components are identified by the concentration results of
Section~\ref{app:rsb-trajectories}. With $\pi$ and $\psi_0$ prescribed,
the local spin means converge to $2\pi(\psi(\tau))-1$, whose sector
averages are given by Eq.~(\ref{eq:rsb-sector-magnetizations}).
Factorization of any fixed number of increments and concentration of
the conditional overlap support the interpretation as a pure-state
component specified by local magnetizations~\cite{KERNEL=000020THEO}.
These properties come from the urn dynamics; grouping histories by
sector counts alone does not establish them. Initial shares producing
identical profiles describe the same future magnetizations.

The kernel has two complementary interpretations. For two independent
continuations, integrating the product of their limiting profiles over
future times gives the overlap in Eq.~(\ref{eq:rsb-overlap-limit});
Eq.~(\ref{eq:rsb-sector-overlap}) recovers it by refining the sectors.
Conversely, averaging the product of two distinct rows over columns
gives the correlation between the corresponding spins~\cite{KERNEL=000020THEO}.
Within a prescribed component this correlation factorizes in the limit;
a mixture can retain correlations through the common initial share of
each history. Two independent histories from the mixed measure instead
sample their shares independently, whereas the quenched mean compares
histories with the same frozen share. A quenched urn function is likewise
shared by both histories of a given system. This explains how the
identical mean trajectories in Figure~\ref{fig:rsb-trajectories}(b,c)
coexist with the different overlap laws in
Figure~\ref{fig:rsb-overlap-densities}(b,c).

The same history measure also determines the CGF. Here $\beta$ weights
the total black count, including the initial population, and $\mathbb E$
is the original HLS expectation at prescribed $\pi$ and $\psi_0$.
The dependence on these inputs remains implicit, as in
Section~\ref{sec:rsb-free-energy}. At finite capacity, multiply each
history's probability by the exponential of $\beta$ times its final
count and normalize the resulting measure. For the continuous uniform
choice of the initial-share parameter, with the initial count rounded
as in Section~\ref{app:rsb-trajectories}, its reweighted distribution is
\begin{equation}
\frac{\mathbb E[\exp(\beta T\psi_T)]\,d\psi_0}
{\displaystyle\int_0^1\mathbb E[\exp(\beta T\psi_T)]\,d\psi_0}.
\label{eq:rsb-tilted-initial-weights}
\end{equation}
Thus each component is weighted by its conditional partition function.
For uniformly sampled admissible initial counts, the same reweighting
applies to their discrete masses.

To compute the conditional partition function entering these weights,
keep $T$, the microscopic initial prefix and $\pi$ fixed, and add one
increment at a time after $S_0$. Conditioning on the past and summing
over the two values of
$\sigma_{n+1}$ gives, for $|S_0|\leq n<T$, the exact identity
\begin{equation}
\mathbb E[\exp(\beta(n+1)\psi_{n+1})]
=\mathbb E\!\left[\exp(\beta n\psi_n)
\{1+\pi(\psi_n)[\exp(\beta)-1]\}\right].
\label{eq:rsb-cavity-increment}
\end{equation}
The factor in braces is the exponential moment of the added ball,
conditional on the current share. Dividing by the preceding moment and
taking the logarithm gives the increment of the log-partition function.
The ratio averages the insertion factor over the preceding histories
with the exponential weighting of their current counts; their
fluctuations remain part of the calculation.
Summing these increments cancels the intermediate terms and, with the
initial-count contribution, recovers the logarithm of the terminal moment.
Division by $T$ and the capacity limit give $\zeta$. Adding one ball is
the urn's insertion of a new site. This provides a direct comparison
with the Constructive Cavity Method (CCM), which evaluates free-energy
increments using the states and weights of the preceding
system~\cite{RSBSPA1,RSBSPA3v2}.

At zero field the insertion factor is one. Multiplying a transition
probability by the exponential of $\beta$ times the added increment
and normalizing the two outcomes divides by the insertion factor,
which generally depends on the history. To recover the globally
reweighted history law, its product over the additions must be restored
and normalized over the full measure. This is the distinction between
local Blanket normalization and the complete Gibbs
measure~\cite{Franchini2024,RSBSPA3v2}.
For the terminal weighting of the original HLS law, the conditioning
also matters. Within a component specified only by its initial share,
the final count generally varies, so the field can change both the
component's weight and its internal history law.
If the initial prefix and every sector count are fixed instead, the
final count is constant within the group, so only the weights between
groups change. In the quenched averages the external laws of initial
shares and urn functions remain prescribed. Figure~\ref{fig:rsb-free-energy-four}
shows the resulting differences between the CGFs.

This interpretation of trajectory components also applies to nonlinear
responses under the stated concentration assumptions; the affine
examples make the profiles explicit. At a fixed affine response with
$b<1$, the profiles in
Eq.~(\ref{eq:rsb-affine-local-mean}) lie on an affine line, with distances given by
Eq.~(\ref{eq:rsb-profile-distance}). The Kernel/Blanket framework also
describes the hierarchical geometry of the Sherrington--Kirkpatrick
construction. There, ultrametricity follows from shared ancestors and
the renewal of fluctuations across branches~\cite{RSBSPA3v2}. Random energy
model limits~\cite{RSBREMLinear} and Poisson--Dirichlet
weights~\cite{RSBREMWeights} depend on the corresponding energy statistics.
For HLS, the preparation and growth law supply the weights directly.
The urn thus makes explicit how the same history measure determines
trajectory components, their overlaps and the increments entering its
free energy.

\par\endgroup
\section*{Acknowledgments}

We thank Alexandros Gelastopoulos (IAS Toulouse) and Pantelis Analytis
(University of Southern Denmark) for stimulating discussions, and
for bringing to my attention the reference \cite{BB}. I am especially
grateful to A. G. for providing Figure~\ref{fig:Vander}.
\section*{Declaration of Generative AI}
Research, writing, and proofreading of the present manuscript were performed with AI assistance (GPT-6 Astra) under strict human guidance and scrutiny. We independently checked all mathematical arguments and take full responsibility for the content of the publication.

\end{document}